%% file: main.tex
\RequirePackage[orthodox]{nag} 
\documentclass[english, a4paper, hidelinks]{scrartcl} 
\input{packages}
\input{commands}
\usepackage[backend=biber,style=numeric,giveninits=true,hyperref=true,isbn=false,url=true,doi=true,maxnames=10]{biblatex}
\numberwithin{equation}{section} 

\title{Steady State Blended Gas Flow on Networks: Existence and Uniqueness of Solutions}
\author{Alena Ulke\footnotemark[1], \;
	Michael Schuster\footnotemark[2], \; 
	Simone Göttlich\footnotemark[1]
}
\date{\today}
\begin{document}
\maketitle{}
\footnotetext[1]{University of Mannheim, Department of Mathematics, 68131 Mannheim, Germany (ulke@uni-mannheim.de, goettlich@uni-mannheim.de)}
\footnotetext[2]{Friedrich-Alexander University Erlangen-Nürnberg, Department of Mathematics, 91058 Erlangen, Germany (michi.schuster@fau.de)}
%
%
%
%
%
%
%
\input{tex/abstract}
\input{tex/introduction}
\input{tex/gas_flow_model}

\input{tex/existence_network_tree}

\input{tex/existence_network_with_cycle}
\input{tex/experiments}
%
%
\input{tex/conclusion}
\printbibliography{}
\end{document}

%% file: packages.tex
\usepackage[utf8]{inputenc} 
\usepackage[T1]{fontenc}  
\usepackage{lmodern} 
\usepackage{pdfpages}
\usepackage{scrhack}
\usepackage[main = english, ngerman]{babel} 
\usepackage{csquotes}
\usepackage{amsmath,amssymb,amsthm,bm}  
\usepackage{mathtools}
\usepackage{blkarray}
\usepackage{tabularx}
\usepackage{booktabs} 
\usepackage{enumitem}
\usepackage{siunitx}  
\usepackage{setspace}
\usepackage[linesnumbered, ruled, algosection]{algorithm2e}
\IncMargin{1em}
\DontPrintSemicolon 
\usepackage{graphicx}
\usepackage{wrapfig}
\usepackage{subcaption}
\usepackage{hyperref}  

\usepackage{bookmark}
\usepackage[noabbrev, capitalise]{cleveref}

\usepackage{xcolor}
\usepackage{grffile}
\usepackage{float} 
\usepackage{tikz, pgf}
\usetikzlibrary{fit,decorations.pathreplacing,patterns,shapes}
\usetikzlibrary{plotmarks}
\usetikzlibrary{arrows.meta, automata}
\usetikzlibrary{calc,intersections}
\usetikzlibrary{external}
\usepackage{multirow}
\usepackage{multicol} 
\usepackage{pgfplots}
\pgfplotsset{compat=newest}
\usepgfplotslibrary{patchplots}

\usepackage{url}
\usepackage{verbatim}
\usepackage[usestackEOL]{stackengine}

\usepackage{a4wide}

%% file: commands.tex
\theoremstyle{plain}
\newtheorem{theorem}{Theorem}[section] 
\newtheorem*{theorem*}{Theorem}
\newtheorem{lemma}[theorem]{Lemma} 
\newtheorem*{lemma*}{Lemma} 
\newtheorem{korollar}[theorem]{Corollary} 
\newtheorem*{korollar*}{Corollary} 
 
\newtheorem*{satz*}{Satz} 
\theoremstyle{definition}
\newtheorem{definition}[theorem]{Definition}
\newtheorem*{definition*}{Definition} 

\newtheorem*{beispiel*}{Beispiel} 
\theoremstyle{remark} 
\newtheorem{remark}[theorem]{Remark}
\newtheorem*{remark*}{Remark}

\newtheorem*{anmerkung*}{Anmerkung}

\newcommand{\R}{\mathbb{R}}
\newcommand{\N}{\mathbb{N}}

\renewcommand{\P}{\mathcal{P}}
\newcommand{\C}{\mathcal{C}}

\newcommand{\G}{\mathcal{G}}

\newcommand{\V}{\mathcal{V}}
\newcommand{\E}{\mathcal{E}}
\newcommand{\x}{\bm{x}}

\newcommand{\q}{\bm{q}}
\newcommand{\p}{\bm{p}}
\newcommand{\bb}{\bm{b}}
\newcommand{\etaa}{\bm{\eta}}

\newcommand{\sigmaa}{\bm{\sigma}}

\newcommand{\cl}{{\mathrm{cl}}}
\newcommand{\crr}{{\mathrm{cr}}}

\newcommand{\inn}{{\mathrm{in}}}
\newcommand{\out}{{\mathrm{out}}}
\newcommand{\const}{{\mathrm{const}}}

\newcommand{\flow}{{\mathrm{flow}}}

\newcommand{\sign}{{\mathrm{sign}}}

\newcommand{\dt}{\frac{\partial}{\partial t}}
\newcommand{\dx}{\frac{\partial}{\partial x}}

\newcommand{\red}[1]{\textcolor{red}{#1}}
\newcommand{\rhong}{\rho_{\mathrm{NG}}}
\newcommand{\rhohh}{\rho_{\mathrm{H}_2}}
\newcommand{\qng}{q_{\mathrm{NG}}}
\newcommand{\qhh}{q_{\mathrm{H}_2}}

\newcommand{\cng}{\sigma_{\mathrm{NG}}}
\newcommand{\chh}{\sigma_{\mathrm{H}_2}}

\newcommand{\fr}{{\mathrm{Fr}}}
\newcommand{\hh}{{\mathrm{H}_2}}
\renewcommand{\ng}{{\mathrm{NG}}}

\definecolor{rwth}   {RGB}{  0  84 159}
\definecolor{rwth-75}{RGB}{ 64 127 183}
\definecolor{rwth-50}{RGB}{142 186 229}
\definecolor{rwth-25}{RGB}{199 221 242}
\definecolor{rwth-10}{RGB}{232 241 250}

\definecolor{rot}   {RGB}{204   7  30}
\definecolor{orange}   {RGB}{246 168   0}
\definecolor{lila}   {RGB}{122 111 172}
\definecolor{bordeaux-50}{RGB}{205 139 135}



%% file: tex/abstract.tex
\begin{abstract}
	\noindent \textbf{Abstract.} 
%
    We prove an existence result for the steady state flow of gas mixtures on networks. The basis of the model are the physical principles of the isothermal Euler equation, coupling conditions for the flow and pressure, and the mixing of incoming flow at nodes. The state equation is based on a convex combination of the ideal gas equations of state for natural gas and hydrogen. 
    We analyze mathematical properties of the model allowing us to prove the existence of solutions in particular for tree-shaped networks and networks with exactly one cycle.
    Numerical examples illustrate the results and explore the applicability of our approach to different network topologies.
\end{abstract}



%% file: tex/introduction.tex
\section{Introduction and Motivation}
The transition to renewable energies is one of the most intensively discussed topics in science, with hydrogen being a key element that is set to play an increasingly important role in the global energy mix. Hydrogen offers numerous advantages. Firstly, it is a clean energy source that can be produced using renewable energy sources such as wind and solar power, and its combustion results only in water and heat, \cite{Bundesregierung2024, EuropeanCommission2020}. Secondly, hydrogen can be utilized to decarbonize sectors that are difficult to electrify, such as heavy industry \cite{Bundesregierung2024, EuropeanCommission2020}. Lastly, hydrogen serves as an energy carrier, making it a valuable resource for energy storage \cite{EuropeanCommission2020, USEnergy2017}. Given these benefits, there is no doubt that hydrogen will play an important role in the near future. 

However, the transition to renewable energies necessitates not only an appreciation of hydrogen's clean energy potential but also a deep understanding of how to integrate it into current energy systems. Even if hydrogen transport networks are planned all over the world (cf. \cite{HydrogenEU2022} for Europe, \cite{HydrogenASIA2023} for Asia, and \cite{HydrogenUSA2023} for U.S.A.), hydrogen is far from completely replacing natural gas as energy source. One key strategy in order to reduce the use of natural gas as energy source as well as to significantly lower the amount of emitted greenhouse gases, is the blending of hydrogen into natural gas networks.

A rigorous mathematical understanding of the hydrogen-natural gas blending processes is essential for the efficiency and safety of existing pipeline networks. Compared to the transport of pure hydrogen, the complexity lies in modelling and understanding the mixing of two gases. 
While the transport of pure hydrogen as well as pure natural gas in pipeline networks can be modeled by the isothermal Euler equations, cf. \cite{BandaHertyKlar2006,Gas-Modell:DomschkeHillerLangTischendorf2017,KochEtAl2015}, a different approach is needed for the modelling of gas mixtures. Since natural gas has been and still is widely used as an energy source for decades almost all over the world, researchers can rely on well-established regulations for natural gas networks and extensive research in this area. 
In \cite{SchmidtSteinbachWillert2014}, the authors discuss stationary states of the isothermal Euler equations and in \cite{GugatHanteHirschLeugering2015} the existence of unique stationary states was shown on arbitrary graphs. In \cite{BandaHertyKlar2006b, BandaHertyKlar2006} the authors analyze the transient model and discuss solutions to Riemann problems on arbitrary graphs including appropriate coupling conditions at the nodes, i.e., conservation of mass and pressure continuity. 
Further, the existence of semi-global Lipschitz continuous solutions of the isothermal Euler Equations on a network with compatible coupling conditions was shown in \cite{GugatUlbrich2018}, based on a fixed point iteration along the characteristic curves for Lipschitz continuous initial and boundary conditions. 

However, applying the isothermal Euler equations for both, hydrogen and natural gas in the same pipeline, does not consider collisions of molecules of different types, which leads to adulterated results \cite{Zlotnik2024}. Based on the analysis of chemical reacting fluid mixtures in \cite{BotheDreyer2015, MalekSoucek2022} we consider a model for the flow of a natural gas and hydrogen mixture under perfect mixing conditions (i.e., both gases have the same velocity and temperature), given by 
\begin{subequations}\label{eq:dynamicEquation}
    \begin{align}
        \dt \rhong(t,x) + \dx \qng(t,x) &= 0, \\
        \dt \rhohh(t,x) + \dx \qhh(t,x) &= 0, \\
        \dt q(t,x) + \dx \left[ p(\rhong(t,x), \rhohh(t,x)) + \frac{q^2(t,x)}{\rho(t,x)} \right] 
            & = - \frac{\lambda_\fr}{2D} \frac{q(t,x) \vert q(t,x) \vert}{\rho(t,x)},  
    \end{align} 
\end{subequations}
where $\rhong, \rhohh$ are the densities of natural gas and hydrogen, $\qng, \qhh$ are the flows of natural gas and hydrogen, $p, \rho, q$ are the pressure, flow and density of the mixture, $\lambda_\fr > 0$ is the pipe friction coefficient and $D$ is the pipe diameter. 
The pressure of the mixture depends on the densities by the state equation derived for example in \cite{Pfetsch2024} for the mixture of ideal gases and in \cite{Bermudez2022} for the mixture of real gases. The latter also assumes a non-isothermal dynamic.
%

The mixing model \eqref{eq:dynamicEquation} on networks also includes coupling conditions at the network junctions. For natural gas transport on networks, conservation of mass and continuity of the pressure are are typical choices, cf. \cite{GugatHanteHirschLeugering2015, GugatSchuster2018, SchusterStrauchGugatLang2022, Zlotnik2022c}. In \cite{BandaHertyKlar2006b}, the authors give an overview of coupling conditions for models of natural gas transport on networks. In the mixing model~\eqref{eq:dynamicEquation}, we also consider perfect mixing at the nodes, which implies that the outgoing hydrogen concentration is equal to the average weighted by the flow of all ingoing concentrations.

For the long time horizon planning of gas networks, dynamic models are often replaced by static models, that represent the steady states, cf. \cite{GugatHanteHirschLeugering2015, GugatSchuster2018, Zlotnik2024, SchmidtSteinbachWillert2014, Zlotnik2022c}. Simulation results for models based on \eqref{eq:dynamicEquation} and corresponding steady state models were presented, for example, in \cite{GuandaliniColbertaldoCampanari2017,Zlotnik2024,WitkowskiEtAl2018}. Considering natural gas and hydrogen as ideal gases, in \cite{GugatGiesselmann2024} the authors present a well-posedness result of \eqref{eq:dynamicEquation} under the assumption, that the hydrogen concentration is sufficiently low. 
To our best knowledge, neither the existence of steady states for the mixing model~\eqref{eq:dynamicEquation}, nor steady states for the mixing on networks were theoretically analyzed yet.
The only result that fits into this line of research is that for tree-shaped networks steady state solutions are unique if they exits, cf. \cite{Zlotnik2024}. However, proving uniqueness becomes more intricate for networks with cycles, since the monotonicity of the pressure, which guarantees the uniqueness for pure natural as transport, is unclear a priori due to the mixing.

In this paper, we present an analytic solution of the steady states of \eqref{eq:dynamicEquation}. Further, we analyze the existence and uniqueness of steady states for the hydrogen blended natural gas flow on networks. We discuss tree-shaped pipeline networks as well as networks that contain cycles. For pure gas flow the existence and uniqueness of steady states was analyzed in \cite{GugatHanteHirschLeugering2015} for ideal gases and in \cite{WintergerstDiss} for real gases. 
In the case of gas mixtures, the gas composition in a pipeline depends on the direction of the flow, which leads to discontinuities in the mixing model if the graph contains cycles. 
We present a novel existence result for steady states of the mixing model for networks with a cycle. The result is based on cutting edges of the cycle such that the model has a solution on the resulting tree-shaped network with new supply and demand nodes.
The challenge here lies in the discontinuity of the composition of the new supply and demand nodes, which underlines the difficulty compared to pure natural gas networks, where this discontinuity does not occur.

The article structures as follows: First, we introduce and motivate the gas flow model for mixtures on networks in \cref{sec: model equations}. We also show that certain properties of the gas flow of a single gas carries over to the gas flow of mixtures.
Then, in \cref{sec: existence}, we establish the existence of steady state solutions for tree-shaped networks and networks with a single cycle.
Finally, we conclude with a numerical study in \cref{sec: experiments}, which analyzes whether the idea of the existence result also extends to networks with arbitrarily many cycles and whether the solutions are unique.

%% file: tex/gas_flow_model.tex
\section{The Model Equations on Networks}\label{sec: model equations}
In this section, we derive a steady-state model for the hydrogen blended natural gas flow in pipeline networks. 
We further define suitable boundary and coupling conditions for the flow, the pressure and the composition of the gas mixture at the junctions of the network. 

Throughout this paper we assume that \(\G = (\V, \E) \) is a connected, directed graph with nodes \( v \in \V \) and edges \( e \in \E \subseteq \V \times \V\). Each edge \( e\) represents a pipe of length \( L^e\) and each node \( v\) represents a junction where pipes meet.
We consider two types of nodes: 
\begin{enumerate}[label=(\roman*)]
    \item \emph{inflow} or \emph{supply nodes} \( v \in \V_{<0}\), where gas enters the network,
    \item \emph{outflow} or \emph{demand nodes} \( v \in \V_{\ge0}\), where gas leaves the network or the mass is conserved,
\end{enumerate}
This classification of nodes satisfies \( \V = \V_{<0} \cup\ \V_{\ge0}\) with $\V_{<0}\cap\ \V_{\ge0} = \emptyset$. For every node $v \in \V$, the set of \emph{incoming} and \emph{outgoing edges} are given by 
\begin{equation*}
    \E_-(v) \coloneqq \{ e\in \E \mid e = (\bullet, v)\} 
    \quad \text{and} \quad 
    \E_+(v) \coloneqq \{ e\in \E \mid e = (v, \bullet)\},
\end{equation*}
respectively. The set containing all edges incident to \( v\) is given by  \( \E(v) = \E_-(v) \cup \E_+(v)\).
For edges \( e \in \E_-(v)\), the node \(v\) is their \emph{head} or \emph{end node} and for edges \( e \in \E_+(v)\), the node \( v\) is their \emph{foot} or \emph{start node}. In the following, we denote the head of an edge by \( h(e)\) and the foot of an edge by \( f(e)\), i.e., we have
\begin{equation*}
    e = (v,w) \quad \Leftrightarrow \quad f(e) = v \quad \text{and} \quad h(e) = w.
\end{equation*}
Further let \(A \in \R^{\vert \V \vert \times \vert \E \vert}\) be the so-called \emph{node edge incidence matrix}, whose entries are defined by
\begin{equation}\label{eq: incidence matrix}
	a(v,e) \coloneqq \begin{cases}
		-1, & \text{if } v = f(e), \text{ i.e., } v \text{ is the start node of } e,\\
		1, & \text{if } v = h(e), \text{ i.e., } v \text{ is the end node of } e,\\
		0, & \text{otherwise}.
	\end{cases}
\end{equation}

\subsection{The Isothermal Euler Equations for Mixtures} 
To model the flow of a mixture through a pipe \( e\), we consider the steady-state hydrogen blended gas flow model corresponding to the transient model \labelcref{eq:dynamicEquation}.
Additionally, we neglect the term \( q_e^2 / \rho_e\) in the differential equation for the pressure \( p_e\) compared to the transient model~\labelcref{eq:dynamicEquation} as the velocity of the gas flow is significantly smaller than the speed of sound. Thus, we are in a subsonic regime since the Mach number is sufficiently small.
Then, the subsonic steady-state hydrogen blended gas flow model reads:
\begin{subequations} \label{eq:stationaryEquation}
    \begin{align}
        \partial_x \vphantom{{\dx}} \qng^e(x) & = 0, && \text{for } x \in [0,L^e] \\
        \partial_x \vphantom{{\dx}} \qhh^e(x) & = 0, && \text{for } x \in [0,L^e]\\
        \partial_x p_e(\rhong^e(x), \rhohh^e(x)) & = - \frac{\lambda_\fr}{2 D} \frac{q_e(x) \vert q_e(x) \vert}{\rho_e(x)}, && \text{for } x \in [0,L^e]\label{eq: pressure ODE}
    \end{align} 
\end{subequations}
where \(\rhong^e\) and \(\rhohh^e\) denote the density and \(\qng^e \) and \(\qhh^e \) the flow of natural gas and hydrogen along the pipe.
The constant \( \lambda_\fr > 0 \) denotes the friction factor and \( D\) the diameter of the pipe. For the reader's convenience we assume the same friction and diameter for every pipe. 
Moreover, \( \rho_e\) describes the density, \( q_e\) the flow, \( \eta_e \) describes the composition, and \( p_e\) the pressure of the mixture along the pipe $e \in \mathcal{E}$. These quantities are defined as
\begin{equation*}
    \rho_e(x) = \rhohh^e(x) + \rhong^e(x), \quad q_e(x) = \qhh^e(x) + \qng^e(x), 
    \quad \eta_e(x) = \frac{\qhh^e(x)}{q_e(x)} 
    = \frac{\rhohh^e(x)}{\rho_e(x)},
\end{equation*}
while the latter holds due to the assumption that both gases have the same velocity.
%
%
For the pressure law we assume that the gases and their mixture are ideal which leads to
\begin{equation} \label{eq:pressureLaw}
    p_e(\rhong^e(x), \rhohh^e(x)) = \left[ \eta_e(x) \chh^2 + (1-\eta_e(x)) \cng^2 \right] \rho_e(x),
\end{equation}
where \( \chh^2 = R_{S,\hh} \, T\) is the speed of sound in hydrogen and \( \cng^2 = R_{S,\ng} \, T\) is the speed of sound in natural gas. Moreover, \(R_{S,\hh}\) and \(R_{S,\ng}\) are the specific gas constants of hydrogen and natural gas, respectively, and \( T \) the temperature of the gas.
Then, the solution to the system of ordinary differential equations~\labelcref{eq:stationaryEquation} describing the gas flow of a hydrogen-natural gas mixture along a pipe, reads:
\begin{subequations}\label{eq: pressure loss pipe}
    \begin{align}
        \qng^e(x) & = \const, \quad
        \qhh^e(x) = \const, \\
        p_e(x) & = \sqrt{ \left[ p_e^2(0) - \frac{\lambda_\fr}{D} \big( \eta_e \chh^2 + (1-\eta_e) \cng^2 \big) q_e \vert q_e \vert x \right] }. 
    \end{align}
\end{subequations}
The flow of hydrogen and natural gas is constant along the pipe and the pressure is uniquely defined by the pressure at the beginning and the end of the pipe. If $p_e(0)$ is given for an edge $e \in \E$, the the critical length of the pipe is given by
\begin{equation}\label{eq: L crit}
    L^e_{\text{crit}} = \frac{p_e^2(0)}{\frac{\lambda_\fr}{D} \big( \eta_e \chh^2 + (1-\eta_e) \cng^2 \big) q_e \vert q_e \vert},
\end{equation}
which means that if $ L^e > L^e_{\text{crit}}$, then a solution of $p_e$ does not exist in $\mathbb{R}$. For the sake of better readability, we omit the dependence of the spatial variable \( x\) in the following. 

\subsection{Boundary and Coupling Conditions for the Nodes}
As gas is injected or withdrawn from the network at the nodes, we introduce the load vector $\bb \in \mathbb{R}^{\vert \V \vert}$ with entries
\begin{equation}\label{eq: loads}
    b_v \, \begin{cases}
        < 0 & \text{if } v \in \V_{<0}, \text{ which means $v$ is a supply node,}\\
    \geq 0 &\text{if } v \in \V_{\ge 0}, \text{ which means $v$ is a demand node,}
    \end{cases}
\end{equation}
describing the supply and demand of a node \(v\). The case \( b_v = 0\) implies, that gas is neither injected nor withdrawn from node $v$.


For the boundary conditions, we provide the load \( b_v\) at each node \( v \in \V\). At supply nodes \( v\in \V_{<0}\), we additionally specify the supply composition \( \zeta_v\). Furthermore, we provide the pressure \( p^{\ast}\) at an arbitrary yet fixed node \( v^{\ast} \in \V\), i.e., we set \( p_{v^\ast} = p^{\ast}\).

For the coupling conditions, we make three assumptions to describe the gas flow across nodes: 
\begin{enumerate}[label=(\roman*)]
    \item The mass of the mixture and the mass of the individual gases are conserved,
    \item the  gas mixes instantaneously and perfectly, and 
    \item  the pressure is continuous across a node.
\end{enumerate}

The first assumption, the conservation of mass at a node $v \in \mathcal{V}$ means, that the amount of gas entering node $v$ must be equal to the amount of gas leaving node $v$ including the supply and demand of the node, respectively, i.e., we have
\begin{equation*}
    \sum_{e \in \E_+(v)} q_e - \sum_{e \in \E_-(v)} q_e = b_v.
\end{equation*}
Using the incidence matrix \(A\) defined in \cref{eq: incidence matrix}, the conservation of mass of the mixture at all nodes \( v \in \V\), is given by
\begin{equation}\label{eq:massConservation}
    \sum_{e \in \E(v)} a(v,e) q_e  = b_v \quad \text{for all} \quad v \in \V 
    \qquad \Leftrightarrow \qquad A \q = \bb,
\end{equation}
where \( \q \in \R^{\vert \E \vert}\) is the vector of (constant) flows at the edges. 

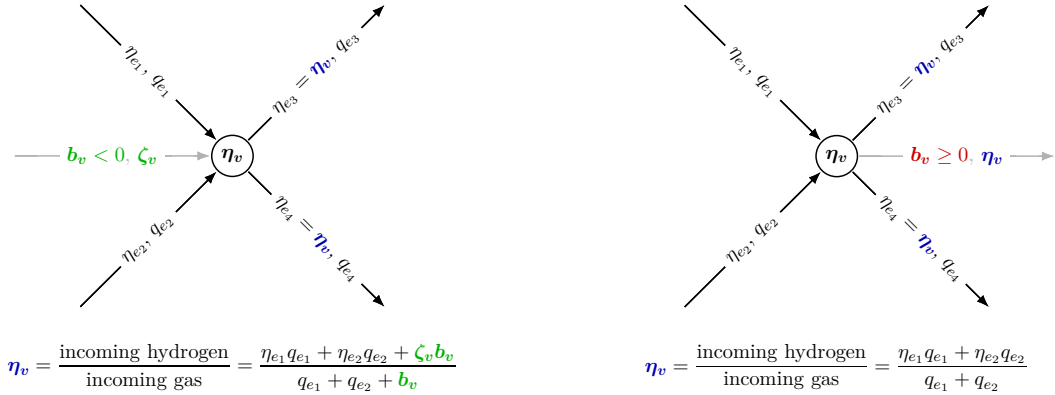
\begin{figure}[ht]
	\centering
    \begin{subfigure}{.475\textwidth}
        \centering
        \input{./tikz/mixing_at_nodes.tex}
    \end{subfigure}\hfill
    \begin{subfigure}{.475\textwidth}
        \centering
        \input{./tikz/consumption_at_node.tex}
    \end{subfigure} 
	\caption{The mixing of gas for a supply node (left) and a demand node (right). The supply and demand of a node can be seen as an invisible pipe with fixed flow.}
	\label{fig: mixing at nodes}
\end{figure}

\begin{remark}
    Due to the conservation of mass and the definition of the incidence matrix \(A\), the loads \(b_v\) must sum up to zero, i.e., we have
    \begin{equation}\label{eq: loads sum to zero}
        \sum_{v \in \V} b_v = 0.
    \end{equation}
    Hence, it is sufficient to provide the loads \(b_v\) for \( \vert \V \vert - 1\) nodes as the missing load is automatically determined by \cref{eq: loads sum to zero}.
\end{remark}
Next, we discuss the perfect mixing at the nodes, which implies that the incoming gas compositions at a node $v$ are weighted by the respective incoming flows and are then averaged and that the composition of all outgoing flows is equal. 
Hence, we introduce the nodal composition \( \eta_v \) of the gas at a node \( v \in \V\): 
\begin{equation*}
    \eta_v = \frac{\text{incoming hydrogen}}{\text{outgoing gas}}
        = \frac{\text{incoming hydrogen}}{\text{incoming gas}},
\end{equation*}
while the latter holds due to the conservation of mass.
In \cref{fig: mixing at nodes}, we provide an example on how to compute the nodal composition at a node \(v\) depending on whether the node is a supply or a demand node.

\begin{remark}
    At a supply node \( v\in \V_{<0}\), the nodal composition \( \eta_v\) differs from the composition \(\zeta_v\) of the supply because the supply at the node mixes with the inflow from incident pipes first, resulting the the nodal composition \( \eta_v\).
    As the individual compositions of the inflow differs for every pipe in general, their weighted average, the composition \( \eta_v\), differs from the supply composition \( \zeta_v\), cf. \cref{fig: mixing at nodes}.
    The relation \( \eta_v = \zeta_v \) only holds for supply nodes \( v \in \V_{<0} \) satisfying \( \vert \E(v) \vert = 1\).
\end{remark}

As the flow direction might differ from the orientation of the edge, it is unclear a priori which pipes carry incoming and which pipes carry outgoing gas. However, we can use the incidence matrix to determine whether a given flow \( q_e \) of an edge \(e \in \E(v) \) is directed in or away from the node \( v \):
\begin{equation}\label{eq: flow vs edge}
	\begin{aligned}
		a(v,e)\ q_e & \le 0 \quad \Leftrightarrow \quad \text{Gas in pipe $e$ flows towards node } v. \\
		a(v,e)\ q_e & \ge 0 \quad \Leftrightarrow \quad \text{Gas in pipe $e$ flows away from } v. 
	\end{aligned}	
\end{equation}
We provide a visualization of the condition \eqref{eq: flow vs edge} in \cref{fig: flow vs edge}. For a number $\alpha$, we define \( \alpha^+ = \max\{\alpha, 0\}\) and \( \alpha^- = \max\{-\alpha, 0\} \). Then, the nodal composition at a node \( v \in \V \) is given by
\begin{align} \label{eq:compositionConservation}
    \eta_v  
    = \frac{\sum_{e \in \E(v)} \eta_e \cdot \big(a(v,e) q_e\big)^+ + \zeta_v b_v^-}{\sum_{e \in \E(v)} \big(a(v,e) q_e\big)^+ + b_v^-
    },
\end{align}
where \( \eta_e\) is the composition along the pipe \( e\).
There is an explicit relation between the nodal composition \( \eta_v\) and the composition along a pipe \( \eta_e\) since the nodal composition \( \eta_v\) is transferred to all pipes \( e \in \E(v)\) carrying outgoing flow. More specifically, we have
\begin{equation}\label{eq: composition edge}
    \eta_e = \begin{cases}
        \eta_{f(e)} & \text{if } q_e \ge 0, \\
        \eta_{h(e)} & \text{if } q_e < 0,
    \end{cases}
    \quad \text{for all} \quad e \in \E,
\end{equation}
which means that the composition along the edge \( e \) is given by the nodal composition at the start node \( f(e)\) of the edge if the flow \(q_e \) is directed away from \( f(e) \), and by nodal composition at the end node \( h(e)\) if the flow \( q_e\) is directed away from \( h(e) \).
The relation~\labelcref{eq: composition edge} allows us to express the computation of the nodal composition \( \eta_v\) fully in terms of nodal compositions, cf. \cref{eq:compositionConservation}:
\begin{equation*}
	\eta_v = \frac{\sum_{e \in \E(v)} \eta_{w_e}\left( a(v,e) q_e \right)^+ + \zeta_v b_v^-}{\sum_{e \in \E(v)} \left( a(v,e) q_e \right)^+ + b_v^-}
		\quad \text{where} \quad w_e = \begin{cases}
			f(e), & \text{if } v = h(e), \\
			h(e), & \text{if } v = f(e).
		\end{cases}
\end{equation*}
The node \( w_e\) of an edge \( e \in \E(v)\) incident to \( v \), is the node that is not \( v\), i.e., if \( v\) is the start node of \(e\) then \( w_e\) is the end node of \( e\) and if \( v\) is the end node then \( w_e\) is the start node.
\begin{figure}[ht]
	\centering
    \begin{subfigure}{0.475\textwidth}
        \centering
        \input{./tikz/edge_vs_flow_direction}
	    \caption{The flow \(q_e\) is directed into a node \( v\), i.e., an inflow, if and only if \( a(v,e) q_e \ge 0\).}
    \end{subfigure} \hfill
    \begin{subfigure}{0.475\textwidth}
        \centering
        \input{./tikz/flow_directed_into_node}
	    \caption{The flow \(q_e\) is directed away from a node \( v\), i.e., an outflow, if and only if \( a(v,e) q_e \le 0\).}
    \end{subfigure} 
    \caption{The four possible relations between flow and edge direction, cf. \cref{eq: flow vs edge}.}
	\label{fig: flow vs edge}
\end{figure}

At last, we discuss the continuity of the pressure across a node. This means, that for all nodes $v \in \mathcal{V}$ and all incident edges \( e, \tilde{e} \in \E(v)\) the following condition holds:
\begin{equation*}
    p_e(x_e) = p_{\tilde{e}}(x_{\tilde{e}}) \quad \text{where} \quad 
    x_e = \begin{cases}
        0, & e \in \E_-(v), \\
        L^e, & e \in \E_+(v).
    \end{cases}
\end{equation*}
As the pressure loss along a pipe is defined by \cref{eq: pressure loss pipe}, we ensure the continuity of the pressure across a node by introducing the nodal pressure \( p_v\) for each node \( v\). Then, \cref{eq: pressure loss pipe} is given by
\begin{equation*}
    p_{h(e)}^2 - p_{f(e)}^2 = -\frac{\lambda_\fr}{D} \left( \eta_e \sigma_\hh^2 + (1- \eta_e) \sigma_\ng^2 \right) q_e \vert q_e \vert L^e 
        \quad \text{for all} \quad e \in \E.
\end{equation*}

As the composition \( \eta_e\) on an edge is discontinuous in \( q_e = 0\), cf. \cref{eq: composition edge}, the pressure \(p_e\) is discontinuous in \( q_e = 0\) as well. However, for pure natural gas transport, the pressure continuously depends on the nodal pressure and on the gas flow, see for example \cite{Gugat2016, WintergerstDiss}. 
As the continuity will be crucial for our further analysis, we rewrite by expressing the pressure loss in terms of the nodal composition instead of the composition on the edges.
First, we rewrite   
\begin{align*}
    \eta_e \sigma_\hh^2 + (1-\eta_e) \sigma_\ng^2
        & = \frac{1}{2} \left( \frac{\vert q_e \vert}{q_e} + 1 \right) \sigma(\eta_{f(e)}) 
              - \frac{1}{2} \left( \frac{\vert q_e \vert}{q_e} - 1 \right) \sigma(\eta_{h(e)}) && \text{for }  e \in \E,\\
        \text{where} \quad \sigma(\eta_v) & \coloneqq \eta_v \chh^2 + (1-\eta_v)\cng^2 && \text{for } v \in \V.            
        \vphantom{\frac{\vert q_e \vert}{q_e}}
\end{align*}
Then, for each edge \( e \in \E\) we define
\begin{equation}\label{eq: sigma}
    \tilde{\sigma}(\eta_{f(e)}, \eta_{h(e)}, q_e) \coloneqq 
        -\frac{\lambda_\fr}{D} L^e \left[\frac{1}{2} \left( \frac{\vert q_e \vert}{q_e} + 1 \right) \sigma(\eta_{f(e)}) 
              - \frac{1}{2} \left( \frac{\vert q_e \vert}{q_e} - 1 \right) \sigma(\eta_{h(e)}) \right],
\end{equation}
which allows to derive an expression of the pressure loss along a pipe in terms of the nodal compositions with removable singularity:
\begin{equation}\label{eq: pressure loss continuous}
    p_{h(e)}^2 - p_{f(e)}^2 = \tilde{\sigma}(\eta_{f(e)}, \eta_{h(e)}, q_e) q_e \vert q_e\vert \quad \text{for all} \quad  e \in \E.
\end{equation}
In the following, we will use \cref{eq: pressure loss continuous} to describe the pressure loss along a pipe.

\subsection{The Full Model on Graphs}
The full steady state gas flow model for gas mixtures consists of the pressure loss along a pipe, the conservation of mass across nodes, the mixing of the gases at nodes, and the boundary conditions defined above. 

For each node \( v \in \V \), we specify its load \( b_v\), which must satisfy \( \sum_{v \in \V} b_v = 0\).
Further, we provide the supply composition \( \zeta_v\) at supply nodes \( v \in \V_{<0}\), and fix a single nodal pressure $p_{v^{\ast}} = p^{\ast}$ at an arbitrary node $v^{\ast} \in \mathcal{V}$. For these given values, the gas flow model for mixtures is given by
\begin{subequations}\label{eq: gas flow}
	\begin{align}
		p_{h(e)}^2 - p_{f(e)}^2 & = \tilde{\sigma}(\eta_{f(e)}, \eta_{h(e)}, q_e) q_e \vert q_e\vert
            && \text{for all } e = (f(e),h(e)) \in \E, \vphantom{\sum_{\E^\inn}} \label{eq: pressure loss}\\
		A \q  & = \bb
            && \text{for all } v \in \V, \label{eq: mass conservation} \\
		\eta_v  & = \frac{ \sum_{e \in \E(v)} \eta_e \cdot (a(v,e) q_e)^+ + \zeta_v b_v^- }{ \sum_{e \in \E(v)} (a(v,e) q_e)^+ + b_v^- }
            && \text{for all } v \in \V. \label{eq: mixing} 
	\end{align}
\end{subequations}
The first equation~\labelcref{eq: pressure loss} describes the pressure loss along the pipes and guarantees the pressure continuity across nodes.
The second equation~\labelcref{eq: mass conservation} ensures the mass conservation at the nodes and describes the injection and withdrawal of gas from the network at the nodes. 
The last equation~\labelcref{eq: mixing} describes the perfect mixing at the nodes and guarantees the conservation of mass for the individual gases, i.e., natural gas and hydrogen, across nodes.
\begin{remark}
A small computation shows that the mixture model~\eqref{eq: gas flow} corresponds to the one presented in~\cite{Pfetsch2024} even though the derivations are different at a first glance.
\end{remark}

\subsection{Properties of the Gas Flow}
For gas flow with only a single gas, the gas cannot flow in cycles in networks without compressor stations, cf. \cite{GugatHanteHirschLeugering2015}. 
This property carries over to the gas flow of mixtures and we will need this property later when proving that the mixture model~\labelcref{eq: gas flow} admits at least one solution.
To prove that the mixture model~\labelcref{eq: gas flow} does not allow circular flows, we determine the pressure loss along a path connecting two nodes first as a cycle is a path with identical start and end node.

\begin{lemma}[Pressure Loss Along a Path] \label{lemma: pressure loss path}
	Let \( \G = (\V, \E) \) be a connected graph and let \( (\P, \E_\P)\) be a path starting at node \( v_0\) and ending at node \( v_k\) where \( \P = \{v_0, \ldots, v_k\}\) and \( \E_\P\) contains the edges connecting two consecutive nodes \( v_i\) and \( v_{i+1}\).
	Then the pressure loss along the path is given by
	\begin{equation*}
		p_{v_k}^2 - p_{v_0}^2 = \sum_{e \in \E_\P}
		\tilde{\sigma}(\eta_{h(e)}, \eta_{f(e)}, q_{e}) a(v_e, e) q_{e} \vert a(v_e, e) q_e \vert,
	\end{equation*}
	where the node \( v_e\) is the node of edge \( e\) that is closer to the start node \(v_0\) of the path, i.e., \( v_e = f(e) \) if  \(e = (v_i, v_{i+1})\) and \(v_e = h(e)\) if \(e = (v_{i+1},v_i)\). 

\end{lemma}
\begin{proof}
    We prove the claim by rewriting the pressure loss between two consecutive nodes \( v_k \) and \( v_{k+1}\) of the path, i.e., the pressure loss along the pipe connecting \( v_k\) and \(v_{k+1}\):
    \begin{align*}
		p_{v_{k+1}}^2 & = \begin{cases} 
			p_{v_{k}}^2 - \tilde{\sigma}(\eta_{v_{k}}, \eta_{v_{k+1}}, q_{e})
			q_{e} \vert q_{e} \vert, & \text{if } e = (v_{k}, v_{k+1}) \\
			p_{v_{k}}^2 + \tilde{\sigma}(\eta_{v_{k+1}}, \eta_{v_{k}}, q_{e})
			q_{e} \vert q_{e} \vert , & \text{if } e = (v_{k+1}, v_{k})
		\end{cases} \\
		& =  \begin{cases} 
			p_{v_{k}}^2 + \tilde{\sigma}(\eta_{f(e)}, \eta_{h(e)}, q_{e}) 	
			a(v_{k},e) q_{e} \vert a(v_{k},e) q_{e} \vert  , & \text{if } e = (v_{k}, v_{k+1}) \\
			p_{v_{k}}^2 + \tilde{\sigma}(\eta_{f(e)}, \eta_{h(e)}, q_{e}) 	
			a(v_{k+1},e) q_{e} \vert a(v_{k+1},e) q_{e} \vert, & \text{if } e = (v_{k+1}, v_{k})
		\end{cases} \\
		& = \hspace{8pt} p_{v_{k}}^2 + \tilde{\sigma}(\eta_{f(e)}, \eta_{h(e)}, q_{e}) a(v_e,e) q_{e} \vert a(v_e,e) q_{e} \vert,
	\end{align*}
    Then, the claim follows by an induction over the path length \( \vert \E_\P \vert\). \qedhere

\end{proof}

\begin{figure}[htbp]
	\centering
	\input{./tikz/flow_along_path.tex}
	\caption{Illustration of a circular flow moving clockwise \glqq starting\grqq{} from node \( v_0 \). The edge direction is marked in black and the flow direction in blue.}
	\label{fig: flow along path}
\end{figure}
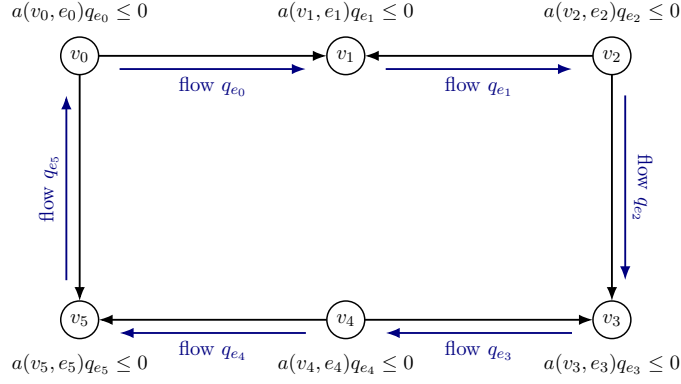

A cycle is a path with identical start and end node. Since gas flowing in cycles means that the gas always flows either away from or into the node when traversing through the cycle starting from a designated start node, we use \cref{eq: flow vs edge} to characterize circular flows, cf. \cref{fig: flow along path} for a visualization. The designated \glqq start node\grqq{} can be any node in the cycle. 

\begin{korollar}\label{corollary: no circular flow}
	Let \( \G = (\V, \E)\) be a connected graph. Then a solution to the gas flow model~\labelcref{eq: gas flow} does not allow circular flows, i.e., there exists no circle \( (\C, \E_\C)\) in the network satisfying one of the following two conditions:
	\begin{enumerate}[label=(\roman*)]
		\item \( a(v_e, e)\ q_e \ge 0 \quad \text{for all} \quad e \in \E_\C 
		\quad \text{and} \quad a(v_e, e)\ q_e > 0 \quad \text{for at least one} \quad e \in \E_\C\), \label{item: condition a}
		\item \( a(v_e, e)\ q_e \le 0 \quad \text{for all} \quad e \in \E_\C
		\quad \text{and} \quad a(v_e, e)\ q_e < 0 \quad \text{for at least one} \quad e \in \E_\C\),\label{item: condition b}
	\end{enumerate}
	where \( v_e \in \V\) is incident to \(e\) and the node closest to the start node while ignoring the edge between start and end node, cf. Lemma~\ref{lemma: pressure loss path}.
\end{korollar}
\begin{proof}
	We assume the contrary and prove for the first case~\labelcref{item: condition a}. A circle \(  (\C, \E_\C) \) is a path with identical start and end node. Thus, applying \cref{lemma: pressure loss path} yields:
	\begin{align*}
		0 = p_{v_k}^2 - p_{v_0}^2 = \sum_{e \in \E_\C} 
		\tilde{\sigma} \underbrace{\left(\eta_{f(e)}, \eta_{h(e)}, q_{e} \right)}_{\ge 0} \underbrace{\vphantom{\left(  \eta_{v_\out^{e}} \right)} a(v_e, e) q_e}_{\ge 0} 
		\underbrace{\vphantom{\left(  \eta_{v_\out^{e}} \right)}\vert a(v_e, e) q_{e} \vert}_{\ge 0} \stackrel{(\ast)}{>} 0,
	\end{align*}
	which is a contradiction. The strict inequality \( (\ast)\) holds because there exits at least one edge \( e \in \E_\C\) satisfying \( a(v_e, e) q_e > 0\).
	The proof for the second case~\labelcref{item: condition b} follows analogously.
\end{proof}

%% file: tikz/mixing_at_nodes.tex
\scalebox{0.7}{
	\begin{tikzpicture}[
		> = {Latex[scale=1]}, 
		node distance=4cm, 
		thick, 
		]

        \pgfmathsetmacro{\L}{3}
        \pgfmathsetmacro{\h}{1}
		
		\tikzstyle{inner state}=[
		fill=white, 
		draw=black,
		circle,
		inner sep=3pt,
		outer sep=1pt,
		minimum size=10pt
		]
		
		\tikzstyle{empty state}=[]
		
		\node[inner state] (0) at (0,0) {$\bm{\eta_v}$};       
		\node[empty state] (1) at (-\L, \L) {}; 
		\node[empty state] (2) at (-\L, -\L) {};  
		\node[empty state] (3) at (\L, \L) {}; 
		\node[empty state] (4) at (\L, -\L) {};
		\node[empty state] (5) at (-1.41*\L, 0) {};     
		
		\path[->] (1) edge node[midway, sloped, below left] {} (0); 
		\path[->] (2) edge node[midway, sloped, above left] {} (0); 
		\path[->, loosely dashed, dash pattern={on 12pt off 12pt}, black!30] (5) edge node[midway, sloped, below right] {} (0); 
	
		\path[->] (0) edge node[midway, sloped, above right] {} (3); 
		\path[->] (0) edge node[midway, sloped, below right] {} (4); 
		
		\begin{scope}[nodes={midway, sloped, fill=white}]
			\path[->] (1) edge node {$\eta_{e_1}, \, q_{e_1}$} (0);
			\path[->] (2) edge node {$\eta_{e_2}, \, q_{e_2}$} (0); 
			
			\path[->] (0) edge node {$\eta_{e_3} = \textcolor{blue!70!black}{\bm{\eta_v}}, \, q_{e_3}$} (3); 
			\path[->] (0) edge node {$\eta_{e_4} = \textcolor{blue!70!black}{\bm{\eta_v}}, \, q_{e_4}$} (4);

			\path[->, black!30] (5) edge node {$\textcolor{green!70!black}{\bm{b_v} < 0}, \, \textcolor{green!70!black}{\bm{\zeta_v}}$} (0);
		\end{scope}

		\node at (0,-\L-\h) {$\textcolor{blue!70!black}{\bm{\eta_v}} =
			\displaystyle{ \frac{\text{incoming hydrogen}}{\text{incoming gas}} }
			= \displaystyle{
				\frac{\eta_{e_1} q_{e_1} + \eta_{e_2} q_{e_2} + \textcolor{green!70!black}{\bm{\zeta_v b_v}}} {q_{e_1} + q_{e_2} + \textcolor{green!70!black}{\bm{b_v}}}
		}$};
	\end{tikzpicture}
}

%% file: tikz/consumption_at_node.tex
\scalebox{0.7}{
	\begin{tikzpicture}[
		> = {Latex[scale=1]}, 
		thick, 
		]

        \pgfmathsetmacro{\L}{3}
        \pgfmathsetmacro{\h}{1}
		
		\tikzstyle{inner state}=[
		fill=white, 
		draw=black,
		circle,
		inner sep=3pt,
		outer sep=1pt,
		minimum size=10pt
		]
		
		\tikzstyle{empty state}=[]
		
		\node[inner state] (0) at (0,0) {$\bm{\eta_v}$};       
		\node[empty state] (1) at (-\L, \L) {}; 
		\node[empty state] (2) at (-\L, -\L) {};  
		\node[empty state] (3) at (\L, \L) {}; 
		\node[empty state] (4) at (\L, -\L) {};
		\node[empty state] (5) at (1.41*\L, 0) {}; 
		
		\path[->] (1) edge node[midway, sloped, below left] {} (0); 
		\path[->] (2) edge node[midway, sloped, above left] {} (0); 
		\path[->, loosely dashed, dash pattern={on 12pt off 12pt}, black!30] (0) edge node[midway, sloped, below right] {} (5); 
		
		\path[->] (0) edge node[midway, sloped, above right] {} (3); 
		\path[->] (0) edge node[midway, sloped, below right] {} (4); 
		
		\begin{scope}[nodes={midway, sloped, fill=white}]
			\path[->] (1) edge node {$\eta_{e_1}, \, q_{e_1}$} (0);
			\path[->] (2) edge node {$\eta_{e_2}, \, q_{e_2}$} (0); 
			
			\path[->] (0) edge node {$\eta_{e_3} = \textcolor{blue!70!black}{\bm{\eta_v}}, \, q_{e_3}$} (3); 
			\path[->] (0) edge node {$\eta_{e_4} = \textcolor{blue!70!black}{\bm{\eta_v}}, \, q_{e_4}$} (4);
			
			\path[->, black!30] (0) edge node {$\textcolor{red!80!black}{\bm{b_v} \ge 0}, \, \textcolor{blue!70!black}{\bm{\eta_v}}$} (5);
		\end{scope}

		\node at (0, -\L-\h) {$\textcolor{blue!70!black}{\bm{\eta_v}} = 
			\displaystyle{ \frac{\text{incoming hydrogen}}{\text{incoming gas}} }
				= \displaystyle{
					\frac{\eta_{e_1} q_{e_1} + \eta_{e_2} q_{e_2}} {q_{e_1} + q_{e_2} }
					}$};

	\end{tikzpicture}
}

%% file: tikz/edge_vs_flow_direction.tex
\scalebox{0.9}{
	\begin{tikzpicture}[
		> = {Latex[scale=1.0]}, 
		line width=1pt, 
		]
		
		\tikzstyle{empty}=[]
		
		\tikzstyle{state}=[
		draw=black,
		circle,
		inner sep=0pt,
		outer sep=0pt,
		minimum size=15pt,
		thick
		]
		
		\pgfmathsetmacro{\L}{4}
		\pgfmathsetmacro{\h}{0.25}
		\pgfmathsetmacro{\x}{4.5}
		
		
		
		\node[state] (v) at (0,0) {\(v\)}; 
		\node[state] (w) at (\L,0) {\(w\)};
		
		\node[empty] (v0) at (\h,-\h) {}; 
		\node[empty] (w0) at (\L-\h,-\h) {};

		\path[->] (v) edge node[midway, above] {$e=(v,w)$} (w);
		\path[->, blue!50!black] (v0) edge node[midway, below] {flow $q_e$} (w0);
		
		\begin{scope}[shift={(0,-2.5)}]
			\node[state] (v) at (0,0) {\(v\)}; 
			\node[state] (w) at (\L,0) {\(w\)};
			
			\node[empty] (v0) at (\h,-\h) {}; 
			\node[empty] (w0) at (\L-\h,-\h) {};

			\path[<-] (v) edge node[midway, above] {$e=(w,v)$} (w);
			\path[->, blue!50!black] (v0) edge node[midway, below] {flow $q_e$} (w0);
			
		\end{scope}
		

	\end{tikzpicture}
}

%% file: tikz/flow_directed_into_node.tex
\scalebox{0.9}{
	\begin{tikzpicture}[
		> = {Latex[scale=1.0]}, 
		line width=1pt, 
		]
		
		\tikzstyle{empty}=[]
		
		\tikzstyle{state}=[
		draw=black,
		circle,
		inner sep=0pt,
		outer sep=0pt,
		minimum size=15pt,
		thick
		]
		
		\pgfmathsetmacro{\L}{4}
		\pgfmathsetmacro{\h}{0.25}
		\pgfmathsetmacro{\x}{4.5}

				
		\node[state] (v) at (0,0) {\(v\)}; 
		\node[state] (w) at (\L,0) {\(w\)};
			
		\node[empty] (v0) at (\h,-\h) {}; 
		\node[empty] (w0) at (\L-\h,-\h) {};

		\path[->] (v) edge node[midway, above] {$e=(v,w)$} (w);
		\path[<-, blue!50!black] (v0) edge node[midway, below] {flow $q_e$} (w0);
			
		
    \begin{scope}[shift={(0,-2.5)}]
		\node[state] (v) at (0,0) {\(v\)}; 
		\node[state] (w) at (\L,0) {\(w\)};
		
		\node[empty] (v0) at (\h,-\h) {}; 
		\node[empty] (w0) at (\L-\h,-\h) {};

		\path[<-] (v) edge node[midway, above] {$e=(v,w)$} (w);
		\path[<-, blue!50!black] (v0) edge node[midway, below] {flow $q_e$} (w0);
    \end{scope}

\end{tikzpicture}
}

%% file: tikz/flow_along_path.tex
\scalebox{0.7}{
	\begin{tikzpicture}[
		> = {Latex[scale=1.0]}, 
		line width=1pt, 
		]
		
		\tikzstyle{empty}=[]
		
		\tikzstyle{state}=[
		draw=black,
		circle,
		inner sep=3pt,
		outer sep=0pt,
		minimum size=15pt,
		thick
		]
		
		\pgfmathsetmacro{\L}{5}
        \pgfmathsetmacro{\Lx}{5}
        \pgfmathsetmacro{\Ly}{5}
		\pgfmathsetmacro{\h}{0.5}
		\pgfmathsetmacro{\dx}{0.75}
		\pgfmathsetmacro{\dy}{0.25}

		\node[state] (v0) at (0,0) {\(v_0\)}; 
		\node[state] (v1) at (\Lx, 0) {\(v_1\)};
  		\node[state] (v2) at (2*\Lx,0) {\(v_2\)};
		\node[state] (v3) at (2*\Lx, -\Ly) {\(v_3\)};
		\node[state] (v4) at (\Lx,-\Ly) {\(v_4\)};
		\node[state] (v5) at (0, -\Ly) {\(v_5\)};

		\node[empty, above] (w0) at (0,\h) {\(a(v_0, e_0) q_{e_0} \le 0 \)}; 
		\node[empty, above] (w1) at  (\Lx, \h) {\(a(v_1, e_1) q_{e_1} \le 0 \)};
		\node[empty, above] (w2) at (2*\Lx, \h) {\(a(v_2, e_2) q_{e_2} \le 0 \)};
        \node[empty, below] (w3) at (2* \Lx, -\Ly -\h) {\(a(v_3, e_3) q_{e_3} \le 0 \)};
		\node[empty, below] (w4) at  (\Lx, -\Ly - \h) {\(a(v_4, e_4) q_{e_4} \le 0 \)};
		\node[empty, below] (w5) at  (0, -\Ly - \h) {\(a(v_5, e_5) q_{e_5} \le 0 \)};

		\path[->] (v0) edge node[midway, above] {} (v1);
        \path[->] (v2) edge node[midway, above] {} (v1);
        \path[->] (v2) edge node[midway, above] {} (v3);
        \draw[->] (v4) edge node[midway, above] {} (v3);
        \draw[->] (v4) edge node[midway, above] {} (v5);
        \draw[->] (v0) edge node[midway, above] {} (v5);

		\path[->, blue!50!black, sloped] (\dx, - \dy) edge node[midway, below] {flow $q_{e_0}$} ({\Lx - \dx, -\dy });	
        \path[->, blue!50!black] (\Lx+\dx,-\dy) edge node[midway, below] {flow $q_{e_1}$} (2*\Lx -\dx,-\dy);	
        \path[->, blue!50!black, sloped] (2*\Lx + \dy, -\dx) edge node[midway, above] {flow $q_{e_2}$} (2* \Lx + \dy, -\Ly +\dx);	
        \path[->, blue!50!black, sloped] (2*\Lx - \dx, -\Ly -\dy) edge node[midway, below] {flow $q_{e_3}$} (\Lx + \dx, - \Ly - \dy);	
        \path[->, blue!50!black, sloped] (\Lx - \dx, -\Ly -\dy) edge node[midway, below] {flow $q_{e_4}$} (\dx, - \Ly - \dy);	
        \path[->, blue!50!black, sloped] (-\dy, -\Ly +\dx) edge node[midway, above] {flow $q_{e_5}$} (-\dy, -\dx);	
             
	\end{tikzpicture}
}

%% file: tex/existence_network_tree.tex
\section{The Existence and Uniqueness of Steady States on Networks}\label{sec: existence}
In this section, we prove the existence of steady states on networks for the mixture model~\labelcref{eq: gas flow}.
The idea follows \cite{Gugat2016,WintergerstDiss} where the authors prove the existence of steady states on networks for the \glqq classic\grqq{} gas flow model involving only a single gas. 
Their proof is based on an induction over the number of cycles in a graph. To apply the induction assumption, one cuts an edge of a cycle to obtain a graph, that has one cycle less than the original one.

Thus, we first show the existence and uniqueness of solutions of the mixture model~\labelcref{eq: gas flow} on tree-shaped networks, i.e., network with no cycles, and establish the existence of steady-states for networks with cycles afterwards.

\subsection{Steady States on Tree-Shaped Networks}
The existence and uniqueness of steady states for tree-shaped networks is, similar to the single gas problem, straightforward as the argument only relies on the fact that for trees, the flow on the network is fully determined by the loads \( b_v\). 
However, we have one additional equation to describe the mixing at nodes. 
Provided that the mixing equation~\labelcref{eq: mixing} admits a unique solution, the proof to show the existence and uniqueness of steady states on tree-shaped networks is analogous to the proof of the steady-state gas flow model with a single gas.
 

We can compute the composition \( \eta_v\) at a node \(v\) if we know its incoming flows and their composition since the nodal composition \( \eta_v\) is given by the ratio of the incoming hydrogen and incoming gas, cf. \cref{eq: mixing}.
As the network is tree-shaped, we can sort the nodes such that gas always flows from nodes early in the ordering to nodes later in the ordering what allows us to compute the nodal compositions \( \eta_v\) in an inductive manner.
Such an ordering of the nodes is called a \emph{topological ordering} of a graph.

\begin{figure}[ht]
	\centering
	\begin{subfigure}[c]{0.4\textwidth}
		\centering
		\input{./tikz/sample_network_circle}
	\end{subfigure} 
	\begin{subfigure}[c]{0.58\textwidth}
		\centering
		\input{./tikz/topological_sorting}
	\end{subfigure}
	\caption{A directed acyclic graph (left) and a topological ordering of its nodes (right).}
	\label{fig: topological ordering}
\end{figure}
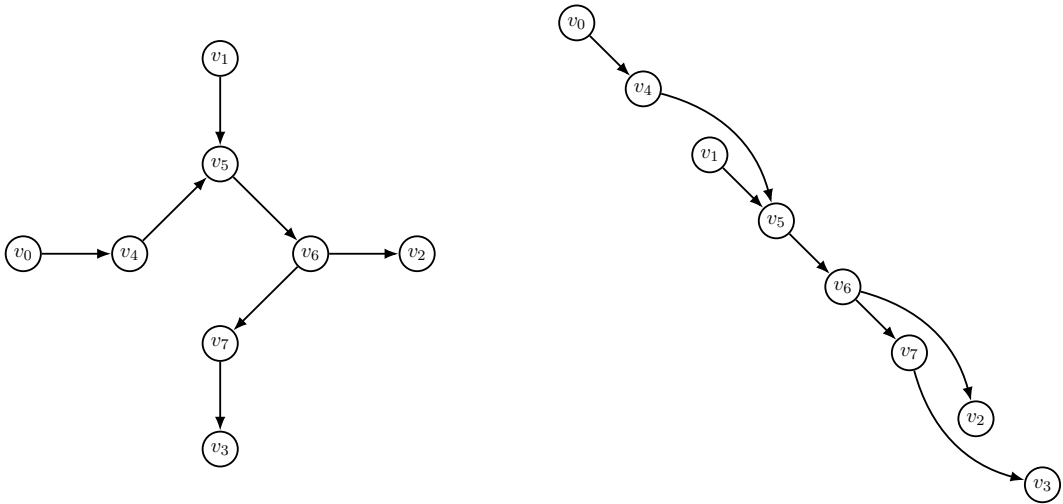

\begin{definition}[{Topological Ordering~\cite{Jensen2008}}]
	Let \( \G = ( \V, \E)\) be a directed, tree-shaped graph. Then a \emph{topological ordering} of a directed graph sorts the nodes of the graph into a sequence in which, for each directed edge \( (v,w)\), the start node \( v\) appears before the end node \( w\), cf. \cref{fig: topological ordering} for a visualization.
\end{definition}

To compute the compositions \( \eta_v \) for fixed flows \(q_e\), we require a topological ordering of the graph \( \G_\flow = (\V, \E_\flow)\), which has the same nodes as \( \G\) but edges according to the flow direction, i.e., its edges are defined by
\begin{align}\label{eq: G flow}
	(v,w) \in \E_\flow \quad :\Leftrightarrow \quad \begin{cases}
		(v,w) \in \E, & \text{if } q_e \ge 0, \\
		(w,v) \in \E, & \text{if } q_e < 0.
	\end{cases}
\end{align}
A topological ordering of the graph \( \G\) implies that pipes start at nodes early in the ordering and end at nodes later in the ordering, as the flow direction may differ from the pipe orientation. This can lead to gas flowing from nodes later in the ordering to modes earlier in the ordering, which is undesired. 
However, a topological ordering of \( \G_\flow\) guarantees that gas always flows from nodes early in the ordering to nodes later in the ordering, which is what we require.

\begin{figure}[ht]
	\centering
	\input{./tikz/G_flow}
	\caption{A directed graph \(\G\) in black, given flow \( q_e \) on the edges in blue (left) and the graph \( \G_\flow\) with edges according to the flow direction (right).}
	\label{fig: G vs G flow}
\end{figure}
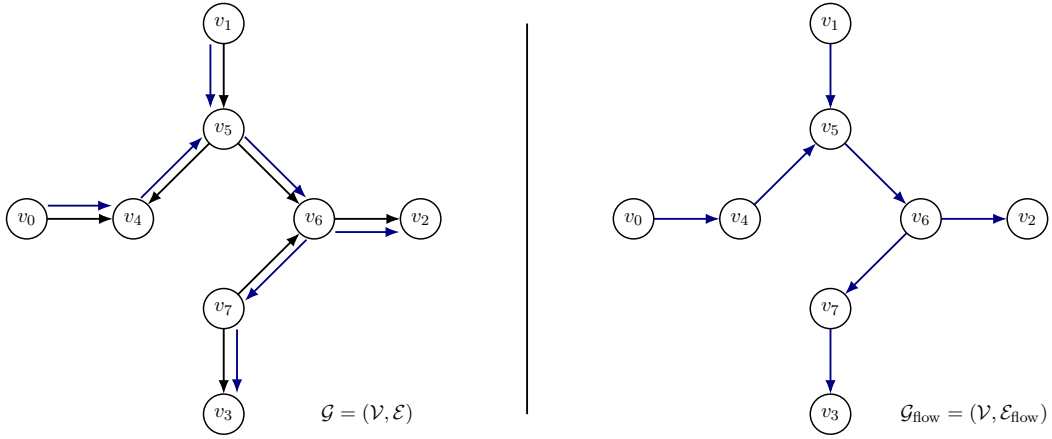

Note that the graph \( \G_\flow\) emerges from the original graph \( \G\) by flipping edges where the flow and edge direction differ. Hence, the graph \( \G_\flow\) depends on the flow \( q_e\) on the edges. We provide an example of the relation between the graph \( \G\) and the graph \( \G_\flow\) in \cref{fig: G vs G flow}.

\begin{remark}
The graph \( \G_\flow\) is tree-shaped because it has the same edge as \( \G\) up to a change in the edge orientation and thus, by \cite[Proposition 2.1.3]{Jensen2008} there exits a topological ordering of the nodes of \( \G_\flow\). 
\end{remark}

Therefore, we  are now equipped to show the existence of a solution to the mixing equation~\labelcref{eq: mixing} given fixed flows \( q_e\) on a network.

\begin{lemma}[Existence of the Composition]\label{lemma: mixing eq existence}
	Let \( \G = (\V, \E)\) be a directed graph and let \( q_e\) be flows on the graph satisfying mass conservation, cf.~\cref{eq: mass conservation}. Then the mixing equation~\labelcref{eq: mixing} has at least one solution \( \eta_v\).
\end{lemma} 

\begin{proof}
	The graph \( \G_\flow = (\V, \E_\flow)\) (cf. \cref{eq: G flow} for the definition) is acyclic.
	Thus, there exits a topological ordering of the nodes of \( \G_\flow\) and we can compute the nodal compositions \( \eta_v\) by induction over the position \(k\) in the topological ordering.
	
	The composition \( \eta_v\) for the node at position \( k=1\) is given by boundary data since the first node in a topological ordering has to be a supply node, i.e., \( b_v < 0\), where \( a(v,e) q_e \le 0\) for all \( e \in \E(v)\), and thus \( \eta_v = \zeta_v \).
		
	Now assume all nodal compositions up to position \(k\) have been determined. Then we compute the composition of the node at position \( k+1\) by invoking
	\begin{align*}
		\eta_v =  \frac{\sum_{e \in \E(v)} \eta_e \cdot \left( a(v,e) q_e \right)^+ + \zeta_v b_v^- }{\sum_{e \in \E(v)} \left( a(v,e) q_e \right)^+ + b_v^-}.
	\end{align*}	
	The incoming flow, i.e., edges with \( (a(v,e) q_e)^+ \neq 0\), only comes from nodes previous in the ordering. Assume otherwise then there exits an edge whose end point comes before the start point in the ordering which is a contradiction to the definition of a topological ordering.
\end{proof}

The proof of the existence of a solution relies on a topological ordering of the nodes. Since this ordering is in general non-unique, the question arises whether the solution depends on the topological ordering used to traverse through the nodes. Fortunately, this is not the case, since the solution to the mixing equation is unique.
	
\begin{lemma}[Uniqueness of the Composition]\label{lemma: mixing eq unique}
	Let \( \G = (\V, \E)\) be a directed, tree-shaped graph and let \( q_e\) be an flows on the graph satisfying the mass conservation, cf.~\cref{eq: mass conservation}. Then the solution \( \eta_v\) to the mixing equation~\labelcref{eq: mixing} is unique.
\end{lemma} 

\begin{proof}
	Follows directly from \cite[Theorem 1]{Zlotnik2024}. \qedhere
	%
		
\end{proof}

Finally, we prove the existence and uniqueness of solution to the steady-state gas flow model~\labelcref{eq: gas flow} for tree-shaped networks.

\begin{theorem}[{Existence of Solutions for Tree-Shaped Networks}]\label{thm: existence tree}
	Let \( \G= (\V, \E)\) be a tree-shaped network. Then the mixture model~\labelcref{eq: gas flow} admits a unique solution.
\end{theorem}

\begin{proof}
	As the network \( \G\) is tree-shaped, the incidence matrix \( A\) has full column rank \cite[Lemma 2.2]{Graphtheory} and thus, the flow equation \(A \q = \bb\) has a unique solution if and only if
	\begin{align*}
		\sum_{v \in \V} b_v = 0.
	\end{align*}
	Since the loads \( b_v \) satisfies this condition by assumption, the flow equation has a unique solution \( \q\) that depends only on the nodal flow rates and the topology of the network.
	Hence, by \cref{lemma: mixing eq existence} and \cref{lemma: mixing eq unique}, the mixing equation~\labelcref{eq: mixing} has a unique solution given the solution \( \q \)  of the flow equation.
 
	Together with the fixed nodal pressure \(p_v\) at the node \( v^{\ast} \), we rewrite \cref{eq: pressure loss}, which describes the pressure loss along the pipes using the incidence matrix. This yields the following linear system: 
	\begin{subequations}
		\begin{align}
			A^\top \p^2 & = \tilde{\sigmaa}(\etaa, \q) \quad \text{with} \quad p_{v^{\ast}} = p^{\ast} \\
			\Leftrightarrow \qquad A_{-v^{\ast}}^\top \p_{-v^{\ast}}^2
			& =  \underbrace{\tilde{\sigmaa}(\etaa, \q) - \left(p^{\ast} \right)^2\bm{a}_{v^{\ast}}}_{\text{fully determined}}, \label{eq: pressure reduced}
		\end{align}
	\end{subequations} 
	where \( A_{-v^{\ast}} \) and \( \p_{-v^{\ast}}\) are the incidence matrix \( A\) and the vector \( \p\), but without the row and entry corresponding \( v^{\ast}\), respectively, and \( \bm{a}_{v^{\ast}}\) denotes the row of \(A\) corresponding to \(v^{\ast}\). 
    Moreover, the function \(  \tilde{\sigmaa} \) is the vector-valued version of the function \(\tilde{\sigma}\) defined in~\cref{eq: sigma}, which means the vector \( \tilde{\sigmaa}(\etaa, \q) \in \R^{\vert \E \vert} \) has the entries:
    \begin{equation*}
        \left[\tilde{\sigmaa}(\etaa, \q) \right]_e = \tilde{\sigma}(\eta_{f(e)}, \eta_{h(e)}, q_e).
    \end{equation*} 
    Then, we obtain the nodal pressures \( p_v\) by solving the linear system~\labelcref{eq: pressure reduced} as the right side of this linear system is fully determined by the topology of the network, the flow \(q_e\) and composition \( \eta_e\) on the edges and the given nodal pressure at node \( v^{\ast}.\) 
    Moreover, the square matrix \( A_{-v^{\ast}}^\top\) is regular, cf. \cite{GotzesHeitschHenrionSchultz2016}, which completes the proof.
\end{proof}

%% file: tikz/sample_network_circle.tex
\scalebox{0.7}{
	\begin{tikzpicture}[
		> = {Latex[scale=1.0]}, 
		line width=1pt, 
		every node/.style={circle, draw, minimum size=10pt, inner sep = 2.5pt}
		]
		
		\pgfmathsetmacro{\L}{2}
		\pgfmathsetmacro{\sqrtL}{1.7}
        \pgfmathsetmacro{\hx}{0.25}
        \pgfmathsetmacro{\hy}{0.5}

		\node (v0) at (-\L,0) {\(v_0\)};
		\node (v4) at (0, 0) {\(v_4\)};
		\node (v1) at (\sqrtL, \L + \sqrtL) {\(v_1\)};
		\node (v5) at (\sqrtL, \sqrtL) {\(v_5\)};
		\node (v7) at (\sqrtL, -\sqrtL) {\(v_7\)};
		\node (v6) at (2*\sqrtL, 0) {\(v_6\)};
		\node (v2) at (\L + 2*\sqrtL, 0) {\(v_2\)};
		\node (v3) at (\sqrtL, - \L - \sqrtL) {\(v_3\)};

		\draw[->] (v0) -- (v4);
		\draw[->] (v4) -- (v5);
		\draw[->] (v1) -- (v5);
		\draw[->] (v5) -- (v6);
		\draw[->] (v6) -- (v7);
		\draw[->] (v7) -- (v3);
		\draw[->] (v6) -- (v2);
\end{tikzpicture}}

%% file: tikz/topological_sorting.tex
\scalebox{0.7}{
\begin{tikzpicture}[
	> = {Latex[scale=1.0]}, 
	line width=1pt, 
	every node/.style={circle, draw, minimum size=10pt, inner sep=2.5pt}
	]
	
	\pgfmathsetmacro{\L}{1.25}

	\node (v0) at (0,0) {\(v_0\)};
	\node (v4) at (\L,-\L) {\(v_4\)};
	\node (v1) at (2*\L,-2*\L) {\(v_1\)};
	\node (v5) at (3*\L,-3*\L) {\(v_5\)};
	\node (v6) at (4*\L,-4*\L) {\(v_6\)};
	\node (v7) at (5*\L,-5*\L) {\(v_7\)};
	\node (v2) at (6*\L,-6*\L) {\(v_2\)};
	\node (v3) at (7*\L,-7*\L) {\(v_3\)};

	\draw[->] (v0) -- (v4);
	\draw[->] (v4) to[bend left] (v5);
	\draw[->] (v1) -- (v5);
	\draw[->] (v5) -- (v6);
	\draw[->] (v6) -- (v7);
	\draw[->] (v7) to[bend right] (v3);
	\draw[->] (v6) to[bend left] (v2);
\end{tikzpicture}}

%% file: tikz/G_flow.tex
\scalebox{0.7}{
	\begin{tikzpicture}[
		> = {Latex[scale=1.0]}, 
		line width=1pt, 
		]
		
		\tikzstyle{empty}=[]
		
		\tikzstyle{state}=[
		draw=black,
		circle,
		inner sep=3.5pt,
		outer sep=0pt,
		minimum size=10pt,
		thick
		]
		
		\pgfmathsetmacro{\L}{2}
		\pgfmathsetmacro{\sqrtL}{1.7}
		\pgfmathsetmacro{\hx}{0.25}
		\pgfmathsetmacro{\hy}{0.25}

		\node[state] (v0) at (-\L,0) {\(v_0\)};
		\node[state] (v4) at (0, 0) {\(v_4\)};
		\node[state] (v1) at (\sqrtL, \L + \sqrtL) {\(v_1\)};
		\node[state] (v5) at (\sqrtL, \sqrtL) {\(v_5\)};
		\node[state] (v7) at (\sqrtL, -\sqrtL) {\(v_7\)};
		\node[state] (v6) at (2*\sqrtL, 0) {\(v_6\)};
		\node[state] (v2) at (\L + 2*\sqrtL, 0) {\(v_2\)};
		\node[state] (v3) at (\sqrtL, - \L - \sqrtL) {\(v_3\)};
		
		\node[align=right] (graph) at (2*\sqrtL+0.5*\L, - \L - \sqrtL) {\(\G=(\V, \E)\)};

		\draw[->] (v0) -- (v4);
		\draw[->] (v5) -- (v4);
		\draw[->] (v1) -- (v5);
		\draw[->] (v5) -- (v6);
		\draw[->] (v7) -- (v6);
		\draw[->] (v7) -- (v3);
		\draw[->] (v6) -- (v2);

		\node[] (w0) at (-\L+\hx,\hy) {};
		\node[] (w4) at (-\hx, \hy) {};
		\node[] (w1) at (\sqrtL-\hx, \L + \sqrtL -\hy) {};
		\node[] (w5) at (\sqrtL-\hx, \sqrtL+\hy) {};
		\node[] (w7) at (\sqrtL+\hx, -\sqrtL-\hy) {};
		\node[] (w6) at (2*\sqrtL+\hx, -\hy) {};
		\node[] (w2) at (\L + 2*\sqrtL -\hx, -\hy) {};
		\node[] (w3) at (\sqrtL+\hx, - \L - \sqrtL + \hy) {};
		
		\node[] (ww4) at (0, \hy) {};
		\node[] (ww5) at (\sqrtL-\hy, \sqrtL) {};
		\node[] (ww7) at (\sqrtL+\hx, -\sqrtL) {};
		\node[] (ww6) at (2*\sqrtL, -\hy) {};
		
		\node[] (www5) at (\sqrtL+\hx, \sqrtL) {};
		\node[] (www6) at (2*\sqrtL, \hy) {};
		
		\draw[->, blue!50!black] (w0) -- (w4);
		\draw[->, blue!50!black] (ww4) -- (ww5);
		\draw[->, blue!50!black] (w1) -- (w5);
		\draw[->, blue!50!black] (www5) -- (www6);
		\draw[->, blue!50!black] (ww6) -- (ww7);
		\draw[->, blue!50!black] (w7) -- (w3);
		\draw[->, blue!50!black] (w6) -- (w2);
		
		\draw[-] (\L + 2*\sqrtL +2 ,\L + \sqrtL) -- (\L + 2*\sqrtL +2 ,-\L - \sqrtL);

		\begin{scope}[shift={(2*\L + 2*\sqrtL + 4, 0)}]
			\node[state] (v0) at (-\L,0) {\(v_0\)};
			\node[state] (v4) at (0, 0) {\(v_4\)};
			\node[state] (v1) at (\sqrtL, \L + \sqrtL) {\(v_1\)};
			\node[state] (v5) at (\sqrtL, \sqrtL) {\(v_5\)};
			\node[state] (v7) at (\sqrtL, -\sqrtL) {\(v_7\)};
			\node[state] (v6) at (2*\sqrtL, 0) {\(v_6\)};
			\node[state] (v2) at (\L + 2*\sqrtL, 0) {\(v_2\)};
			\node[state] (v3) at (\sqrtL, - \L - \sqrtL) {\(v_3\)};
			
			\node[align=right] (graph) at (2*\sqrtL+0.5*\L, - \L - \sqrtL) {\(\G_\flow=(\V, \E_\flow)\)};

			\draw[->, blue!50!black] (v0) -- (v4);
			\draw[->, blue!50!black] (v4) -- (v5);
			\draw[->, blue!50!black] (v1) -- (v5);
			\draw[->, blue!50!black] (v5) -- (v6);
			\draw[->, blue!50!black] (v6) -- (v7);
			\draw[->, blue!50!black] (v7) -- (v3);
			\draw[->, blue!50!black] (v6) -- (v2);
			
		\end{scope}
\end{tikzpicture}}

%% file: tex/existence_network_with_cycle.tex
\subsection{Steady States on Networks with One Cycle}\label{sec: existence single cycle}
In this section we discuss the existence of solutions to the mixture model~\labelcref{eq: gas flow} for networks with cycles. In particular, we prove the main result \cref{thm: existence cycle}, the existence of solutions for networks with a single cycle. 
Moreover, we also comment on whether the result extends to networks with active elements such as compressors.

\begin{theorem}[Existence of Solutions for Networks with One Cycle]\label{thm: existence cycle}
	Let \( \G = (\V, \E) \) be a directed graph with exactly one cycle. Then the mixture model~\labelcref{eq: gas flow} admits at least one solution.
\end{theorem}

The key idea to prove \Cref{thm: existence cycle} is to cut an edge of the cycle, which results in a tree-shaped network for which we know that a unique solution exists (\cref{thm: existence tree}).
The cut creates two new nodes that require boundary data, for which we introduce two parameters.
As our aim is to show that the mixture model~\labelcref{eq: gas flow} admits a solution, we derive additional constraints that a solution for the cut networks must satisfy in order to be also a solution to the original one, and prove that there exist parameters that satisfy these constraints.

Before we prove \cref{thm: existence cycle}, we derive an equivalent characterization of the solvability of the mixture model~\labelcref{eq: gas flow}. At the end of this section, we utilize this characterization to conclude that the mixture model~\labelcref{eq: gas flow} admits at least one solution.
We start by introducing the notion of a \emph{cut graph} to formalize the process of cutting edges of a graph.
         
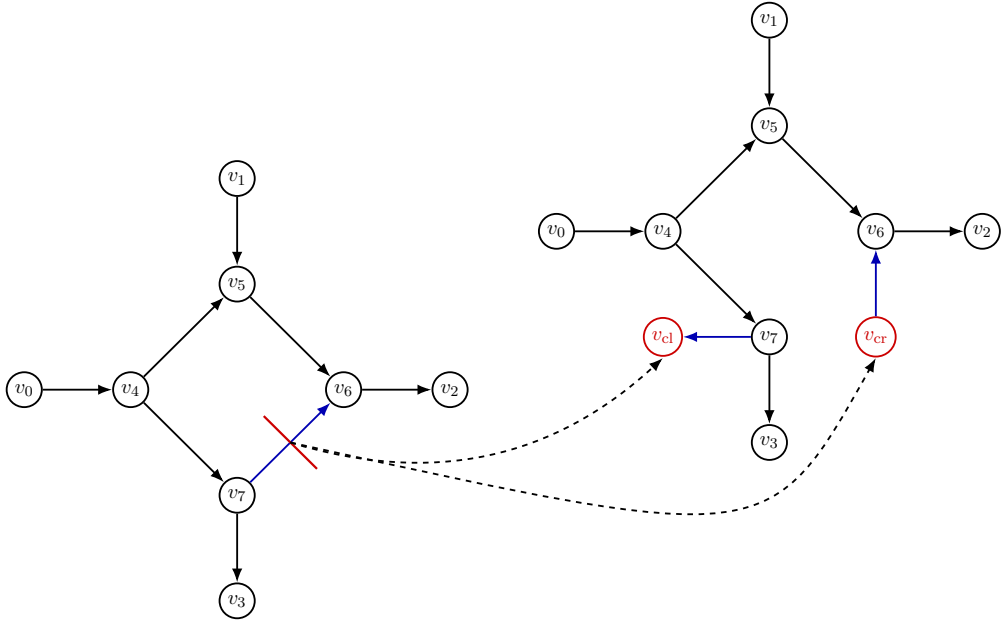
\begin{figure}[ht]
	\centering
	\input{./tikz/cut_graph}
	\caption{A graph with a single circle (left) and the resulting cut graph when cutting the blue edge (right).}
	\label{fig: cut graph}
\end{figure}

\begin{definition}[{Cut Graph, cf. \cite{Gugat2016, WintergerstDiss}}]
	Let \( \G \) be a directed, connected graph and let \( e^c = ( f(e^c), h(e^c)) \) be an edge of \( \G\). Then the so-called \emph{cut graph \( \G^c = (\V^c, \E^c)\)} of \(\G\) with respect to \( e^c\) is defined by
	\begin{align*}
		\V^c = \V \cup \{ v_\cl, v_\crr\},
		\quad
		\E^c = \left(\E \setminus \{ e^c\} \right) \cup \{ e_\cl, e_\crr\},
	\end{align*}
	where the nodes \( v_\cl, v_\crr \notin \V\) are generated by cutting the edge \( e^c\) which splits into the two new edges \( e_\cl = (f(e^c), v_\cl)\) and \( e_\crr = (v_\crr, h(e^c))\), see \cref{fig: cut graph} for a visualization.
\end{definition}

Since our goal is to find solutions to the mixture model on the cut graph \( \G^c\) that are also solutions to the original graph \( \G\), we assume that the boundary data for all old nodes \( v\in \V\) remains the same.
At the new nodes \( v \in \{v_\cl, v_\crr\}\), solutions for the cut graph \( \G^c\) -- indicated by the superscript \(c\) -- must satisfy three conditions, in order to be solutions for the original graph \( \G\) as well:
\begin{multicols}{3}
\begin{enumerate}[label=(c.\roman*)]
	\item \(b_{v_\cl}^c = - b_{v_\crr}^c,\)\label{item: c1}
	\item \(\eta_{v_\cl}^c = \eta_{v_\crr}^c,\) \label{item: c2}
	\item \( p_{v_\cl}^c = p_{v_\crr}^c.\) \label{item: c3}
\end{enumerate}
\end{multicols}
\noindent
The first condition~\labelcref{item: c1} guarantees that the amount of gas leaving node \( {v_\cl}\) is equal to the amount of entering node \( v_\crr\). One can also think of this condition as an invisible pipe connecting the nodes \( v_\cl\) and \( v_\crr\), ensuring that the gas flow through the cut is constant.
The condition~\labelcref{item: c2} and the condition~\labelcref{item: c3} ensure that the composition and the pressure are constant through the cut, respectively.

Next, we introduce appropriate boundary data for the newly introduced nodes \( v_\cl\) and \( v_\crr\)  to guarantee the problem on the cut graph is well-defined and to ensure that the condition~\labelcref{item: c1} is always satisfied.
Therefore, we introduce the parameter \( \lambda \in \R\) and set:
\begin{align*}
	b_{v_\cl}^c = \lambda \quad \text{and} \quad b_{v_\crr}^c  = - \lambda.
\end{align*}
Further, we introduce a second parameter \( \mu \in [0,1]\) to model the composition of the supply at the nodes \( v_\cl\) and \(v_\crr\) depending on the sign of \( \lambda\).
Then, the full set of boundary data for the gas flow model on the cut graph then reads: 
\begin{align*}
	p_{v^{\ast}}^c = p^{\ast},  \quad
	b^c_v = \begin{cases}
		b_v, & v \in \V, \\
		-\lambda, & v = v_\crr, \\
		\lambda, & v = v_\cl,
	\end{cases}  \quad \text{and} \quad 
	\zeta_v^c = \begin{cases}
		\zeta_v, & v \in \V_{<0}, \\
		\mu, & v = v_\crr \text{ and } \lambda \le 0, \\
		\mu, & v = v_\cl \text{ and } \lambda > 0.
	\end{cases}
\end{align*}
We also have to set the length, diameter, and friction factor of the edges \( e \in \{e_\cl, e_\crr\}\) generated by the cut. We set the edge length to half the length of the cut edge \( e^c\), i.e., we have \( L^{e} = \frac{1}{2} L^{e^c}\) for \( e \in \{e_\cl, e_\crr\}\).
Since we assume that all edges have the same diameter and the same friction factor, this property transfers to the two new edges.

\medbreak
We rewrite the condition of constant composition and constant pressure through the cut by using the functions:
\begin{subequations}\label{eq: Hp and  eta}
\begin{align}
	H_p \colon \R \times [0,1] \rightarrow \R, \quad
	H_p(\lambda, \mu) & = p_{v_\crr}(\lambda, \mu)^2 - p_{v_\cl}^c(\lambda, \mu)^2 \\
	H_\eta \colon \R \times [0,1] \rightarrow \R, \quad 
	H_{\eta}(\lambda, \mu) &= \eta_{v_\crr}^c(\lambda, \mu) - \eta_{v_\cl}^c(\lambda, \mu),
\end{align}
\end{subequations}
which measure the difference between the pressure and the composition at the new nodes, respectively. Notice that if the pipes are sufficiently short, cf. \cref{eq: L crit}, the condition \( H_p(\lambda, \mu) = 0\) is equivalent to \( p_{v_\crr}(\lambda, \mu) = p_{v_\cl}(\lambda, \mu) \).
Then, the task of finding parameters \( \lambda\) and \(\mu\), which satisfy the constraints~\labelcref{item: c1},~\labelcref{item: c2}, and~\labelcref{item: c3}, becomes:
\begin{align*}
	\text{Find} \quad (\lambda, \mu) \in \R \times [0,1] \quad \text{such that} \quad 
	H_p (\lambda, \mu) = 0 \quad \text{and} \quad 
	H_{\eta} (\lambda, \mu) = 0,
\end{align*}
which means we have to find a common root of \( H_p\) and \( H_{\eta}\). 
With this derivation at hand, we obtain an equivalent characterization of the solvability of the steady-state gas flow model~\labelcref{eq: gas flow}, which we capture in the following lemma.

\begin{lemma}\label{lemma: solvability 2}
	Let \( \G = (\V, \E)\) be a directed graph. Then the following holds:
	\begin{center}
		\Centerstack[c]{The gas flow model~\labelcref{eq: gas flow} \\ has a solution.}
		\( \quad \Leftrightarrow \quad \) 
		\Centerstack[c]{The functions \( H_p\) and \( H_{\eta}\) \\ have a common root.}
	\end{center}
\end{lemma}
\begin{proof}
	Assume the steady-state gas flow model~\labelcref{eq: gas flow} has a solution \( q_e, \eta_v, p_v\) and let \( e^c\) be the cut edge. Then, by definition of \( H_p\) and \( H_{\eta}\), these two functions have the common root:
	\begin{align*}
		\lambda^{\ast} = q_{e^c} \quad \text{and} \quad 
		\mu^{\ast} = \begin{cases}
			\eta_{f(e^c)}, & \text{if } q_{e^c} \ge 0, \\
			\eta_{h(e^c)}, & \text{if } q_{e^c} < 0.
		\end{cases}
	\end{align*}
	Now suppose that \( H_p\) and \( H_{\eta}\) have a common root \( (\lambda^{\ast}, \mu^{\ast})\). Then, we recover a solution to the gas flow model~\labelcref{eq: gas flow} by setting:
	\begin{align*}
		q_e & = \begin{cases}
			q_e^c(\lambda^{\ast}, \mu^{\ast}), & e \in \E \setminus \{ e^c\}, \\
			\lambda^{\ast}, & e = e^c, 
		\end{cases}
		\quad 
		\eta_v =\eta_v^c(\lambda^{\ast}, \mu^{\ast}),
		\quad 
		p_v = p_v^c(\lambda^{\ast}, \mu^{\ast}). \qedhere
	\end{align*}
\end{proof}

In case the functions \( H_p\) and \( H_{\eta}\) have a common root, the solvability of the mixture model~\labelcref{eq: gas flow}, i.e., \cref{thm: existence cycle}, follows immediately with \cref{lemma: solvability 2}.
Therefore, our next goal is to show that the functions \( H_p\) and \( H_{\eta}\) have a common root.

\begin{lemma}\label{lemma: common root}
	Let \( \G= (\V,\E)\) be a directed graph with one cycle and let \( e^c \in \E\) be an edge belonging to the cycle. Further, let \( \G^c\) be the cut graph obtained by cutting the edge \(e^c\). 
	Then, there exists a tuple \( (\lambda^{\ast}, \mu^{\ast}) \in \R \times [0,1]\) such that
	\begin{align*}
		H_p(\lambda^{\ast}, \mu^{\ast})  = 0 \quad \text{and} \quad 
		H_{\eta}(\lambda^{\ast}, \mu^{\ast})  = 0,
	\end{align*}
	where the functions \( H_p\) and \( H_{\eta}\) are defined in \cref{eq: Hp and  eta}, i.e., these two functions have at least one common root.
\end{lemma}
%
\noindent
The idea of the  proof of \cref{lemma: common root} is as follows:
\begin{enumerate}[label=(\roman*)]
	\item For a fixed \( \lambda\), we show that the function \( H_{\eta}\) admits a unique root \( \mu_{\eta}(\lambda) \), which means \( H_{\eta}(\lambda, \mu_{\eta}(\lambda)) = 0\) holds. \label{item: step 1}
	\item Then, we restrict the function \( H_p\) to the root curve \( \mu_{\eta}(\lambda) \) resulting in the scalar function \( g(\lambda) = H_p(\lambda, \mu_{\eta}(\lambda))\). Then, we apply the  intermediate value theorem to the function \(g\) to show that it admits at least one root, i.e., we have to show that 
	\begin{enumerate}[label=(\alph*)]
		\item the function \( g \) is continuous, and that
		\item there exit values \( \lambda^-\) and \(\lambda^+ \) such that \(g(\lambda^-) \le 0\) and \(g(\lambda^+) \ge 0.\)
	\end{enumerate} \label{item: step 2}
\end{enumerate}
In \cref{fig: overview}, we provide an overview of the intermediate results required and group them into two categories corresponding to the two steps of the proof. In the following, we postpone the detailed proof to show the intermediate results first.

\begin{figure}[ht]
	\centering
	\input{./tikz/overview}
	\caption{Overview and dependence of intermediate results to show \cref{lemma: common root}.}
	\label{fig: overview}
\end{figure}


To show the two intermediate steps of the proof, we analyze the properties of the nodal composition \( \eta_v^c\) and the flow \( q_e^c\) on the edges, because both functions, \( H_p\) and \(H_{\eta}\), can be written in terms of \(\eta_v^c\) and \( q_e^c\).
For the function \(H_{\eta}\), this holds by definition. 
The function \( H_p\) describes the pressure loss between the node \(v_\cl\) and the node \( v_\crr\). Hence, we invoke \cref{lemma: pressure loss path} to express \( H_p\) in terms of  \(\eta_v^c\) and \( q_e^c\), which results in:
\begin{align}\label{eq: H_p alternative}
	H_p(\lambda, \mu) = \sum_{e \in \E_\P} \tilde{\sigma}( \eta_{f(e)}^c(\lambda, \mu) , \eta_{h(e)}^c(\lambda, \mu) , q_e^c(\lambda)) a^c(v_e,e) q_e^c (\lambda) \vert a^c(v_e, e) q_e^c (\lambda) \vert,
\end{align}
where \( \E_\P \subseteq \E^c\) is the set of edges that belong to the path connecting \(v_\cl\) and \( v_\crr\) and \( v_e\) denotes the node of the edge \( e\) that is closer to the start node of the path.

Thus, besides analyzing the properties of the flow \( q_e^c\) and the nodal composition \( \eta_v^c\), we also show how these properties carry over to the functions \( H_p\) and \( H_{\eta}\), and how they contribute to the proof of \cref{lemma: common root}.

\subsubsection{Properties of the Flow and the Nodal Composition}\label{sec: properties flow}
We start by showing that the flow \( q_e^c\) on the edges of the former cycle of the network depends linearly on \( \lambda\) as this property is essential to that the function \( H_{\eta}\) admits a unique root, and that the function \(  \lambda \mapsto H_p(\lambda, \mu_{\eta}(\lambda))\) satisfies the requirements on the intermediate value theorem.

\begin{lemma}\label{lemma: flow linear}
	Let \( \G = (\V, \E)\) be a graph with one cycle and let \( e^c \in \E \) be the cut edge belonging to the cycle \( (\C, \E_\C) \). Assume \( \lambda\) and \( \mu\) are fixed.
	Then the flow \( q_e^c\) on the cut graph \( \G^c \) depends linearly on \( \lambda\) and can be written as:
	\begin{equation*}
		a^c(v_e, e) q_e^c(\lambda) = \gamma_e - \lambda
		\quad \text{for all} \quad e \in \left( \E_\C \setminus \{e^c\} \right) \cup \{e_\cl, e_\crr\},
	\end{equation*}
	where the constants \( \gamma_e \) depend only on the loads \( b_v\) and \( v_e\) is the node of the edge \(e\) that is closer to the node \( v_\crr\) (cf. \cref{lemma: pressure loss path} for the definition of \( v_e\)).
\end{lemma}
\begin{proof}    
	Since the graph \( \G\) has only one cycle that is cut, the cut graph \( \G^c\) is tree-shaped. Therefore, there exits exactly one path connecting the nodes \( v_\cl\) and \( v_\crr\). We denote the path as a connected sub-graph \( (\P, \E_\P)\) of \( \G^c\), where
	\begin{align}\label{eq: path as graph}
		\P = \C \cup \{ v_\cl, v_\crr\} \quad \text{and} \quad 
		\E_\P = \E_\C \setminus \{e^c \} \cup \{e_\cl, e_\crr\}.
	\end{align}
	The proof follows part of the proof of \cite[Theorem 1]{Gugat2016}. Without loss of generality we assume that the path starts at node \( v_\crr\) and ends at node \( v_\cl\).  
	Every edge \( e \notin \E_\P\) is part of a sub-tree of \( \G^c\) with one node in \( \P\). As the loads \( b_v^c\) for nodes \( v \notin \P\) are fixed and independent of \( \lambda\), the flow on the edge \( e \notin \E_\P\) is uniquely determined and independent of \(\lambda\). Hence, only the flow values along the path \( (\P, \E_\P)\) can change when \( \lambda\) changes and we can redefine the loads for the nodes on the path:
	\begin{equation*}
		b_v^\P \coloneqq b^c_v - \sum_{e \in \E(v) \setminus \E_\P} a^c(v,e) q_e^c \quad \text{for} \quad v \in \P.
	\end{equation*}
	Then, the conservation of mass, cf. \cref{eq: mass conservation}, becomes \( A^\P \q^\P = \bb^\P\) where \( A^\P\) is the incidence matrix of the sub-graph of the path and has a lower triangular structure:
	\begin{align*}
		A^\P = \begin{pmatrix}
			a^c(v_1, e_1) & & & & \\
			a^c(v_2, e_1) & a^c(v_2, e_2) & & & \\
			& a^c(v_3, e_2) & \ddots & & \\
			& &  \ddots & a^c(v_{p-2}, e_{p-2}) & \\
			& & &  a^c(v_{p-1}, e_{p-2}) & a^c(v_{p-1}, e_{p-1}) \\
			& & & & a^c(v_p, e_{p-1})
		\end{pmatrix},
	\end{align*}
	where \( p = \vert \P \vert = \vert \E_\P \vert + 1\). The numbering of the nodes indicates the order in which they are traversed by the path, i.e., \( v_1 = v_\crr\) and \( v_p = v_\cl\). An edge \( e_i\) connects the nodes \( v_i\) and \( v_{i+1}\). The vector \( \q^\P\) has entries corresponding to the flow values of the edges \( e \in \E_P\) and the vector \( \bb^\P\) entries corresponding to the redefined loads \( b_v^\P\). 
	Then, we apply forward substitution to the linear system and obtain:
	\begin{align*}
		a^c(v_i, e_i) q^c_{e_i} & = \sum_{k=1}^i b_{v_k}^\P \quad \text{for} \quad i = 1, \ldots, p - 1.
	\end{align*}
	As the path starts at node \( v_1 = v_\crr\), we have \( a^c(v_1,e_1) q^c_{e_1} = b_{v_1}^\P = - \lambda \). Further, the redefined loads \( b_v^\P\) are independent of \( \lambda\) for \( v \notin\{ v_1 = v_\cl, v_p = v_\crr \}\), which leads to
	\begin{align*}
		a^c(v_i, e_i) q^c_{e_i} & = - \lambda + \underbrace{\sum_{k=2}^i \left[b_{v_k} - \sum_{e \in \E(v_k) \setminus \E_\P} a(v_k,e) q_e^c \right]}_{\eqqcolon \gamma_{e_i} \text{ which is independent of } \lambda} \quad \text{for} \quad i = 1, \ldots p - 1\\
		\Leftrightarrow \quad a^c(v_e, e) q_e^c & = \gamma_e - \lambda \quad \text{for} \quad e \in \E_\P = \left( \E_\C \setminus \{e^c\} \right) \cup \{e_\cl, e_\crr\}. \qedhere
	\end{align*}
\end{proof}

With this representation of \( q_e^c\) for edges of the former cycle and the knowledge that the flow is independent of \( \lambda\) for edges not belonging to the former cycle, it follows directly that the flow is continuous in \( \lambda\).

\begin{korollar}[Continuity of the Flow]\label{korollar: flow continuous}
	Let the requirements of \cref{lemma: flow linear} be satisfied. Then the flow \( q_e^c\) on the cut graph \( \G^c\) is continuous in \( \lambda\) for all edges \( e \in \E^c\).
\end{korollar}
\begin{proof}
	The claim follows directly from the proof of \cref{lemma: flow linear} as \( q_e^c\) is independent of \( \lambda\) if \( e \notin \E_\C \cup \{e_\cl, e_\crr \}\) and linear in \( \lambda\) if \( e \in \E_\C \cup \{e_\cl, e_\crr \}\).
\end{proof}

The continuity of \( q_e^c\) with respect to \( \lambda\) allows us to show that the nodal compositions \( \eta_v^c\) are also continuous. 
However, the roles of the nodes \( v_\cl\) and \( v_\crr\) depend on the sign of \( \lambda\).
That is, for \( \lambda > 0\), the node \( v_\cl\) is a demand node and \( v_\crr\) a supply node, while for \(\lambda<0\), the roles are reversed.
One can compare the switching of the sign to a supply station that suddenly becomes a consumer, in which case we anticipate the change, and thus, the composition at these nodes to be discontinuous for \( \lambda = 0\). 

\begin{lemma}[Continuity of the Composition  with respect to \(\lambda\)]\label{lemma: composition continuous}
	Let \( \G\) be a connected graph with a single cycle, \( e^c \in \E \) an edge belonging to the cycle, \( \G^c\) the resulting cut graph. Assume \( \mu \in [0,1]\) is arbitrary, yet fixed.
	Then:
	\begin{enumerate}[label=(\roman*)]
		\item The nodal compositions \( \eta_v^c\) are continuous for \( \lambda \in \R \setminus \{0\}\).
		\item The nodal compositions \( \eta_v^c\) are also continuous in \( \lambda = 0\) for nodes \( v \in \V \setminus \{ v_\cl, v_\crr \}\).
	\end{enumerate}
\end{lemma}
\begin{proof}
	We prove the claim by using the concept of a topological ordering to traverse through the nodes of \(\G^c \) as the graph \( \G_\flow^c \) with then same nodes as cut graph \( \G^c\) but edges according to the flow direction, is tree-shaped and thus, has a topological ordering of its nodes. 
	
	The flow along the path connecting the nodes \( v_\cl\) and \( v_\crr\) only changes when \( \lambda\) changes. because the cut graph \( \G^c\) is tree-shaped.
	\cref{lemma: flow linear} allows to determine flow direction of an edge belonging to this path depending on \( \lambda\) as well as for which value of \( \lambda\) the flow of an edge becomes zero. We sort these roots in ascending order:
	\begin{equation*}
		\gamma_{\min} = \tilde{\gamma}_1 < \cdots < \tilde{\gamma}_{\vert \E_\P\vert} = \gamma_{\max}
		\quad \text{where} \quad \{ \gamma_e \mid e \in \E_\P\} = \{ \tilde{\gamma}_i \mid i = 1, \ldots, \vert \E_\P\vert\}.
	\end{equation*}
	Then, for a fixed \( i\), the flow directions are the same for all \( \lambda \in (\tilde{\gamma}_i, \tilde{\gamma}_{i+1})\). Thus, the topological ordering differs for each sub-interval \(  (\tilde{\gamma}_i, \tilde{\gamma}_{i+1})\) as the flow directions depend on \( \lambda\).
	Now, we show the continuity of \( \eta_v^c \) with respect to \( \lambda\) on \( (\tilde{\gamma}_i, \tilde{\gamma}_{i+1})\) by an induction over the position \( k\) in the topological ordering for a fixed, yet arbitrary index  \( i \in \{ 1, \ldots, \vert \E_\P\vert\}\).
	
	The first node \( v\) in the ordering has to be a supply node where the composition is given through boundary data, i.e., \( \eta_v^c = \zeta_v^c\) is constant and hence continuous in \( \lambda\).
	
	In the induction step, we assume that all compositions up to position \( k \) are continuous on the sub-interval and that \( v\) is the node at position \( k+1\). Then the composition \( \eta_v^c\) is given by:
	\begin{equation*}
		\eta_v^c = \frac{\sum_{e \in \E(v)} \eta_{w_e}^c\left( a^c(v,e) q_e^c \right)^+ + \zeta_v^c (b_v^c)^-}{\sum_{e \in \E(v)} \left( a^c(v,e) q_e^c \right)^+ + (b_v^c)^-}
		\quad \text{where} \quad w_e = \begin{cases}
			f(e), & \text{if } v = h(e), \\
			h(e), & \text{if } v = f(e).
		\end{cases}
	\end{equation*}
	Because the topological ordering is based on the graph with edge in the flow direction, the inflow comes from nodes previous in the ordering. Thus, by the induction assumption we know that the compositions \( \eta_{w_e}^c\) and the flows \(q_e^c\) are continuous on the sub-interval and since \( a^c(v,e)q_e^c \neq 0\) by definition of the sub-interval, we also never divide by zero. Hence, the compositions \( \eta_v^c\) are continuous for \( \lambda \neq \tilde{\gamma}_i\). The continuity for the intervals \( (-\infty, \tilde{\gamma}_0)\) and \((\tilde{\gamma}_{\vert \E_\P\vert}, \infty)\) follows analogously.
	
	It remains to prove the continuity in \( \lambda = \tilde{\gamma}_i\). 
	Let \( \tilde{e} \in \E\) be an edge with \( q_{\tilde{e}}(\tilde{\gamma}_i) = 0\). Then, the continuity of \( \eta_v^c\) in \( \lambda = \tilde{\gamma}_i\) for nodes \( v\) not incident to \( e\) follows directly. 
	For nodes incident to \( \tilde{e}\) two cases occur:
	\begin{enumerate}[label=(\roman*)]
		\item \( a^c(v, \tilde{e}) q_{\tilde{e}}^c < 0 \) for \( \lambda < \tilde{\gamma}_i\) \(\quad \Rightarrow \quad \)
		\( a^c(v, \tilde{e}) q_{\tilde{e}}^c > 0 \) for \( \lambda > \tilde{\gamma}_i\), and
		\item \( a^c(v, \tilde{e}) q_{\tilde{e}}^c > 0 \) for \( \lambda < \tilde{\gamma}_i\) \(\quad \Rightarrow \quad \)
		\( a^c(v, \tilde{e}) q_{\tilde{e}}^c < 0 \) for \( \lambda > \tilde{\gamma}_i\).
	\end{enumerate}
	We demonstrate how to handle the first case since the second one is analogous. Let \( v\) be a node incident to \( \tilde{e}\), then the following holds:
	\begin{align*}
		\lim_{\lambda \rightarrow \tilde{\gamma}_i-}  \eta_v^c(\lambda, \mu)  
		& = \lim_{\lambda \rightarrow \tilde{\gamma}_i-} 
		\frac{\sum_{e \in \E(v) \setminus \{ \tilde{e} \}} \eta_{w_e}^c\left( a^c(v,e) q_e^c \right)^+ + \zeta_v^c (b_v^c)^-}{\sum_{e \in \E(v) \setminus \{ \tilde{e} \}} \left( a^c(v,e) q_e^c \right)^+ (b_v^c)^-} \\
		&= 	\frac{\sum_{e \in \E(v) \setminus \{ \tilde{e} \}} \eta_{w_e}^c\left( a^c(v,e) q_e^c \right)^+ + \zeta_v^c (b_v^c)^-}{\sum_{e \in \E(v) \setminus \{ \tilde{e} \}} \left( a^c(v,e) q_e^c \right)^+ + (b_v^c)^-} \\
		& = \lim_{\lambda \rightarrow \tilde{\gamma}_i+} 
		\frac{\sum_{e \in \E(v)} \eta_{w_e}^c\left( a^c(v,e) q_e^c \right)^+ + \zeta_v^c(b_v^c)^- }{\sum_{e \in \E(v)} \left( a^c(v,e) q_e^c \right)^+ +(b_v^c)^-} 
		= \lim_{\lambda \rightarrow \tilde{\gamma}_i+} \eta_v^c(\lambda,\mu),
	\end{align*}
	   which proves the continuity of \( \eta_v^c\) for \( v \in \V \setminus \{ v_\cl, v_\crr \}\). Notice, that the nodes \(v_\crr\) and \( v_\cl\) lead to a special case since both nodes are boundary nodes and their load \( b_v^c\) can change depending on the value of \( \lambda\). In particular, both nodes can switch between supply and demand node when  \( \lambda\) switches its sign. Hence, the compositions \(\eta_v^c\) are discontinuous in \( \lambda = 0\) in general for \( v \in\{ v_\cl, v_\crr \}\).
\end{proof}

The continuity of the composition with respect to \( \mu\) follows from the fact that each nodal composition \( \eta_v\) is either independent of \( \mu\) or depends linearly on it.

\begin{lemma}[Continuity of the Composition  with respect to \(\mu\)]\label{lemma: composition linear}
	Let \( \G\) be a connected graph with a single cycle, \( e^c \in \E \) an edge belonging to the cycle, \( \G^c\) the resulting cut graph. Moreover, let \( u \in \V_{<0}^c \) be a supply node with supply composition \( \zeta_u^c\) and let \( q_e^c\) be a flow on the network satisfying \cref{eq: mass conservation}. Then for every node \( v \in \V^c\), one of the two cases apply:
	\begin{enumerate}[label=(\roman*)]
		\item The nodal composition \( \eta_v^c\) depends linearly on \(\zeta_u^c\).
		\item The nodal composition \( \eta_v^c\) is independent of \(\zeta_u^c\).
	\end{enumerate}
	Further, the nodal compositions are well-defined, meaning \( \eta_v^c \in [0,1]\) for all \( v \in \V^c\), if the boundary data satisfies \( \zeta_v^c \in [0,1]\) for every supply node \(v \in \V_{<0}^c\).
\end{lemma}
\begin{proof}
	There exists a topological ordering of the nodes of \( \G_\flow^c \) since we cut the only cycle of \( \G\) and we can prove the claim by induction over the position \( k\) in the ordering. 
	
	The node at position \( k = 1\) must be a supply node where \( b_v^c < 0\) holds. Therefore, \( \eta_v^c\) is given directly by boundary data, which implies that \( \eta_v\) is independent of  \(\zeta_{u}^c\) if \( v \neq u\) and linear in  \(\zeta_{u}^c\) if \( v = u\). By choice of the boundary data, we have \( \eta_v^c \in [0,1]\).
	
	Now, assume that the claim holds for all nodes up to position \( k\) in the ordering. Let \( v\) be the node at position \( k+1\). Then, two cases can apply: (a) \(v\) is a supply node, or (b) \(v\) is not. The case (a) is analogous to the case \( k=1\). In case (b), the composition can be computed by 
	\begin{equation*}
		\eta_v^c = \frac{\sum_{e \in \E(v)} \eta_{w_e}^c\left(a^c(v,e) q^c_e \right)^+ +\zeta_v^c (b_v^c)^- }{\sum_{e \in \E(v)} \left( a^c(v,e) q_e^c \right)^+ + (b_v^c)^-}
		\quad \text{where} \quad w_e = \begin{cases}
			f(e), & \text{if } v = h(e), \\
			h(e), & \text{if } v = f(e).
		\end{cases}
	\end{equation*}
	As the gas flows from nodes early in the ordering to nodes later in the ordering, we know by the induction assumption that the compositions \( \eta_{w_e}^c\) are linear or independent of  \(\zeta_w^c\) and \( \eta_{w_e}^c \in [0,1]\). 
\end{proof}

\subsubsection{Properties of the Functions \( H_p\) and \(H_{\eta}\)}\label{sec: properties Hp}
The next step is to extend the results that hold for the flow \( q_e^c\) and composition \( \eta_v^c\) to the functions \( H_p\) and \( H_{\eta}\) since both functions depend on these variables, cf. \cref{eq: Hp and eta} and \cref{eq: H_p alternative}.
First, we show that the function \( H_p\) is continuous in \( \lambda\) and \( \mu\) as the flow and the composition are continuous too.

\begin{lemma}[Properties of \( H_p\)]\label{lemma: Hp continuous}
	The function \( H_p \colon \R \times [0,1] \rightarrow \R\) is 
	\begin{enumerate}[label=(\roman*)]
		\item linear in \( \mu\), i.e., continuous and monotone in \( \mu\), 
		\item continuous in \( \lambda\), and
		\item constant in \( \mu\) if \( \lambda = 0\).
	\end{enumerate} 
\end{lemma}
\begin{proof}
	The function \( H_p\) describes the pressure loss along the path connecting the nodes \( v_\cl\) and \( v_\crr\). Hence, we invoke \cref{lemma: pressure loss path} to obtain an alternative expression of \( H_p\):
	\begin{align*}\label{eq: H_p alternative}
		H_p(\lambda, \mu) = \sum_{e \in \E_\P} \tilde{\sigma}( \eta_{f(e)}^c(\lambda, \mu) , \eta_{h(e)}^c(\lambda, \mu) , q_e^c(\lambda)) a^c(v_e,e) q_e^c (\lambda) \vert a^c(v_e, e) q_e^c (\lambda) \vert,
	\end{align*}
	where \( v_e\) is the node of the edge \( e\) that is closer to node \( v_\cl\).
	Since the cut graph \( \G^c\) is tree-shaped, the flow \( q_e^c\) is fully determined by the loads \( b_v^c\). Thus, the flow \( q_e^c\) along the edges is only dependent on the parameter \(\lambda\) and independent of the parameter \(\mu\).

	\begin{enumerate}[label=(\roman*)]
	\item The linearity of \( H_p\) with respect to \( \mu\) follows directly from the linearity of the composition \( \eta_v\) with respect to \( \mu \), cf. \cref{lemma: composition linear}.
	
	\item The argument to show the continuity with respect to \( \lambda\) is similarly. The continuity in \( \lambda \neq 0\) follows from the continuity of the flow (\cref{korollar: flow continuous}) and the continuity of the composition (\cref{lemma: composition continuous}).
	In order to prove also the continuity in \( \lambda = 0\), we rewrite the function \( H_p\) as
	\begin{multline*}
		H_p(\lambda, \mu) = \underbrace{\sum_{\substack{e \in \E_\P \\ e \neq e_\cl, e_\crr}} \tilde{\sigma}( \eta_{f(e)}^c(\lambda, \mu), \eta_{h(e)}^c(\lambda, \mu), q_e^c (\lambda)) a^c(v_e,e) q_e^c (\lambda) \vert a^c(v_e, e) q_e^c (\lambda) \vert}_{\text{continuous in } \lambda = 0 \text{ by \cref{korollar: flow continuous} and  \cref{lemma: composition continuous}}} \\  
		+\underbrace{\tilde{\sigma}(\eta_{v_\crr}^c(\lambda), \eta_{h(e_\crr)}^c(\lambda, \mu), \lambda) \lambda \vert  \lambda \vert 
			+ \tilde{\sigma}(\eta_{f(e_\cl)}^c(\lambda, \mu) , \eta_{v_\cl}^c (\lambda), \lambda) \lambda \vert \lambda \vert}_{\eta_{v_\cl}^c \text{ and } \eta_{v_\crr}^c \text{ are in general discontinuous in } \lambda = 0} .
	\end{multline*}
	Although, the compositions \( \eta_{v_\cl}^c\) and \( \eta_{v_\crr}^c\) are discontinuous at \( \lambda = 0\) by \cref{lemma: composition continuous}, the function \( H_p\) is continuous at \( \lambda = 0\) because terms of the form \( \tilde{\sigma}(\eta_{f(e)}^c, \eta_{h(e)}^c, \lambda) \lambda \vert \lambda \vert\) are continuous at \( \lambda = 0\) even if the compositions \( \eta_{f(e)}^c\) and \( \eta_{h(e)}^c\) are not. This is because the function \( \tilde{\sigma}\) is bounded from above which leads to
	\begin{align*}
		\bigg \vert  \tilde{\sigma}(\eta_{f(e)}^c, \eta_{h(e)}^c, \lambda) \lambda \vert \lambda \vert  \bigg\vert 
		\le \frac{\lambda_\fr L^e}{D} \max\{\sigma_\hh^2, \sigma_\ng^2\} \lambda^2 \longrightarrow 0
		\quad \text{for} \quad \lambda \longrightarrow 0.
	\end{align*}
	Hence, the discontinuity vanishes for \( \lambda = 0\) and \( H_p\) is continuous at every \( \lambda \).
	
	\item At last,  we prove that \( H_p\) is constant in \( \mu \)  if \( \lambda = 0\). 
	As the flow is fully determined by the loads, the dependence on the parameter \( \mu\) enters through the mixing equation. Specifically, through the term \( \mu \lambda^2\) when computing the nodal composition \( \eta_v^c\). Hence, if \( \lambda = 0\) the nodal compositions are independent of \( \mu\) and we get \( 	H_p(0, \mu)  \equiv \const. \) \qedhere
	 \end{enumerate}
\end{proof}

Next, we identify for which values of \( \lambda\) the flow is either directed from node \( v_\cl\) to node \( v_\crr\) or vice versa. As the function \( H_p\) describes the pressure loss along the path connecting \( v_\cl\) and \( v_\crr\), this allows us to determine the sign of \( H_p\) for these values of \( \lambda\).

\begin{lemma}[Unidirectional Flow]\label{lemma: unidirectional flow}
	Assume the requirements of \cref{lemma: flow linear} are satisfied and define \(\gamma_{\min} = \min_{e \in \C \setminus \{ e^c\}} \gamma_e \) and \(\gamma_{\max} = \max_{e \in \C \setminus \{ e^c\}} \gamma_e\).
	Then:
	\begin{enumerate}[label=(\roman*)]
		\item The gas flows from \( v_\crr\) to \( v_\cl\) if and only if \( \lambda \ge \gamma_{\max}\),  which implies \(H_p(\lambda, \mu) \ge 0\).
		\item The gas flows from \( v_\cl\) to \( v_\crr\) if and only if \( \lambda \le \gamma_{\min} \), which implies \(H_p(\lambda, \mu) \le 0\).
	\end{enumerate}  
\end{lemma}
\begin{proof} 
	Since the requirements of \cref{lemma: flow linear} are satisfied, the cut graph \( \G^c\) is tree-shaped and thus there exits only one path connecting the nodes \( v_\cl\) and \( v_\crr\). Consider the sub-graph of \( \G^c \) that forms the path between \( v_\crr \) and \( v_\cl \) which we denote by 
	\begin{align*}
		(\P, \E_\P) \quad \text{where} \quad \P & = \C  \cup \{ v_\cl, v_\crr\} 
		\quad \text{and} \quad 
		\E_{\P}  = \left( \E_{\C} \setminus \{e^c\} \right) \cup \{e_\cl, e_\crr\}.
	\end{align*}
	We assume without loss of generality that the path starts at \( v_\crr\) and ends at \( v_\cl\).
	Then the flow being directed from \( v_\crr\) to \(v_\cl\) translates into  \(  a^c(v_e, e) q_e^c \le 0\) and being directed in the opposite direction translates into \(a^c(v_e, e) q_e^c \ge 0\) for all \( e \in \E_{\P}\).
	Consider the first case. Then, by \cref{lemma: flow linear} we get
	\begin{equation*}
		a^c(v_e, e) q_e^c = \gamma_e - \lambda \le \gamma_e - \gamma_{\max} \le 0
		\quad \text{for all } e \in \E_\P \quad \Leftrightarrow \quad \lambda \ge \lambda_{\max}
	\end{equation*}
	which means that the gas flows from node \( v_\crr\) to node \( v_\cl\). The node \( v_e\) is the node of the edge \( e\) which is closer to node \( v_\crr\).
	Now, we apply \cref{lemma: pressure loss path} which yields:
	\begin{equation*}
		H_p(\lambda, \mu) = p_{v_\crr}^c(\lambda,\mu)^2 - p_{v_\cl}^c(\lambda, \mu)^2 \ge 0.
	\end{equation*}
	The second case follows analogously.
\end{proof}

From \cref{lemma: composition linear}, we know that the composition \( \eta_v^c\) depends linearly on \( \mu\). As the function \( H_{\eta}\) is the difference between two nodal compositions, it is clear that it also depends linearly on \( \mu\). We exploit this linearity in the next lemma to prove that \( H_{\eta}\) has a unique root.

\begin{lemma}[Root Curves of  \( H_\eta\)]\label{lemma: root curve H eta}
	The function \( H_{\eta} \colon \R \times [0,1]\) with 
	\begin{equation*}
		H_{\eta}(\lambda, \mu) = \eta_{v_\crr}^c(\lambda, \mu) - \eta_{v_\cl}^c(\lambda, \mu),
	\end{equation*}
	is linear in \( \mu\).	Furthermore, for a fixed \( \lambda\), the function \( H_{\eta}\) admits a unique root \( \mu_{\eta}(\lambda)\). For \( \lambda \in \left[\gamma_{\min}, \gamma_{\max}\right]\) the root is given by:
	\begin{equation*}
		\mu_{\eta}(\lambda) = \begin{cases}
			\eta_{v_\crr}^c(\lambda), & \lambda < 0, \\
			\eta_{v_\cl}^c(\lambda), & \lambda \ge 0.
		\end{cases}
	\end{equation*}
\end{lemma}
\begin{proof}
	The linearity of \( H_{\eta}\) for a fixed \( \lambda\) follows immediately from \cref{lemma: composition linear}. Due to the linearity it is also clear that \( H_{\eta}(\lambda, \mu) = 0\) has a unique solution for a fixed \( \lambda\). 
	The goal now is to determine the explicit solution to this equation for \( \lambda \in \left[\gamma_{\min}, \gamma_{\max}\right]\). We distinguish between four cases:
	\begin{multicols}{4}
		\begin{enumerate}[label=(\roman*)]
			\item \( \lambda = \gamma_{\min}\),
			\item \( \lambda \in (\gamma_{\min}, 0)\),
			\item \( \lambda \in [0, \gamma_{\max})\),
			\item \( \lambda = \gamma_{\max}\).
		\end{enumerate}
	\end{multicols}
	
	First, we consider the case \( \lambda \in [0, \gamma_{\max})\). In this case, we know that \( v_\crr\) is a supply node with composition \( \eta_{v_\crr}^c = \mu\) and that \( v_\cl\) is a demand node with composition \( \eta_{v_\cl}^c(\lambda, \mu)\) due to the sign of \( \lambda\).   
	Then, the composition \( \eta_{v_\cl}^c\) is independent of \( \mu\).
	
	Assume otherwise, then the flow along the path connecting \( v_\crr\) and \( v_\cl\) must be directed from \( v_\crr\) to \( v_\cl\), which contradicts the assumption \( \lambda \le \gamma_{\max}\) due to \cref{lemma: unidirectional flow}. 
	Thus, we determine the root of \( H_{\eta}\) for a fixed \( \lambda\) by
	\begin{align*}
		H_{\eta}(\lambda, \mu) = 0 \quad \Leftrightarrow \quad 
		\mu - \eta_{v_\cl}^c(\lambda) = 0 \quad \Leftrightarrow \quad 
		\mu_{\eta}(\lambda) = \eta_{v_\cl}^c(\lambda).
	\end{align*}
	The case \( \lambda \in (\gamma_{\min}, 0)\) follows analogously.
	
	Next, we consider the case \( \lambda = \gamma_{\max}\). Then the flow of at least one of the edges of the path connecting \( v_\crr \) and \( v_\cl\), is zero.
	Furthermore, \( a(v_e,e)q_e^c \le 0\) holds for all edges \( e\) of the path. 
	In this scenario, no flow form node \( v_\crr\) can reach node \( v_\cl\) since at least one edge along the path has flow zero, even though -- technically -- the flow is directed from \( v_\crr\) to \( v_\cl\). 
	Thus, the composition \( \eta_{v_\cl}^c\) is again independent of \( \mu\) and we can compute the root as for the case  \( \lambda \in [0, \gamma_{\max})\). 
	Again, the case \( \lambda = \gamma_{\min}\) follows analogously.
\end{proof}

The structure of the root curve \( \mu_{\eta} \) of \( H_{\eta}\) for \( \lambda \in \left[\gamma_{\min}, \gamma_{\max}\right]\) allows us directly to show that it is continuous since the root curve is given by the nodal composition either of the node \( v_\cl \) or the node \( v_\crr\) depending on the sign of \( \lambda\).

\begin{lemma}[Continuity of the Root Curve \( \mu_{\eta}\)]\label{lemma: root curve continuous}
	The root curve \( \mu_{\eta} \colon \left[ \gamma_{\min}, \gamma_{\max} \right] \rightarrow [0,1] \) defined in \cref{lemma: root curve H eta} is continuous in \( \lambda \in \left[ \gamma_{\min}, \gamma_{\max} \right] \setminus 0\).
\end{lemma}
\begin{proof}
	The claim follows directly from \cref{lemma: composition continuous}.
\end{proof}

Finally, we are equipped to prove that the two functions \( H_p \) and \( H_{\eta}\) have a common root (\cref{lemma: common root}). We recall the structure of the idea of the proof and indicate where the results from \cref{sec: properties flow} and \cref{sec: properties Hp} come into play.
\begin{enumerate}[label=(\roman*)]
	\item For a fixed \( \lambda\), \cref{lemma: root curve H eta} shows that the function \( H_{\eta}\) has a unique root \( \mu_{\eta}(\lambda) \).
	\item Then, we restrict the function \(H_p\) to the root curve \( \mu_{\eta}(\lambda) \) and apply the intermediate value theorem to the restricted function \( g(\lambda) = H_p(\lambda, \mu_{\eta}(\lambda))\), which requires that
	\begin{enumerate}[label=(\alph*)]
		\item the function \( g \) is continuous (\cref{lemma: Hp continuous,lemma: root curve continuous}), and that
		\item there exit values \( \lambda^-\) and \(\lambda^+ \) such that \(g(\lambda^-) \le 0\) and \(g(\lambda^+) \ge 0\) (\cref{lemma: unidirectional flow}).
	\end{enumerate} 
\end{enumerate}
In the following, we provide the detailed proof of \cref{lemma: common root}.

\begin{proof}[Proof of \cref{lemma: common root}]
	We distinguish between to cases: The case where \( \lambda = 0 \) lies on the boundary of \( [\gamma_{\min}, \gamma_{\max}]\) and the case where \( \lambda = 0 \) lies on the boundary of the interval.
	
	First, we consider the case where \( \lambda = 0\) lies in the interior of \( [\gamma_{\min}, \gamma_{\max}]\). 
	From \cref{lemma: root curve H eta}, it follows that the equation \( H_{\eta}(\lambda, \mu_{\eta}(\lambda)) = 0 \) has a unique solution \( \mu_{\eta}(\lambda) \) for a fixed \( \lambda \in [\gamma_{\min}, \gamma_{\max}]\). The idea of the proof is now to use the root curve \( \mu_{\eta}\) to define an auxiliary function:
	\begin{align*}
		g \colon [\gamma_{\min}, \gamma_{\max}] \rightarrow \R, 
		\quad g(\lambda) = H_p(\lambda, \mu_{\eta}(\lambda)),
	\end{align*}
	and show that \(g\) has at least one root in \( [\gamma_{\min}, \gamma_{\max}]\) by applying the intermediate value theorem. Because once we have a root \( \lambda^{\ast}\) of \( g\), we set \( \mu^{\ast} = \mu_{\eta}(\lambda^\ast)\).
	All that remains, is to verify the requirements of the intermediate value theorem:
	\begin{multicols}{2}
		\begin{enumerate}[label=(\alph*)]
			\item \(g\) is continuous on \( [\gamma_{\min}, \gamma_{\max}]\), \label{item: continuous}
			\item \( g(\gamma_{\min}) \le 0\) and \( g(\gamma_{\max}) \ge 0\). \label{item: opposite sign}
		\end{enumerate}
	\end{multicols}
	We start by proving the continuity of \(g\).
	From \cref{lemma: root curve continuous}, it follows that the root curve  \( \mu_{\eta}\) is continuous on \( [\gamma_{\min}, \gamma_{\max}] \setminus \{0\}\). Furthermore, we know that \( H_p\) is continuous (\cref{lemma: Hp continuous}), which immediately provides the continuity of \( g\) for \( \lambda \neq 0\). 
	Thus, it remains to show the continuity of \( g\) in \( \lambda = 0\).
	Recall, that from \cref{lemma: Hp continuous}, it also follows that \( H_p(0,\cdot)\) is constant, i.e., there exits an \( \alpha \in \R\) such that \(H_p(0, \mu) = \alpha \) for all \( \mu \in [0,1]\) which leads to:
	\begin{align*}
		\lim_{\lambda \rightarrow 0-} g(\lambda) 
		& = \lim_{\lambda \rightarrow 0-}  H_p(\lambda, \mu_{\eta}(\lambda)) 
		= H_p(0, \lim_{\lambda \rightarrow 0-} \mu_{\eta}(\lambda))
		= \alpha \\
		\alpha & =  H_p(0, \lim_{\lambda \rightarrow 0+}  \mu_{\eta}(\lambda))
		= \lim_{\lambda \rightarrow 0+}  H_p(\lambda, \mu_{\eta}(\lambda)) 
		= \lim_{\lambda \rightarrow 0+} g(\lambda).
	\end{align*} 
	
	After we have established the continuity of \( g\), requirement~\labelcref{item: opposite sign} follows immediately from \cref{lemma: unidirectional flow}, which completes the proof for this case.

	Next, we consider the case where \( \lambda = 0\) lies on the boundary of \( [\gamma_{\min}, \gamma_{\max}]\). Then either \( \gamma_{\min}= 0\) or \( \gamma_{\max} = 0\). Let us assume that \( \gamma_{\min}= 0\) as the other case follow analogously.
	Since \( H_p\) is continuous on \( [\gamma_{\min}, \gamma_{\max}] \) and constant in \( \lambda = 0\), there exits an \( \varepsilon > 0\) independent of \( \mu\) such that 
	\begin{equation*}
		\sign(H_p(0,\mu)) = \sign(H_p(\varepsilon, \mu)).
	\end{equation*}
	Hence, we can apply the same argumentation to the interval  \( [\varepsilon, \gamma_{\max}]\) as for the case where \( \lambda = 0\) lies in the interior of the interval \( [\gamma_{\min}, \gamma_{\max}] \).
\end{proof}

After we have proved \cref{lemma: common root}, we automatically obtain that the mixture model~\labelcref{eq: gas flow} has at least one solution, which proves our main result, \cref{thm: existence cycle}. 

\begin{proof}[Proof of \cref{thm: existence cycle}]
    Due to \cref{lemma: solvability 2}, the solvability of the mixture model~\labelcref{eq: gas flow} is equivalent to the functions \( H_p\) and \(H_{\eta}\) having at least one common root, which follows from \cref{lemma: common root} and thus proves the claim.
\end{proof}

\begin{remark}
    So far, we have considered (passive) networks with at most one cycle, but the argumentation can also be applied to networks with active elements such as compressors, as long as they are not in parallel or form cycles. Furthermore, the compressor ratio must be independent of the composition and flow.

    For tree-shaped networks, the flows \( q_e\) and the compositions \( \eta_v\) are uniquely determined by the conservation of mass and the mixing at the nodes, regardless of whether the network contains active elements or not.
    Thus, a solution for the pressure \(p_v\) exists and is also unique for networks with compressors.
    To prove the existence of a solution for single cycle networks, we use the function \( H_p\). For networks with one cycle and compressors outside the cycle, the structure of \( H_p\) is the same as in \cref{eq: H_p alternative}. Therefore, proving the existence of solutions is identical to proving \cref{thm: existence cycle}.
    If compressors are inside the cycle, the compressor ratios enter the function \( H_p\), i.e., some summands are multiplied by the compressor ratios. Which summands are affected depends on the position of the compressors in the cut cycle. As the compressor ratios are constant, again, the idea of the proof of \cref{thm: existence cycle} still applies.
\end{remark}

%% file: tikz/cut_graph.tex
\scalebox{0.7}{
	\begin{tikzpicture}[
		> = {Latex[scale=1.0]}, 
		line width=1pt, 
		every node/.style={circle, draw, minimum size=0pt, inner sep=2.5pt}
		]
		
		\pgfmathsetmacro{\L}{2}
		\pgfmathsetmacro{\sqrtL}{2}
		\pgfmathsetmacro{\h}{0.5}

		\node (v1) at (-\L,0) {\(v_0\)};
		\node (v5) at (0, 0) {\(v_4\)};
		\node (v2) at (\sqrtL, \L + \sqrtL) {\(v_1\)};
		\node (v6) at (\sqrtL, \sqrtL) {\(v_5\)};
		\node (v8) at (\sqrtL, -\sqrtL) {\(v_7\)};
		\node (v7) at (2*\sqrtL, 0) {\(v_6\)};
		\node (v3) at (\L + 2*\sqrtL, 0) {\(v_2\)};
		\node (v4) at (\sqrtL, - \L - \sqrtL) {\(v_3\)};

		\draw[->] (v1) -- (v5);
		\draw[->, blue!70!black] (v8) -- (v7);
		\draw[->]  (v5) -- (v6);
		\draw[->] (v5) -- (v8);
		\draw[->] (v2) -- (v6);
		\draw[->] (v6) -- (v7);
		\draw[->] (v8) -- (v4);
		\draw[->] (v7) -- (v3);
		
		\draw[red!80!black, very thick] (1.5*\sqrtL - \h, -0.5*\sqrtL + \h) -- (1.5*\sqrtL + \h, -0.5*\sqrtL - \h);

		\draw[->, dashed] (1.5*\sqrtL, -0.5*\sqrtL) to[bend right] (2* \L + 2* \sqrtL + 2, 
			0.5* \sqrtL - 0.4);
			
		\draw[->, dashed] (1.5*\sqrtL, -0.5*\sqrtL)  .. controls (2* \L + 3* \sqrtL + 2, -1.5*\sqrtL) ..  (2* \L + 4* \sqrtL + 2, 
		0.5* \sqrtL - 0.4);

		\begin{scope}[shift={({2* \L + 2* \sqrtL + 2}, 1.5*\sqrtL)}]
			\node (v1) at (-\L,0) {\(v_0\)};
			\node (v5) at (0, 0) {\(v_4\)};
			\node (v2) at (\sqrtL, \L + \sqrtL) {\(v_1\)};
			\node (v6) at (\sqrtL, \sqrtL) {\(v_5\)};
			\node (v8) at (\sqrtL, -\sqrtL) {\(v_7\)};
			\node (v7) at (2*\sqrtL, 0) {\(v_6\)};
			\node (v3) at (\L + 2*\sqrtL, 0) {\(v_2\)};
			\node (v4) at (\sqrtL, - \L - \sqrtL) {\(v_3\)};
			\node[red!80!black] (vl) at (0, - \sqrtL) {\(v_\cl\)};
			\node[red!80!black]  (vr) at (2*\sqrtL,  -\sqrtL) {\(v_\crr\)};
			
			\draw[->] (v1) -- (v5);
			\draw[->]  (v5) -- (v6);
			\draw[->] (v5) -- (v8);
			\draw[->] (v2) -- (v6);
			\draw[->] (v6) -- (v7);
			\draw[->] (v8) -- (v4);
			\draw[->] (v7) -- (v3);
			\draw[->, blue!70!black] (v8) -- (vl);
			\draw[->, blue!70!black] (vr) -- (v7);
		\end{scope}
\end{tikzpicture}}

%% file: tikz/overview.tex
\scalebox{0.75}{
	\begin{tikzpicture}[
		> = {Latex[scale=1.0]}, 
		]
		
		\tikzstyle{dataset} = [rectangle, draw, text width = 2in, rounded corners = 5pt, thick]
		
		\pgfmathsetmacro{\hx}{6}
		\pgfmathsetmacro{\hy}{2}
		\pgfmathsetmacro{\hyy}{1.75}
		
		\pgfmathsetmacro{\d}{0.3}
		\pgfmathsetmacro{\dx}{1.5}
		\pgfmathsetmacro{\dxx}{0.5}
		\pgfmathsetmacro{\dy}{1.5}
		\pgfmathsetmacro{\dyy}{0.25}
		
		\draw[rounded corners = 5pt, thick] (0, \hyy + \dyy) rectangle (\hx, \dyy + \hy + \hyy) node  [rectangle split, rectangle split parts = 2, midway] {For fixed \(\lambda \in \R \) the function \nodepart{second} \( H_{\eta}(\lambda, \mu)\) has a unique root \( \mu_{\eta}(\lambda)\).};
			
		\draw[rounded corners = 5pt, thick] (0,0) rectangle (\hx, \hyy)  node [midway, align=right] {\(\bullet\) \cref{lemma: composition linear}  \hspace{7pt} \(\quad \bullet \) \cref{lemma: root curve H eta}};
		
		\draw[rounded corners = 5pt, thick] (\hx + \dx,\hyy + \dyy) rectangle (2*\hx+ \dx,  \dyy + \hy + \hyy) node  [rectangle split, rectangle split parts = 2,midway] {The function \(g = H_p(\cdot, \mu_{\eta}(\cdot))\)  \nodepart{second} is continuous.};
		
		\draw[rounded corners = 5pt, thick] (\hx + \dx, 0) rectangle (2*\hx+\dx, \hyy) 
		node  [rectangle split, rectangle split parts = 2, midway]
		{\(\bullet \) \cref{korollar: flow continuous} 
			\(\quad \bullet \) \cref{lemma: composition continuous} 
			\nodepart{second} \(\bullet \) \cref{lemma: Hp continuous} \hspace{7pt}
			\(\quad \bullet \) \cref{lemma: root curve continuous}};
		
		\draw[rounded corners = 5pt, thick] (\dx + \dxx+2*\hx,\hyy + \dyy) rectangle (\dx + \dxx+3*\hx,  \dyy + \hy + \hyy) node  [rectangle split, rectangle split parts = 2, midway] {There exit values \( \lambda^{\pm}\) such that \nodepart{second} \(g(\lambda^-) \le 0 \) and \(g(\lambda^+) \ge 0 \).};
		
		\draw[rounded corners = 5pt, thick, align=left] (2*\hx + \dx + \dxx, 0) rectangle (3*\hx+\dx+\dxx, \hyy)
		node  [midway, align=left] {\(\bullet\) \cref{lemma: flow linear}  \( \quad \bullet\) \cref{lemma: unidirectional flow}};
		
		\draw[rounded corners = 5pt, thick] (0.5 * \hx,-\dy) rectangle (\dx + \dxx + 2.5*\hx, -\dy-\hyy) node  [midway] {Apply the \emph{intermediate value theorem} to \( g(\lambda) = H_p(\lambda, \mu_{\eta}(\lambda))\).};
		
		\draw[rounded corners = 5pt, thick] (\hx + \dx,-2*\dy - \hyy) rectangle (2*\hx + \dx, -2*\dy - \hyy - \hy) node  [rectangle split, rectangle split parts = 2, midway] {\textbf{\cref{lemma: common root}} \nodepart{second} \( H_p\) and \( H_{\eta} \) have a common root};
		
		\node[above] at ({\hx / 2} , {\dyy + \hy + \hyy + \d}) {\textbf{Step~\labelcref{item: step 1}}};
		\draw[rounded corners = 5pt, very thick, dashed]
		(-\d, - \d) rectangle (\hx + \d, {\dyy + \hy + \hyy + \d});
		
		\node[above] at (\dxx / 2 + \dx + 2*\hx + \d, {\dyy + \hy + \hyy + \d}) {\textbf{Step~\labelcref{item: step 2}}};
		\draw[rounded corners = 5pt, very thick, dashed]
		(\hx + \dx -\d, -\d) rectangle (\dxx + \dx+3*\hx + \d,  \dyy + \hy + \hyy + \d);

		\draw[->, thick] (0.5 * \hx, -\d) -- (\hx + 0.5*\dx, -\dy);
		\draw[->, thick] (1.5 * \hx + \dx, -\d) -- (1.5*\hx +\dx, -\dy);
		\draw[->, thick] (2.5*\hx + \dx + \dxx , -\d) -- (2*\hx + \dx + \dxx/2, -\dy);
		\draw[->, thick] (1.5 * \hx + \dx, -\dy - \hyy) -- (1.5*\hx +\dx, -2*\dy +- \hyy);
	\end{tikzpicture}
}

%% file: tex/experiments.tex
\section{Numerical Examples}\label{sec: experiments}
So far, we restricted our theoretical investigation to networks with at most one cycle and skipped networks with multiple cycles such as the diamond network shown in \cref{fig: diamond network}.
The idea to prove the existence of a solution for these kind of networks is to use an induction over the number of cycles in the network.
In \cref{sec: existence single cycle}, we have presented the first induction step going from networks without cycles to networks with one cycle.
However, networks with multiple cycles pose new challenges when performing the induction since the cut graph of such networks is no longer tree-shaped.

One consequence is that the flow \(q_e^c\) on the edges of the cut graph is not fully determined by the loads \(b_v^c\), i.e., it also depends on the supply composition \( \zeta_v^c\)and the fixed pressure \(p^{\ast}\).
As the flow \(  q_e^c\) now depends on both parameters, \( \lambda\) and \( \mu\), the composition \( \eta_v^c\) at a node \(v\) becomes non-linear in \(\mu\).
The function \(H_{\eta}\) inherits the non-linearity and it is becomes unclear whether the equation \( H_{\eta}(\lambda, \mu) = 0\) admits a (unique) solution for every \( \lambda \in \R\).
Thus, the analysis of the non-linearity and its influence on the solvability of \( H_{\eta}(\lambda, \mu) = 0\) remains subject to further research.

\medbreak
Therefore, we dedicate this section to a numerical study, which purpose is twofold.
On the one hand, we visualize the idea of the proof of \cref{thm: existence cycle} by computing the functions \( H_p(\lambda, \mu)\) and \( H_{\eta}(\lambda, \mu)\) and calculating the root curve \( \mu_{\eta}(\lambda)\) for a sample network with exactly one cycle.
On the other hand, we apply the idea of the proof to a network with multiple cycles to showcase that the idea might also be applicable to networks with multiple cycles instead of  networks with one cycle only.
In particular, we compute the restriction \( H_p(\lambda, \mu_{\eta}(\lambda))\) and prove, again numerically, that this function has a unique root, which implies with \cref{lemma: solvability 2} that the mixture model~\labelcref{eq: gas flow} admits at least one solution.

\subsection{A Network with One Cycle}\label{sec: experiments one cycle}
In this section, we compute the functions \( H_p \) and \( H_{\eta}\), and the root curve \( \mu_{\eta}\) for a sample network with only one cycle. Furthermore, we link the theoretical results of \cref{sec: existence single cycle} to our numerical example.
As a sample network we use the (cut) network displayed in \cref{fig: cut graph}.
This network has exactly one cycle, two supply nodes, \( v_0\) and \( v_1\), and two demand nodes, \( v_2\) and \( v_3\). At the remaining nodes, the mass is conserved, that is, i.e., \( b_v = 0\). In \cref{tab: network data} we provide the specific values of the boundary data.

\renewcommand{\arraystretch}{1.25}
\begin{table}[ht]
\caption{Boundary data at the nodes} 
\centering 
\begin{tabular}{cccccc} %
\toprule
    \multicolumn{2}{c}{supply} & & \multirow{2}{*}{demand} & & 
        \multirow{2}{*}{pressure}  \\ 
    \cmidrule{1-2} 
    {composition} & {load} & &  \\
    \midrule
    \(\zeta_{v_0} = \frac{3}{4}\) & \( b_{v_0} = 4\) & &
        \(b_{v_2} = - 2\) & &
        \multirow{2}{*}{\( p_{v_0} = p^{\ast} = 60\)} \\
    \(\zeta_{v_1} = \frac{1}{4}\) & \( b_{v_1} = 4 \) & & 
         \(b_{v_3} = -6\)  & &  \\
\bottomrule
\end{tabular}
\label{tab: network data}
\end{table}

\begin{figure}[ht]
	\centering
	\begin{subfigure}[t]{0.475\textwidth}
		\centering
		\includegraphics[width=\textwidth]{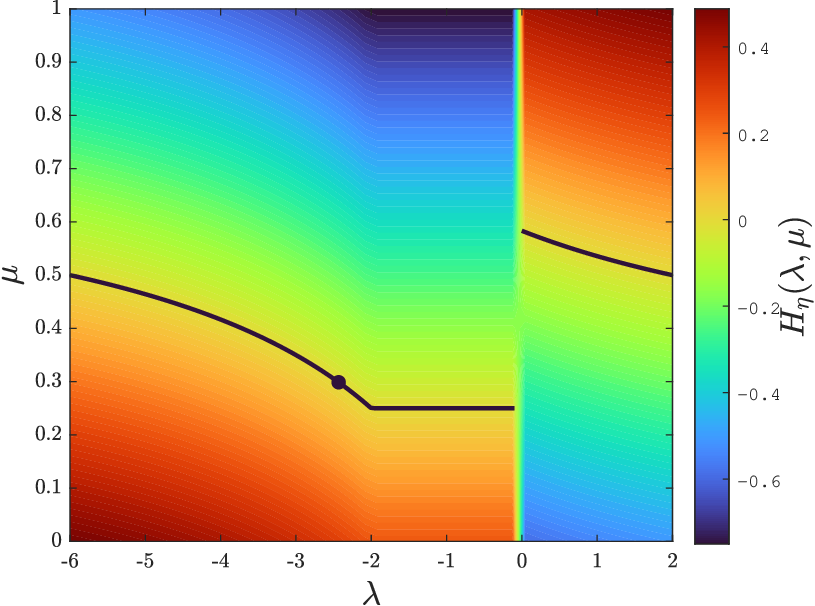}
	\end{subfigure} \hfill
	\begin{subfigure}[t]{0.475\textwidth}
		\centering
		\includegraphics[width=\textwidth]{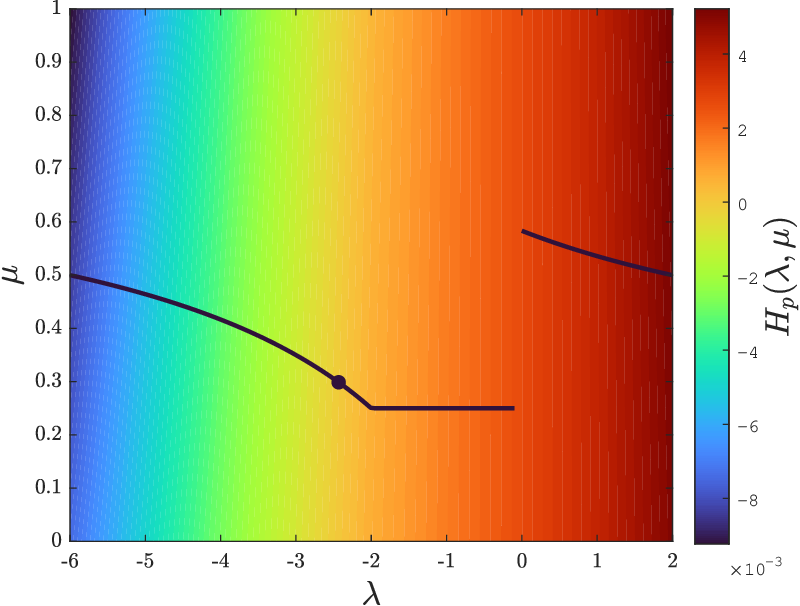}
	\end{subfigure} 
	\caption{The functions \(H_{\eta}(\lambda,\mu)\) (left) and \( H_{p}(\lambda, \mu)\) (right) for the sample network in \cref{fig: cut graph}. The black dot is the common root \( (\lambda^{\ast}, \mu^{\ast})\), and the black line the root curve \( \mu_{\eta}(\lambda)\).}
	\label{fig: visual proof}
\end{figure}

We can narrow down the range of \(\lambda\), for which the function \( H_p\) can become zero, to the interval \( [\gamma_{\min}, \gamma_{\max}] = [-6, 2]\) due to \cref{lemma: unidirectional flow} and the boundary data given in \cref{tab: network data}.
Thus, we evaluate the functions \( H_p\) and \( H_{\eta}\) on the domain  \( [-6,2] \times [0,1]\).
As the evaluation of these functions at a fixed point \(  (\lambda, \mu)\) requires the solution of the mixture model~\labelcref{eq: gas flow}, we discretize the intervals \( [-6,2]\) and \( [0,1]\) by uniform grids of size \( N_{\lambda}=50\) and \( N_{\mu}=51\), respectively, and solve the mixture model~\labelcref{eq: gas flow} only at these grid points.
To obtain a solution to the non-linear system corresponding to the mixture model~\labelcref{eq: gas flow}, we invoke the Levenberg-Marquardt method.
The resulting functions \( H_p\) and \( H_{\eta}\) are displayed in \cref{fig: visual proof}.

First, we observe that the function \( H_{\eta}\) is discontinuous in \( \lambda = 0\) for any value of \(\mu \in [0,1]\), cf. \cref{lemma: composition continuous} and \cref{lemma: root curve H eta}. This discontinuity is due to the sudden switch of roles at \( \lambda = 0\) at the nodes \( v_\cl\) and \(v_\crr\). For \( \lambda > 0\), the node \( v_\cl\) is a demand node and for \( \lambda <0 \) it is a supply node. For the node \( v_\crr\) it is vice versa.
Moreover, we see that the function \( H_p\) is continuous in the entire domain, cf. \cref{lemma: Hp continuous} and that \( H_p(\gamma_{\min}, \mu) \le 0\) and \(H_p(\gamma_{\max}, \mu) \ge 0\) for any value of \( \mu \in [0,1]\), cf. \cref{lemma: unidirectional flow}.

Besides the functions \( H_p\) and \( H_{\eta}\), \cref{fig: visual proof} also shows the root curve \( \mu_{\eta}\).
As the cut network in \cref{fig: cut graph} is tree-shaped, the flows \(q_e^c\) on the edges are fully determined by the loads \( b_v^c\). Thus, the computation of the nodal compositions \( \eta_v^c\) is straightforward and by \cref{lemma: root curve H eta} we can compute the root curve \( \mu_{\eta}\) explicitly:
\begin{equation}\label{eq: root curve}
    \mu_{\eta}(\lambda) = \begin{cases}
        \dfrac{\zeta_{v_0} \left[ b_{v_1} + b_{v_2} + \lambda\right] - \zeta_{v_1} b_{v_1}}{b_{v_2} + \lambda}, 
            & \text{if } \lambda \in [b_{v_3}, b_{v_0} + b_{v_3}], \\[12pt]
        \zeta_{v_1}, & \text{if }\lambda \in[b_{v_0} + b_{v_3}, 0], \\[10pt]
        \dfrac{\zeta_{v_0} b_{v_0} - \zeta_{v_1} \left[ b_{v_0} + b_{v_3} - \lambda\right]}{\lambda - b_{v_3}}, 
            & \text{if } \lambda \in [0, -b_{v_2}].
    \end{cases}    
\end{equation}
Notice that \( b_{v_0} + b_{v_3} = - (b_{v_1} + b_{v_2}) \) holds due to mass conservation and that \( b_{v_0} + b_{v_3} \le 0\) holds for the boundary data given in \cref{tab: network data}. Also note that the root curve \( \mu_{\eta}\) is discontinuous in \(\lambda = 0\) (unless \( b_{v_0} = 0\)), which is also visible in \cref{fig: visual proof} and backed-up by \cref{lemma: root curve continuous}.

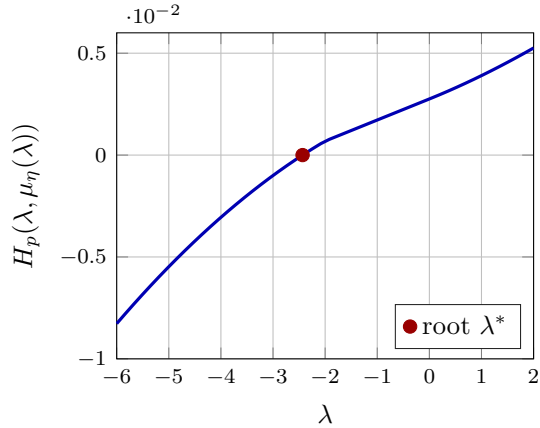
\begin{figure}[ht]
	\centering
    \begin{subfigure}[t]{0.475\textwidth}
		\centering
		\input{./tikz/restriction_cycle}
	\end{subfigure}
	\caption{The function \(H_p\) restricted to the root curve \( \mu_{\eta}(\lambda)\) for the sample network in \cref{fig: cut graph}, resulting in the scalar function \( g(\lambda) =  H_p(\lambda,\mu_{\eta}(\lambda))\) from the proof of \cref{thm: existence cycle}.}
	\label{fig: restriction cycle}
\end{figure}

With this explicit representation, we can evaluate the function \( H_p\) along the root curve \( \mu_{\eta}\), i.e., the function \(  H_p(\lambda,\mu_{\eta}(\lambda))\), by solving the mixture model~\labelcref{eq: gas flow} on the cut graph. Again, to evaluate the function, we discretize the interval \( [-6,2]\) by a uniform grid of size \( N_{\lambda} = 50\) and invoke the Levenberg-Marquardt method.
The function \( H_p(\lambda,\mu_{\eta}(\lambda))\) is displayed in \cref{fig: restriction cycle}.

First of all, \cref{fig: restriction cycle} shows that the restriction \( H_p(\lambda,\mu_{\eta}(\lambda))\) is continuous in \( \lambda = 0\) even though the root curve \( \mu_{\eta}\) is not, which follows from \cref{lemma: common root}. Moreover, \cref{fig: restriction cycle} indicates that the function  \( H_p(\lambda,\mu_{\eta}(\lambda))\) is strictly monotonically increasing, which implies that its root \( \lambda^{\ast}\) is unique.
Then, \cref{lemma: solvability 2} would imply that solutions to the mixture model~\labelcref{eq: gas flow} are unique, which would be analogous to the \glqq classic\grqq{} gas flow model involving only one gas instead of a mixture of multiple gases.

\subsection{A Network with Multiple cycles} \label{sec:multipleCycles}
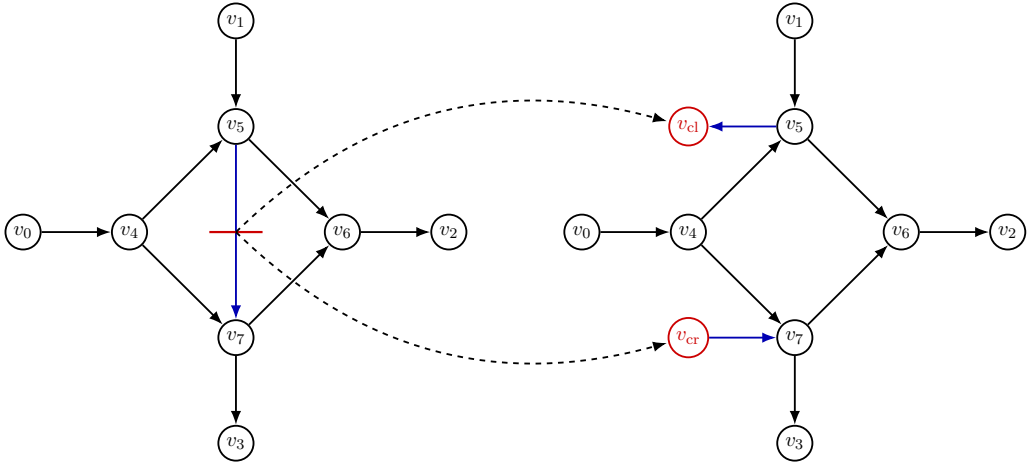
\begin{figure}[ht]
	\centering
	\input{./tikz/diamond_network_cut}
	\caption{A graph with multiple, intertwined cycles (left) and the resulting cut graph when cutting the blue, inner edge (right).}
	\label{fig: diamond network}
\end{figure}

The subject of this section is to analyze whether the idea of the proof of \cref{thm: existence cycle} is also applicable to networks with multiple cycles.
Therefore, we perform the same computations as in \cref{sec: experiments one cycle} but for the example of the diamond network shown in  \cref{fig: diamond network}.
This network has the same supply and demand nodes as the network shown in \cref{fig: cut graph} but an additional edge connecting the nodes \( v_5\) and \( v_7\). 
We also assume the same boundary data as before.

Again, we can narrow down the range of \( \lambda\), for which the function \( H_p\) becomes zero. However, this time we have to invoke \cite[Lemma 14]{WintergerstDiss}, which yields the interval:
\begin{equation*}
    \bigg[ - \sum_{v \in \V} \vert b_v \vert, \, \sum_{v \in \V} \vert b_v \vert \bigg] = [-8,8].
\end{equation*}

\begin{figure}[ht]
	\centering
    \begin{subfigure}[t]{0.475\textwidth}
		\centering
		\input{./tikz/composition_nonlinear_diamond}
	\end{subfigure} 
	\caption{The nodal composition \( \eta_{v_4}^c(\lambda,\mu)\) as a function of \( \mu\) for different values of \( \lambda\). The value of \( \lambda\) is encoded by the color bar.}
    \label{fig: composition nonlinear}
\end{figure}
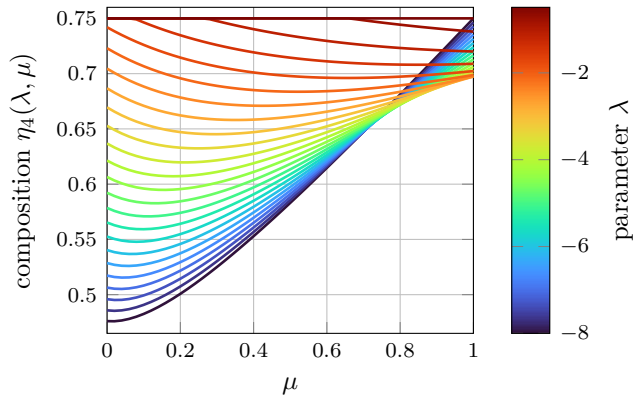

\medbreak
Before we delve into the analysis of the functions \( H_p\) and \( H_{\eta}\), we emphasize the main difficulty when transferring the idea of the proof of \cref{thm: existence cycle} to networks with multiple cycles:
The non-linear dependence of the composition \( \eta_v^c\) on \( \mu\), which makes the solving \( H_{\eta}(\lambda, \mu) = 0\) with respect to \( \mu\) more challenging.
The non-linearity of the composition becomes evident in \cref{fig: composition nonlinear}, which shows the composition \( \eta_{v_4}^c\) as a function of \( \mu\) for different values of \( \lambda \in [-8,8] \). The color bar encodes the values of the parameter \( \lambda\). Since the composition \( \eta_{v_4}^c \) is constant for \( \lambda \ge -0.5 \), specifically \( \eta_{v_4}^c(\lambda,\mu) \equiv 0.75\), we omit these values in the plot.

\begin{figure}[ht]
	\centering
	\begin{subfigure}[t]{0.475\textwidth}
		\centering
		\includegraphics[width=\textwidth]{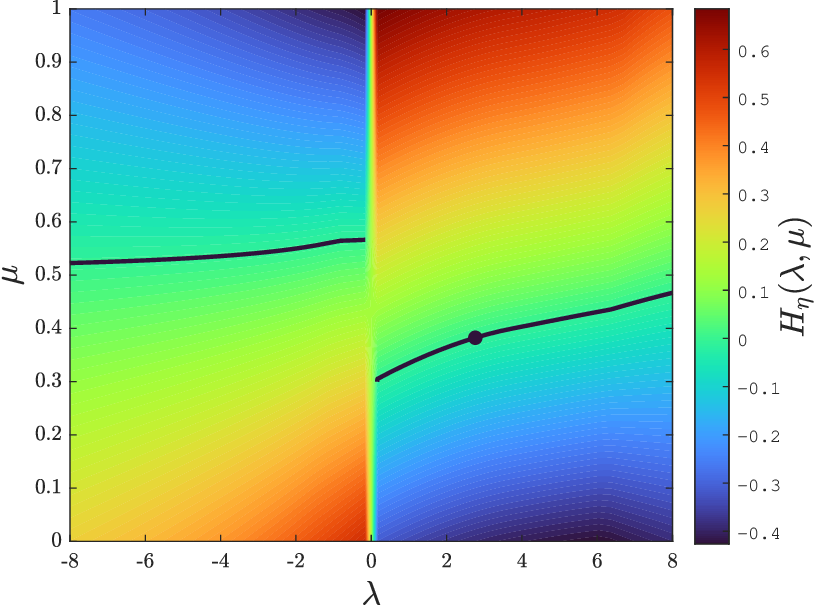}
	\end{subfigure} \hfill
	\begin{subfigure}[t]{0.475\textwidth}
		\centering
		\includegraphics[width=\textwidth]{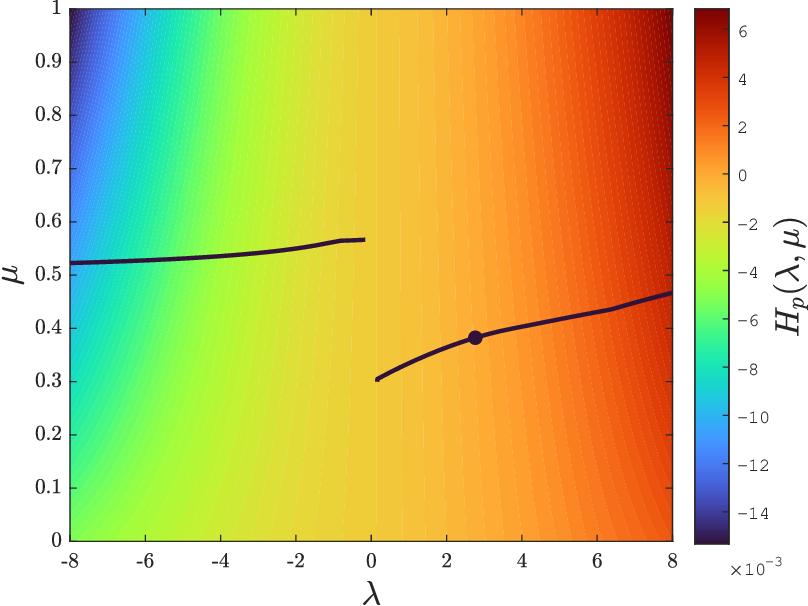}
	\end{subfigure}
	\caption{The functions \(H_{\eta}(\lambda,\mu)\) (left) and \( H_p(\lambda, \mu)\) (right) for the network in \cref{fig: diamond network}. The black dot is the common root \( (\lambda^{\ast}, \mu^{\ast})\), and the black line the root curve \( \mu_{\eta}(\lambda)\).}
	\label{fig: visual proof diamond}
\end{figure}

\medbreak
Next, we analyze the functions \( H_p\) and \( H_{\eta}\) for the sample (cut) network shown in \cref{fig: diamond network}.
Both functions are displayed in \cref{fig: visual proof diamond}.
Similar to the network with one cycle, we observe that the function \( H_{\eta}\) is discontinuous in \( \lambda = 0\) for any value of \( \mu \in[0,1]\) and the function \( H_p\) has opposite signs for \( \lambda = -8\) and \( \lambda = 8\).

In \cref{fig: visual proof diamond}, we also depict the root curve \( \mu_{\eta}\) of the function \(H_{\eta}\). 
The cut graph of the diamond network shown in \cref{fig: cut graph} is not tree-shaped. This makes the computation of the exact root curve more intricate and we compute an approximation of the exact root curve \( \mu_{\eta}\) instead by extracting the coordinates of the level line \( H_{\eta}(\lambda, \mu) = 0\) from the plot using Python's \texttt{get\_paths()} function for contour plots.
As for the network with one cycle, we see that the root curve \( \mu_{\eta}\) is discontinuous in \( \lambda = 0\), indicating that \cref{lemma: composition continuous} also hold for networks with multiple cycles. Moreover, \cref{fig: visual proof diamond} reveals that \( H_{\eta}(\lambda, \mu) = 0\) has a unique solution for \( \lambda \in [-8,8] \) fixed.

\begin{figure}[ht]
    \centering
    \begin{subfigure}[t]{0.475\textwidth}
		\centering 
		\input{./tikz/restriction_diamond}
	\end{subfigure}
	\caption{The functions \(H_p(\lambda,\mu)\) (left) and \( H_{\eta}(\lambda, \mu)\) (right) for the network in \cref{fig: diamond network}. The black dot is the common root \( (\lambda^{\ast}, \mu^{\ast})\), and the black line the root curve \( \mu_{\eta}(\lambda)\).}
	\label{fig: restriction diamond}
\end{figure}
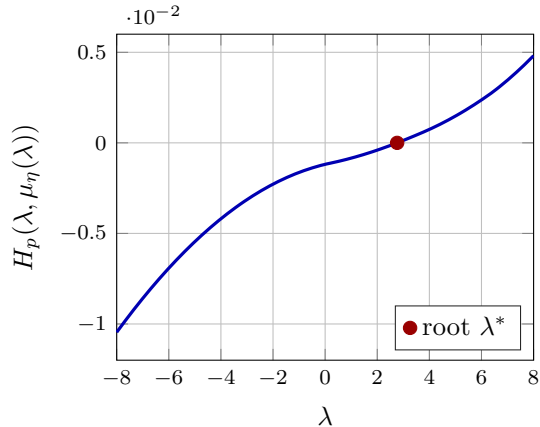

The approximation of the root curve \( \mu_{\eta} \) allows us to compute the function \( H_p(\lambda, \mu_{\eta}(\lambda))\) by evaluating the function \( H_p\) at the corresponding coordinates. We provide the result in \cref{fig: restriction diamond}.
Again the function  \( H_p(\lambda, \mu_{\eta}(\lambda))\) is continuous in \( \lambda = 0\) despite the discontinuity of the root curve \( \mu_{\eta}\). Furthermore, we observe that  \( H_p(\lambda, \mu_{\eta}(\lambda))\) is strictly monotonically increasing with a unique root \( \lambda^{\ast}\), implying that the solution of the mixture model~\labelcref{eq: gas flow} is unique, cf. \cref{lemma: solvability 2}.

In summary, the numerical experiments suggest that the idea of the proof of \cref{thm: existence cycle} also applies to networks with multiple cycles. However, to successfully extend the proof, we need to rigorously prove that \( H_{\eta}(\lambda,\mu) = 0\) admits a unique solution for these types of networks.

%% file: tikz/restriction_cycle.tex
%
\definecolor{mycolor1}{rgb}{0.00000,0.44700,0.74100}%
\definecolor{mycolor2}{rgb}{0.85000,0.32500,0.09800}%
\begin{tikzpicture}

\begin{axis}[%
width=0.742\textwidth, 
height=0.585\textwidth, 
at={(0\textwidth,0\textwidth)},
scale only axis,
xmin=-6,
xmax=2,
xlabel style={font=\color{white!15!black}, font=\small},
xlabel={\(\lambda\)},
xtick={-6,-5,-4,-3,-2,-1,0,1,2},
ymin=-0.01,
ymax=0.006,
ylabel style={font=\color{white!15!black}, font=\small},
ylabel={\(H_p(\lambda, \mu_{\eta}(\lambda))\)},
ticklabel style = {font=\scriptsize},
axis background/.style={fill=white},
xmajorgrids,
ymajorgrids,
legend style={at={(0.97,0.03)}, anchor=south east, legend cell align=left, align=left, draw=white!15!black, font = \small}
]
\addplot [very thick, color=blue!70!black, forget plot]
  table[row sep=crcr]{%
-6	-0.00827262440289989\\
-5.83673469387755	-0.00779090521538084\\
-5.6734693877551	-0.00731914307293202\\
-5.51020408163265	-0.00685733181711612\\
-5.3469387755102	-0.00640546474207326\\
-5.18367346938776	-0.00596353453233833\\
-5.02040816326531	-0.00553153319191146\\
-4.85714285714286	-0.00510945196328638\\
-4.69387755102041	-0.00469728123456115\\
-4.53061224489796	-0.0042950104326451\\
-4.36734693877551	-0.00390262790002627\\
-4.20408163265306	-0.00352012075208741\\
-4.04081632653061	-0.00314747471127852\\
-3.87755102040816	-0.00278467391365744\\
-3.71428571428571	-0.00243170068225496\\
-3.55102040816327	-0.00208853526046937\\
-3.38775510204082	-0.00175515549694083\\
-3.22448979591837	-0.0014315364714097\\
-3.06122448979592	-0.00111765004808539\\
-2.89795918367347	-0.000813464339745273\\
-2.73469387755102	-0.000518943060926258\\
-2.57142857142857	-0.000234044742528638\\
-2.40816326530612	4.12782281205004e-05\\
-2.24489795918367	0.000307080788715552\\
-2.08163265306122	0.000563426666984057\\
-1.91836734693878	0.000772427617623528\\
-1.75510204081633	0.000941130632568843\\
-1.59183673469388	0.00110983364751727\\
-1.42857142857143	0.00127853666246436\\
-1.26530612244898	0.00144723967741101\\
-1.10204081632653	0.0016159426923581\\
-0.938775510204082	0.00178464570730519\\
-0.775510204081633	0.00195334872225184\\
-0.612244897959184	0.00212205173719904\\
-0.448979591836735	0.00229075475214602\\
-0.285714285714286	0.00245945776709355\\
-0.122448979591837	0.00262816078203976\\
0.0408163265306118	0.002797059203691\\
0.204081632653061	0.0029704133272892\\
0.367346938775509	0.00314985372553878\\
0.530612244897958	0.00333525190049677\\
0.693877551020408	0.00352649189060306\\
0.857142857142857	0.00372346877825658\\
1.02040816326531	0.00392608740563283\\
1.18367346938775	0.00413426126561922\\
1.3469387755102	0.00434791154062353\\
1.51020408163265	0.00456696626675734\\
1.6734693877551	0.00479135960472021\\
1.83673469387755	0.00502103120181374\\
2	0.00525592563210076\\
};

\addplot[only marks, mark=*, mark options={}, mark size=2.5000pt, color=red!60!black, fill=red!60!black] table[row sep=crcr]{%
x	y\\
-2.43300609834716	0\\
};
\addlegendentry{root \( \lambda^{\ast} \)}

\end{axis}
\end{tikzpicture}%

%% file: tikz/diamond_network_cut.tex
\scalebox{0.7}{
	\begin{tikzpicture}[
		> = {Latex[scale=1.0]}, 
		line width=1pt, 
		every node/.style={circle, draw, minimum size=0pt, inner sep=2.5pt}
		]
		
		\pgfmathsetmacro{\L}{2}
		\pgfmathsetmacro{\sqrtL}{2}
		\pgfmathsetmacro{\h}{0.5}

		\node (v1) at (-\L,0) {\(v_0\)};
		\node (v5) at (0, 0) {\(v_4\)};
		\node (v2) at (\sqrtL, \L + \sqrtL) {\(v_1\)};
		\node (v6) at (\sqrtL, \sqrtL) {\(v_5\)};
		\node (v8) at (\sqrtL, -\sqrtL) {\(v_7\)};
		\node (v7) at (2*\sqrtL, 0) {\(v_6\)};
		\node (v3) at (\L + 2*\sqrtL, 0) {\(v_2\)};
		\node (v4) at (\sqrtL, - \L - \sqrtL) {\(v_3\)};

		\draw[->] (v1) -- (v5);
		\draw[->] (v8) -- (v7);
		\draw[->]  (v5) -- (v6);
		\draw[->] (v5) -- (v8);
		\draw[->] (v2) -- (v6);
		\draw[->] (v6) -- (v7);
	    \draw[->, blue!70!black] (v6) -- (v8);
		\draw[->] (v8) -- (v4);
		\draw[->] (v7) -- (v3);
		
		\draw[red!80!black, very thick] (\sqrtL - \h, 0) -- (\sqrtL + \h, 0);
				
		\draw[->, dashed] (\sqrtL, 0) to[bend left] (2* \L + 2* \sqrtL + 2.1, \sqrtL + 0.1);			
		\draw[->, dashed] (\sqrtL, 0) to[bend right] (2* \L + 2* \sqrtL + 2.1, -\sqrtL - 0.1);	;

		\begin{scope}[shift={({2* \L + 2* \sqrtL + 2.5}, 0*\sqrtL)}]
			\node (v1) at (-\L,0) {\(v_0\)};
			\node (v5) at (0, 0) {\(v_4\)};
			\node (v2) at (\sqrtL, \L + \sqrtL) {\(v_1\)};
			\node (v6) at (\sqrtL, \sqrtL) {\(v_5\)};
			\node (v8) at (\sqrtL, -\sqrtL) {\(v_7\)};
			\node (v7) at (2*\sqrtL, 0) {\(v_6\)};
			\node (v3) at (\L + 2*\sqrtL, 0) {\(v_2\)};
			\node (v4) at (\sqrtL, - \L - \sqrtL) {\(v_3\)};
			\node[red!80!black] (vr) at (0, -\sqrtL) {\(v_\crr\)};
			\node[red!80!black]  (vl) at (0,  \sqrtL) {\(v_\cl\)};
			
			\draw[->] (v1) -- (v5);
			\draw[->]  (v5) -- (v6);
			\draw[->] (v5) -- (v8);
			\draw[->] (v2) -- (v6);
			\draw[->] (v6) -- (v7);
			\draw[->] (v8) -- (v7);
			\draw[->] (v8) -- (v4);
			\draw[->] (v7) -- (v3);
			\draw[->, blue!70!black] (v6) -- (vl);
			\draw[->, blue!70!black] (vr) -- (v8);
		\end{scope}
\end{tikzpicture}}

%% file: tikz/composition_nonlinear_diamond.tex
%
\definecolor{mycolor1}{rgb}{0.18995,0.07176,0.23217}%
\definecolor{mycolor2}{rgb}{0.23611,0.19817,0.52574}%
\definecolor{mycolor3}{rgb}{0.26504,0.31891,0.75013}%
\definecolor{mycolor4}{rgb}{0.27673,0.43398,0.90534}%
\definecolor{mycolor5}{rgb}{0.27027,0.54356,0.99060}%
\definecolor{mycolor6}{rgb}{0.21753,0.65329,0.98229}%
\definecolor{mycolor7}{rgb}{0.14189,0.75850,0.90002}%
\definecolor{mycolor8}{rgb}{0.09438,0.84893,0.78669}%
\definecolor{mycolor9}{rgb}{0.12514,0.91546,0.67941}%
\definecolor{mycolor10}{rgb}{0.24504,0.96345,0.54588}%
\definecolor{mycolor11}{rgb}{0.41027,0.99246,0.39992}%
\definecolor{mycolor12}{rgb}{0.57299,0.99819,0.27752}%
\definecolor{mycolor13}{rgb}{0.69570,0.97595,0.21283}%
\definecolor{mycolor14}{rgb}{0.80595,0.92375,0.20471}%
\definecolor{mycolor15}{rgb}{0.90030,0.84924,0.22071}%
\definecolor{mycolor16}{rgb}{0.96636,0.76098,0.22726}%
\definecolor{mycolor17}{rgb}{0.99476,0.65977,0.19817}%
\definecolor{mycolor18}{rgb}{0.99074,0.53471,0.14665}%
\definecolor{mycolor19}{rgb}{0.95969,0.40451,0.09039}%
\definecolor{mycolor20}{rgb}{0.90658,0.29010,0.04620}%
\definecolor{mycolor21}{rgb}{0.83369,0.20106,0.02176}%
\definecolor{mycolor22}{rgb}{0.73827,0.12585,0.00783}%
\definecolor{mycolor23}{rgb}{0.62024,0.06411,0.00409}%
\definecolor{mycolor24}{rgb}{0.47960,0.01583,0.01055}%
\begin{tikzpicture}

\begin{axis}[%
width=0.652\textwidth,
height=0.585\textwidth,
at={(0\textwidth,0\textwidth)},
scale only axis,
point meta min=-8,
point meta max=-0.489795918367347,
xmin=0,
xmax=1,
xlabel style={font=\color{white!15!black}, font=\small},
xlabel={\(\mu\)},
ymin=0.465,
ymax=0.76,
ylabel style={font=\color{white!15!black}, font=\small},
ytick={0.5,0.55,0.6,0.65,0.7,0.75},
ylabel={composition \( \eta_4(\lambda, \mu)\)},
ticklabel style = {font=\scriptsize},
axis background/.style={fill=white},
xmajorgrids,
ymajorgrids,
colormap={mymap}{[1pt] rgb(0pt)=(0.18995,0.07176,0.23217); rgb(1pt)=(0.19483,0.08339,0.26149); rgb(2pt)=(0.19956,0.09498,0.29024); rgb(3pt)=(0.20415,0.10652,0.31844); rgb(4pt)=(0.2086,0.11802,0.34607); rgb(5pt)=(0.21291,0.12947,0.37314); rgb(6pt)=(0.21708,0.14087,0.39964); rgb(7pt)=(0.22111,0.15223,0.42558); rgb(8pt)=(0.225,0.16354,0.45096); rgb(9pt)=(0.22875,0.17481,0.47578); rgb(10pt)=(0.23236,0.18603,0.50004); rgb(11pt)=(0.23582,0.1972,0.52373); rgb(12pt)=(0.23915,0.20833,0.54686); rgb(13pt)=(0.24234,0.21941,0.56942); rgb(14pt)=(0.24539,0.23044,0.59142); rgb(15pt)=(0.2483,0.24143,0.61286); rgb(16pt)=(0.25107,0.25237,0.63374); rgb(17pt)=(0.25369,0.26327,0.65406); rgb(18pt)=(0.25618,0.27412,0.67381); rgb(19pt)=(0.25853,0.28492,0.693); rgb(20pt)=(0.26074,0.29568,0.71162); rgb(21pt)=(0.2628,0.30639,0.72968); rgb(22pt)=(0.26473,0.31706,0.74718); rgb(23pt)=(0.26652,0.32768,0.76412); rgb(24pt)=(0.26816,0.33825,0.7805); rgb(25pt)=(0.26967,0.34878,0.79631); rgb(26pt)=(0.27103,0.35926,0.81156); rgb(27pt)=(0.27226,0.3697,0.82624); rgb(28pt)=(0.27334,0.38008,0.84037); rgb(29pt)=(0.27429,0.39043,0.85393); rgb(30pt)=(0.27509,0.40072,0.86692); rgb(31pt)=(0.27576,0.41097,0.87936); rgb(32pt)=(0.27628,0.42118,0.89123); rgb(33pt)=(0.27667,0.43134,0.90254); rgb(34pt)=(0.27691,0.44145,0.91328); rgb(35pt)=(0.27701,0.45152,0.92347); rgb(36pt)=(0.27698,0.46153,0.93309); rgb(37pt)=(0.2768,0.47151,0.94214); rgb(38pt)=(0.27648,0.48144,0.95064); rgb(39pt)=(0.27603,0.49132,0.95857); rgb(40pt)=(0.27543,0.50115,0.96594); rgb(41pt)=(0.27469,0.51094,0.97275); rgb(42pt)=(0.27381,0.52069,0.97899); rgb(43pt)=(0.27273,0.5304,0.98461); rgb(44pt)=(0.27106,0.54015,0.9893); rgb(45pt)=(0.26878,0.54995,0.99303); rgb(46pt)=(0.26592,0.55979,0.99583); rgb(47pt)=(0.26252,0.56967,0.99773); rgb(48pt)=(0.25862,0.57958,0.99876); rgb(49pt)=(0.25425,0.5895,0.99896); rgb(50pt)=(0.24946,0.59943,0.99835); rgb(51pt)=(0.24427,0.60937,0.99697); rgb(52pt)=(0.23874,0.61931,0.99485); rgb(53pt)=(0.23288,0.62923,0.99202); rgb(54pt)=(0.22676,0.63913,0.98851); rgb(55pt)=(0.22039,0.64901,0.98436); rgb(56pt)=(0.21382,0.65886,0.97959); rgb(57pt)=(0.20708,0.66866,0.97423); rgb(58pt)=(0.20021,0.67842,0.96833); rgb(59pt)=(0.19326,0.68812,0.9619); rgb(60pt)=(0.18625,0.69775,0.95498); rgb(61pt)=(0.17923,0.70732,0.94761); rgb(62pt)=(0.17223,0.7168,0.93981); rgb(63pt)=(0.16529,0.7262,0.93161); rgb(64pt)=(0.15844,0.73551,0.92305); rgb(65pt)=(0.15173,0.74472,0.91416); rgb(66pt)=(0.14519,0.75381,0.90496); rgb(67pt)=(0.13886,0.76279,0.8955); rgb(68pt)=(0.13278,0.77165,0.8858); rgb(69pt)=(0.12698,0.78037,0.8759); rgb(70pt)=(0.12151,0.78896,0.86581); rgb(71pt)=(0.11639,0.7974,0.85559); rgb(72pt)=(0.11167,0.80569,0.84525); rgb(73pt)=(0.10738,0.81381,0.83484); rgb(74pt)=(0.10357,0.82177,0.82437); rgb(75pt)=(0.10026,0.82955,0.81389); rgb(76pt)=(0.0975,0.83714,0.80342); rgb(77pt)=(0.09532,0.84455,0.79299); rgb(78pt)=(0.09377,0.85175,0.78264); rgb(79pt)=(0.09287,0.85875,0.7724); rgb(80pt)=(0.09267,0.86554,0.7623); rgb(81pt)=(0.0932,0.87211,0.75237); rgb(82pt)=(0.09451,0.87844,0.74265); rgb(83pt)=(0.09662,0.88454,0.73316); rgb(84pt)=(0.09958,0.8904,0.72393); rgb(85pt)=(0.10342,0.896,0.715); rgb(86pt)=(0.10815,0.90142,0.70599); rgb(87pt)=(0.11374,0.90673,0.69651); rgb(88pt)=(0.12014,0.91193,0.6866); rgb(89pt)=(0.12733,0.91701,0.67627); rgb(90pt)=(0.13526,0.92197,0.66556); rgb(91pt)=(0.14391,0.9268,0.65448); rgb(92pt)=(0.15323,0.93151,0.64308); rgb(93pt)=(0.16319,0.93609,0.63137); rgb(94pt)=(0.17377,0.94053,0.61938); rgb(95pt)=(0.18491,0.94484,0.60713); rgb(96pt)=(0.19659,0.94901,0.59466); rgb(97pt)=(0.20877,0.95304,0.58199); rgb(98pt)=(0.22142,0.95692,0.56914); rgb(99pt)=(0.23449,0.96065,0.55614); rgb(100pt)=(0.24797,0.96423,0.54303); rgb(101pt)=(0.2618,0.96765,0.52981); rgb(102pt)=(0.27597,0.97092,0.51653); rgb(103pt)=(0.29042,0.97403,0.50321); rgb(104pt)=(0.30513,0.97697,0.48987); rgb(105pt)=(0.32006,0.97974,0.47654); rgb(106pt)=(0.33517,0.98234,0.46325); rgb(107pt)=(0.35043,0.98477,0.45002); rgb(108pt)=(0.36581,0.98702,0.43688); rgb(109pt)=(0.38127,0.98909,0.42386); rgb(110pt)=(0.39678,0.99098,0.41098); rgb(111pt)=(0.41229,0.99268,0.39826); rgb(112pt)=(0.42778,0.99419,0.38575); rgb(113pt)=(0.44321,0.99551,0.37345); rgb(114pt)=(0.45854,0.99663,0.3614); rgb(115pt)=(0.47375,0.99755,0.34963); rgb(116pt)=(0.48879,0.99828,0.33816); rgb(117pt)=(0.50362,0.99879,0.32701); rgb(118pt)=(0.51822,0.9991,0.31622); rgb(119pt)=(0.53255,0.99919,0.30581); rgb(120pt)=(0.54658,0.99907,0.29581); rgb(121pt)=(0.56026,0.99873,0.28623); rgb(122pt)=(0.57357,0.99817,0.27712); rgb(123pt)=(0.58646,0.99739,0.26849); rgb(124pt)=(0.59891,0.99638,0.26038); rgb(125pt)=(0.61088,0.99514,0.2528); rgb(126pt)=(0.62233,0.99366,0.24579); rgb(127pt)=(0.63323,0.99195,0.23937); rgb(128pt)=(0.64362,0.98999,0.23356); rgb(129pt)=(0.65394,0.98775,0.22835); rgb(130pt)=(0.66428,0.98524,0.2237); rgb(131pt)=(0.67462,0.98246,0.2196); rgb(132pt)=(0.68494,0.97941,0.21602); rgb(133pt)=(0.69525,0.9761,0.21294); rgb(134pt)=(0.70553,0.97255,0.21032); rgb(135pt)=(0.71577,0.96875,0.20815); rgb(136pt)=(0.72596,0.9647,0.2064); rgb(137pt)=(0.7361,0.96043,0.20504); rgb(138pt)=(0.74617,0.95593,0.20406); rgb(139pt)=(0.75617,0.95121,0.20343); rgb(140pt)=(0.76608,0.94627,0.20311); rgb(141pt)=(0.77591,0.94113,0.2031); rgb(142pt)=(0.78563,0.93579,0.20336); rgb(143pt)=(0.79524,0.93025,0.20386); rgb(144pt)=(0.80473,0.92452,0.20459); rgb(145pt)=(0.8141,0.91861,0.20552); rgb(146pt)=(0.82333,0.91253,0.20663); rgb(147pt)=(0.83241,0.90627,0.20788); rgb(148pt)=(0.84133,0.89986,0.20926); rgb(149pt)=(0.8501,0.89328,0.21074); rgb(150pt)=(0.85868,0.88655,0.2123); rgb(151pt)=(0.86709,0.87968,0.21391); rgb(152pt)=(0.8753,0.87267,0.21555); rgb(153pt)=(0.88331,0.86553,0.21719); rgb(154pt)=(0.89112,0.85826,0.2188); rgb(155pt)=(0.8987,0.85087,0.22038); rgb(156pt)=(0.90605,0.84337,0.22188); rgb(157pt)=(0.91317,0.83576,0.22328); rgb(158pt)=(0.92004,0.82806,0.22456); rgb(159pt)=(0.92666,0.82025,0.2257); rgb(160pt)=(0.93301,0.81236,0.22667); rgb(161pt)=(0.93909,0.80439,0.22744); rgb(162pt)=(0.94489,0.79634,0.228); rgb(163pt)=(0.95039,0.78823,0.22831); rgb(164pt)=(0.9556,0.78005,0.22836); rgb(165pt)=(0.96049,0.77181,0.22811); rgb(166pt)=(0.96507,0.76352,0.22754); rgb(167pt)=(0.96931,0.75519,0.22663); rgb(168pt)=(0.97323,0.74682,0.22536); rgb(169pt)=(0.97679,0.73842,0.22369); rgb(170pt)=(0.98,0.73,0.22161); rgb(171pt)=(0.98289,0.7214,0.21918); rgb(172pt)=(0.98549,0.7125,0.2165); rgb(173pt)=(0.98781,0.7033,0.21358); rgb(174pt)=(0.98986,0.69382,0.21043); rgb(175pt)=(0.99163,0.68408,0.20706); rgb(176pt)=(0.99314,0.67408,0.20348); rgb(177pt)=(0.99438,0.66386,0.19971); rgb(178pt)=(0.99535,0.65341,0.19577); rgb(179pt)=(0.99607,0.64277,0.19165); rgb(180pt)=(0.99654,0.63193,0.18738); rgb(181pt)=(0.99675,0.62093,0.18297); rgb(182pt)=(0.99672,0.60977,0.17842); rgb(183pt)=(0.99644,0.59846,0.17376); rgb(184pt)=(0.99593,0.58703,0.16899); rgb(185pt)=(0.99517,0.57549,0.16412); rgb(186pt)=(0.99419,0.56386,0.15918); rgb(187pt)=(0.99297,0.55214,0.15417); rgb(188pt)=(0.99153,0.54036,0.1491); rgb(189pt)=(0.98987,0.52854,0.14398); rgb(190pt)=(0.98799,0.51667,0.13883); rgb(191pt)=(0.9859,0.50479,0.13367); rgb(192pt)=(0.9836,0.49291,0.12849); rgb(193pt)=(0.98108,0.48104,0.12332); rgb(194pt)=(0.97837,0.4692,0.11817); rgb(195pt)=(0.97545,0.4574,0.11305); rgb(196pt)=(0.97234,0.44565,0.10797); rgb(197pt)=(0.96904,0.43399,0.10294); rgb(198pt)=(0.96555,0.42241,0.09798); rgb(199pt)=(0.96187,0.41093,0.0931); rgb(200pt)=(0.95801,0.39958,0.08831); rgb(201pt)=(0.95398,0.38836,0.08362); rgb(202pt)=(0.94977,0.37729,0.07905); rgb(203pt)=(0.94538,0.36638,0.07461); rgb(204pt)=(0.94084,0.35566,0.07031); rgb(205pt)=(0.93612,0.34513,0.06616); rgb(206pt)=(0.93125,0.33482,0.06218); rgb(207pt)=(0.92623,0.32473,0.05837); rgb(208pt)=(0.92105,0.31489,0.05475); rgb(209pt)=(0.91572,0.3053,0.05134); rgb(210pt)=(0.91024,0.29599,0.04814); rgb(211pt)=(0.90463,0.28696,0.04516); rgb(212pt)=(0.89888,0.27824,0.04243); rgb(213pt)=(0.89298,0.26981,0.03993); rgb(214pt)=(0.88691,0.26152,0.03753); rgb(215pt)=(0.88066,0.25334,0.03521); rgb(216pt)=(0.87422,0.24526,0.03297); rgb(217pt)=(0.8676,0.2373,0.03082); rgb(218pt)=(0.86079,0.22945,0.02875); rgb(219pt)=(0.8538,0.2217,0.02677); rgb(220pt)=(0.84662,0.21407,0.02487); rgb(221pt)=(0.83926,0.20654,0.02305); rgb(222pt)=(0.83172,0.19912,0.02131); rgb(223pt)=(0.82399,0.19182,0.01966); rgb(224pt)=(0.81608,0.18462,0.01809); rgb(225pt)=(0.80799,0.17753,0.0166); rgb(226pt)=(0.79971,0.17055,0.0152); rgb(227pt)=(0.79125,0.16368,0.01387); rgb(228pt)=(0.7826,0.15693,0.01264); rgb(229pt)=(0.77377,0.15028,0.01148); rgb(230pt)=(0.76476,0.14374,0.01041); rgb(231pt)=(0.75556,0.13731,0.00942); rgb(232pt)=(0.74617,0.13098,0.00851); rgb(233pt)=(0.73661,0.12477,0.00769); rgb(234pt)=(0.72686,0.11867,0.00695); rgb(235pt)=(0.71692,0.11268,0.00629); rgb(236pt)=(0.7068,0.1068,0.00571); rgb(237pt)=(0.6965,0.10102,0.00522); rgb(238pt)=(0.68602,0.09536,0.00481); rgb(239pt)=(0.67535,0.0898,0.00449); rgb(240pt)=(0.66449,0.08436,0.00424); rgb(241pt)=(0.65345,0.07902,0.00408); rgb(242pt)=(0.64223,0.0738,0.00401); rgb(243pt)=(0.63082,0.06868,0.00401); rgb(244pt)=(0.61923,0.06367,0.0041); rgb(245pt)=(0.60746,0.05878,0.00427); rgb(246pt)=(0.5955,0.05399,0.00453); rgb(247pt)=(0.58336,0.04931,0.00486); rgb(248pt)=(0.57103,0.04474,0.00529); rgb(249pt)=(0.55852,0.04028,0.00579); rgb(250pt)=(0.54583,0.03593,0.00638); rgb(251pt)=(0.53295,0.03169,0.00705); rgb(252pt)=(0.51989,0.02756,0.0078); rgb(253pt)=(0.50664,0.02354,0.00863); rgb(254pt)=(0.49321,0.01963,0.00955); rgb(255pt)=(0.4796,0.01583,0.01055)},
colorbar,
colorbar style={ylabel style={font=\color{white!15!black}, font=\small}, ylabel={parameter \(\lambda\)}}
]
\addplot [color=mycolor1, line width=1.0pt, forget plot]
  table[row sep=crcr]{%
0	0.476062281231506\\
0.02	0.475949739412606\\
0.04	0.476661459047031\\
0.06	0.47806983514205\\
0.08	0.480072309775905\\
0.1	0.482585514073741\\
0.12	0.485540980809584\\
0.14	0.488881957355459\\
0.16	0.492561003024984\\
0.18	0.496538154455982\\
0.2	0.500779508292561\\
0.22	0.505256114473297\\
0.24	0.509943103509495\\
0.26	0.514818992000703\\
0.28	0.519865125316853\\
0.3	0.52506522684747\\
0.32	0.530405030778862\\
0.34	0.535871980882107\\
0.36	0.541454981870982\\
0.38	0.547144192928314\\
0.4	0.552930855286631\\
0.42	0.558807147485625\\
0.44	0.564766063258356\\
0.46	0.570801308023651\\
0.48	0.576907210759126\\
0.5	0.583078648652898\\
0.52	0.589310982423214\\
0.54	0.595600000584512\\
0.56	0.601941871248763\\
0.58	0.608333100299697\\
0.6	0.614770494978062\\
0.62	0.621251132078432\\
0.64	0.627772330090323\\
0.66	0.634331624724454\\
0.68	0.640926747353801\\
0.7	0.647555605972309\\
0.72	0.654216268334808\\
0.74	0.660906946992083\\
0.76	0.667625985977144\\
0.78	0.674371848933974\\
0.8	0.681143108509666\\
0.82	0.687938436855808\\
0.84	0.694756597106113\\
0.86	0.70159643571521\\
0.88	0.708456875558754\\
0.9	0.715336909708029\\
0.92	0.722235595803348\\
0.94	0.729152050960133\\
0.96	0.736085447149728\\
0.98	0.743035007004142\\
1	0.75\\
};
\addplot [color=mycolor2, line width=1.0pt, forget plot]
  table[row sep=crcr]{%
0	0.485800592544729\\
0.02	0.48539127874882\\
0.04	0.485797621838295\\
0.06	0.486895962820553\\
0.08	0.48858643694571\\
0.1	0.490787513318953\\
0.12	0.493431972210373\\
0.14	0.496463896997699\\
0.16	0.499836394215592\\
0.18	0.503509844071676\\
0.2	0.507450542812776\\
0.22	0.511629638231607\\
0.24	0.516022287034949\\
0.26	0.520606981937279\\
0.28	0.525365009890922\\
0.3	0.530280012575361\\
0.32	0.535337627314121\\
0.34	0.540525191756773\\
0.36	0.545831499495312\\
0.38	0.551246596652361\\
0.4	0.556761611645109\\
0.42	0.562368611979391\\
0.44	0.568060483195838\\
0.46	0.573830826070826\\
0.48	0.57967386893936\\
0.5	0.585584392606796\\
0.52	0.591557665789886\\
0.54	0.597589389403934\\
0.56	0.603675648313485\\
0.58	0.609812869405534\\
0.6	0.615997785039327\\
0.62	0.622227401085219\\
0.64	0.628498968894123\\
0.66	0.634809960644959\\
0.68	0.641158047604499\\
0.7	0.647541080905959\\
0.72	0.65395707451234\\
0.74	0.660404190080181\\
0.76	0.666880723480939\\
0.78	0.673385092771983\\
0.8	0.679915827438533\\
0.82	0.686471558752573\\
0.84	0.693051011115762\\
0.86	0.699652994271113\\
0.88	0.706276396283409\\
0.9	0.712920177201266\\
0.92	0.719583363324847\\
0.94	0.726265042012759\\
0.96	0.732964356969888\\
0.98	0.739680503964986\\
1	0.74641272693298\\
};
\addplot [color=mycolor3, line width=1.0pt, forget plot]
  table[row sep=crcr]{%
0	0.495904823789265\\
0.02	0.495197859205416\\
0.04	0.495296194742843\\
0.06	0.496080432852309\\
0.08	0.49745366766051\\
0.1	0.499336427043424\\
0.12	0.501662920865058\\
0.14	0.504378218096557\\
0.16	0.507436095141909\\
0.18	0.510797376273252\\
0.2	0.514428639688424\\
0.22	0.518301198543047\\
0.24	0.522390291115449\\
0.26	0.526674431684575\\
0.28	0.531134886102837\\
0.3	0.535755244984969\\
0.32	0.540521073952322\\
0.34	0.545419625176703\\
0.36	0.550439598045316\\
0.38	0.555570939457179\\
0.4	0.560804676300304\\
0.42	0.56613277421791\\
0.44	0.571548017973268\\
0.46	0.577043909655411\\
0.48	0.582614581697002\\
0.5	0.588254722249346\\
0.52	0.593959510913814\\
0.54	0.599724563190786\\
0.56	0.605545882297069\\
0.58	0.611419817236181\\
0.6	0.617343026194835\\
0.62	0.623312444492664\\
0.64	0.6293252564378\\
0.66	0.635378870544054\\
0.68	0.641470897650396\\
0.7	0.647599131553798\\
0.72	0.653761531824944\\
0.74	0.65995620852506\\
0.76	0.666181408582944\\
0.78	0.672435503625508\\
0.8	0.678716979084092\\
0.82	0.685024424423172\\
0.84	0.691356524358831\\
0.86	0.697712050951978\\
0.88	0.704089856476302\\
0.9	0.710488866973831\\
0.92	0.716908076421976\\
0.94	0.723346541445421\\
0.96	0.729803376514394\\
0.98	0.736277749577918\\
1	0.74276887808675\\
};
\addplot [color=mycolor4, line width=1.0pt, forget plot]
  table[row sep=crcr]{%
0	0.506392916046655\\
0.02	0.505388387965906\\
0.04	0.505176897872502\\
0.06	0.50564363259809\\
0.08	0.506694921634767\\
0.1	0.508253584591029\\
0.12	0.510255455433839\\
0.14	0.512646750306253\\
0.16	0.515382049345976\\
0.18	0.518422731661044\\
0.2	0.521735749035648\\
0.22	0.525292655811166\\
0.24	0.529068834608325\\
0.26	0.533042873267554\\
0.28	0.537196059638268\\
0.3	0.54151196900614\\
0.32	0.545976124928375\\
0.34	0.550575718677681\\
0.36	0.555299375809991\\
0.38	0.560136960873028\\
0.4	0.565079413177936\\
0.42	0.570118608018465\\
0.44	0.575247238853159\\
0.46	0.580458716847068\\
0.48	0.585747084860378\\
0.5	0.591106943516856\\
0.52	0.596533387418045\\
0.54	0.602021949915186\\
0.56	0.607568555128655\\
0.58	0.613169476129087\\
0.6	0.618821298376339\\
0.62	0.624520887660904\\
0.64	0.630265361913892\\
0.66	0.63605206635169\\
0.68	0.641878551503972\\
0.7	0.647742553742254\\
0.72	0.653641977983137\\
0.74	0.659574882288085\\
0.76	0.665539464121454\\
0.78	0.67153404806216\\
0.8	0.677557074792671\\
0.82	0.68360709121308\\
0.84	0.689682741548373\\
0.86	0.695782759334408\\
0.88	0.701905960182937\\
0.9	0.708051235238737\\
0.92	0.714217545252808\\
0.94	0.72040391520502\\
0.96	0.72660942941768\\
0.98	0.73283322710851\\
1	0.739074498337631\\
};
\addplot [color=mycolor5, line width=1.0pt, forget plot]
  table[row sep=crcr]{%
0	0.517283606260942\\
0.02	0.515982755425962\\
0.04	0.515460613423526\\
0.06	0.515607279547196\\
0.08	0.516332604572047\\
0.1	0.517561944348182\\
0.12	0.519232963761254\\
0.14	0.521293199946412\\
0.16	0.523698182339757\\
0.18	0.526409966561655\\
0.2	0.529395979613812\\
0.22	0.532628101904746\\
0.24	0.53608193131252\\
0.26	0.53973618851714\\
0.28	0.543572232945799\\
0.3	0.547573666048742\\
0.32	0.551726004060664\\
0.34	0.5560164064517\\
0.36	0.560433449316156\\
0.38	0.564966935255811\\
0.4	0.569607733080231\\
0.42	0.574347642007309\\
0.44	0.579179276103825\\
0.46	0.584095965531932\\
0.48	0.589091671817576\\
0.5	0.594160914871776\\
0.52	0.599298709905815\\
0.54	0.604500512710033\\
0.56	0.609762172030563\\
0.58	0.615079887992599\\
0.6	0.620450175693117\\
0.62	0.625869833228396\\
0.64	0.631335913538636\\
0.66	0.636845699548388\\
0.68	0.642396682161307\\
0.7	0.647986540734047\\
0.72	0.653613125709478\\
0.74	0.659274443135637\\
0.76	0.664968640835812\\
0.78	0.670693996027871\\
0.8	0.676448904218711\\
0.82	0.682231869223185\\
0.84	0.688041494176897\\
0.86	0.693876473429264\\
0.88	0.699735585217902\\
0.9	0.705617685037847\\
0.92	0.711521699629927\\
0.94	0.717446621521885\\
0.96	0.723391504063821\\
0.98	0.72935545690654\\
1	0.735337641877362\\
};
\addplot [color=mycolor6, line width=1.0pt, forget plot]
  table[row sep=crcr]{%
0	0.52859639021015\\
0.02	0.527001831017375\\
0.04	0.526169415194648\\
0.06	0.525994484688079\\
0.08	0.526390704997731\\
0.1	0.527286223225491\\
0.12	0.528620753881537\\
0.14	0.53034334118647\\
0.16	0.532410621479655\\
0.18	0.534785460018819\\
0.2	0.537435871320231\\
0.22	0.540334156540146\\
0.24	0.543456208644569\\
0.26	0.546780948486934\\
0.28	0.55028986389574\\
0.3	0.553966630468115\\
0.32	0.557796797656251\\
0.34	0.56176752739654\\
0.36	0.56586737529955\\
0.38	0.570086106528924\\
0.4	0.574414540118406\\
0.42	0.57884441673117\\
0.44	0.583368285844085\\
0.46	0.587979409107572\\
0.48	0.592671677238382\\
0.5	0.597439538284809\\
0.52	0.602277935489242\\
0.54	0.607182253282703\\
0.56	0.612148270196277\\
0.58	0.61717211767746\\
0.6	0.622250243965242\\
0.62	0.627379382313515\\
0.64	0.632556522964168\\
0.66	0.637778888363633\\
0.68	0.643043911193233\\
0.7	0.64834921484756\\
0.72	0.653692596048435\\
0.74	0.659072009326776\\
0.76	0.664485553142357\\
0.78	0.669931457443257\\
0.8	0.675408072493754\\
0.82	0.680913858822283\\
0.84	0.686447378160618\\
0.86	0.692007285262078\\
0.88	0.697592320500853\\
0.9	0.703201303166807\\
0.92	0.708833125380693\\
0.94	0.714486746563809\\
0.96	0.720161188404051\\
0.98	0.725855530267141\\
1	0.731568905007764\\
};
\addplot [color=mycolor7, line width=1.0pt, forget plot]
  table[row sep=crcr]{%
0	0.540351456104047\\
0.02	0.538467435720841\\
0.04	0.537326581405393\\
0.06	0.536829806206734\\
0.08	0.536894888735087\\
0.1	0.537453031269878\\
0.12	0.538446227711282\\
0.14	0.539825227457783\\
0.16	0.541547943521815\\
0.18	0.54357819577399\\
0.2	0.545884709795012\\
0.22	0.548440312672034\\
0.24	0.55122128197461\\
0.26	0.554206814918135\\
0.28	0.557378592601304\\
0.3	0.560720420027419\\
0.32	0.56421792696646\\
0.34	0.567858317989929\\
0.36	0.571630162499678\\
0.38	0.575523217479159\\
0.4	0.579528277168242\\
0.42	0.583637045008201\\
0.44	0.587842024100495\\
0.46	0.592136423130024\\
0.48	0.596514075264381\\
0.5	0.600969367987914\\
0.52	0.605497182188366\\
0.54	0.610092839103227\\
0.56	0.614752053967541\\
0.58	0.619470895395956\\
0.6	0.624245749688114\\
0.62	0.629073289374961\\
0.64	0.633950445429521\\
0.66	0.638874382653568\\
0.68	0.64384247782464\\
0.7	0.648852300248875\\
0.72	0.653901594416221\\
0.74	0.658988264497562\\
0.76	0.66411036045952\\
0.78	0.669266065603371\\
0.8	0.674453685360531\\
0.82	0.679671637199229\\
0.84	0.684918441515881\\
0.86	0.690192713400895\\
0.88	0.695493155182497\\
0.9	0.700818549664159\\
0.92	0.706167753981481\\
0.94	0.711539694013357\\
0.96	0.716933359289897\\
0.98	0.722347798346395\\
1	0.727782114478367\\
};
\addplot [color=mycolor8, line width=1.0pt, forget plot]
  table[row sep=crcr]{%
0	0.552569578474767\\
0.02	0.550402282146847\\
0.04	0.548956583097977\\
0.06	0.548139287372922\\
0.08	0.547872586610676\\
0.1	0.548091008842495\\
0.12	0.548739066811411\\
0.14	0.549769424494455\\
0.16	0.55114145332734\\
0.18	0.55282008480211\\
0.2	0.554774890807681\\
0.22	0.556979340664381\\
0.24	0.559410196475168\\
0.26	0.562047017655732\\
0.28	0.564871752315786\\
0.3	0.567868398234152\\
0.32	0.57102271998103\\
0.34	0.574322011630291\\
0.36	0.577754896713944\\
0.38	0.581311158773257\\
0.4	0.584981597182599\\
0.42	0.588757903955106\\
0.44	0.592632558052135\\
0.46	0.59659873436216\\
0.48	0.600650225027403\\
0.5	0.604781371207161\\
0.52	0.608987003697499\\
0.54	0.613262391094676\\
0.56	0.617603194407386\\
0.58	0.622005427200849\\
0.6	0.626465420501811\\
0.62	0.630979791813924\\
0.64	0.635545417692608\\
0.66	0.640159409411291\\
0.68	0.64481909131999\\
0.7	0.649521981554994\\
0.72	0.654265774806945\\
0.74	0.659048326895558\\
0.76	0.663867640933782\\
0.78	0.668721854893565\\
0.8	0.673609230410306\\
0.82	0.678528142684405\\
0.84	0.683477071356476\\
0.86	0.688454592248451\\
0.88	0.693459369876188\\
0.9	0.698490150650806\\
0.92	0.703545756695936\\
0.94	0.708625080216775\\
0.96	0.71372707836431\\
0.98	0.718850768544668\\
1	0.723995224129201\\
};
\addplot [color=mycolor9, line width=1.0pt, forget plot]
  table[row sep=crcr]{%
0	0.565271958578599\\
0.02	0.56282986987867\\
0.04	0.561085037545126\\
0.06	0.559950469999106\\
0.08	0.559353068944968\\
0.1	0.559230962267924\\
0.12	0.559531428620677\\
0.14	0.560209266017603\\
0.16	0.561225497430671\\
0.18	0.562546334872849\\
0.2	0.564142343686675\\
0.22	0.565987763305493\\
0.24	0.56805995134205\\
0.26	0.570338925647559\\
0.28	0.572806984772667\\
0.3	0.575448391604813\\
0.32	0.578249108244418\\
0.34	0.581196572692503\\
0.36	0.584279509853886\\
0.38	0.58748777085759\\
0.4	0.590812195865174\\
0.42	0.594244496456606\\
0.44	0.597777154409997\\
0.46	0.6014033342698\\
0.48	0.60511680756076\\
0.5	0.608911886877211\\
0.52	0.612783368378328\\
0.54	0.616726481464655\\
0.56	0.62073684461095\\
0.58	0.624810426494249\\
0.6	0.628943511691019\\
0.62	0.633132670328855\\
0.64	0.637374731170874\\
0.66	0.641666757688213\\
0.68	0.646006026740629\\
0.7	0.650390009539465\\
0.72	0.654816354612906\\
0.74	0.659282872532066\\
0.76	0.663787522189154\\
0.78	0.668328398446781\\
0.8	0.672903721001179\\
0.82	0.677511824322378\\
0.84	0.682151148551773\\
0.86	0.686820231252431\\
0.88	0.691517699920379\\
0.9	0.696242265176206\\
0.92	0.700992714565963\\
0.94	0.705767906908675\\
0.96	0.710566767135046\\
0.98	0.715388281568273\\
1	0.720231493603411\\
};
\addplot [color=mycolor10, line width=1.0pt, forget plot]
  table[row sep=crcr]{%
0	0.578479992865184\\
0.02	0.575774319433331\\
0.04	0.573738612152969\\
0.06	0.572292371412056\\
0.08	0.571367497509637\\
0.1	0.570905991634467\\
0.12	0.570858150058225\\
0.14	0.571181132359181\\
0.16	0.571837816308796\\
0.18	0.572795874623594\\
0.2	0.574027025000418\\
0.22	0.575506416642074\\
0.24	0.577212125145779\\
0.26	0.579124734061885\\
0.28	0.581226986255244\\
0.3	0.583503491851812\\
0.32	0.58594048233805\\
0.34	0.588525602522437\\
0.36	0.591247733727912\\
0.38	0.594096842878906\\
0.4	0.597063853163926\\
0.42	0.600140532758735\\
0.44	0.60331939873476\\
0.46	0.606593633788879\\
0.48	0.60995701384216\\
0.5	0.613403844887633\\
0.52	0.616928907737408\\
0.54	0.620527409539885\\
0.56	0.624194941118655\\
0.58	0.627927439333527\\
0.6	0.631721153787196\\
0.62	0.635572617303262\\
0.64	0.639478619686433\\
0.66	0.643436184346915\\
0.68	0.647442547430785\\
0.7	0.651495139148398\\
0.72	0.65559156703542\\
0.74	0.659729600917051\\
0.76	0.663907159376638\\
0.78	0.668122297555943\\
0.8	0.672373196136647\\
0.82	0.676658151371771\\
0.84	0.680975566052117\\
0.86	0.685323941306996\\
0.88	0.689701869150703\\
0.9	0.694108025696794\\
0.92	0.698541164971392\\
0.94	0.703000113264706\\
0.96	0.707483763966942\\
0.98	0.711991072840802\\
1	0.716521053688127\\
};
\addplot [color=mycolor11, line width=1.0pt, forget plot]
  table[row sep=crcr]{%
0	0.592214944729106\\
0.02	0.58926012225442\\
0.04	0.586944858913114\\
0.06	0.58519540806206\\
0.08	0.583948941551897\\
0.1	0.583151601189401\\
0.12	0.58275695310733\\
0.14	0.582724751083679\\
0.16	0.583019939086176\\
0.18	0.5836118407395\\
0.2	0.58447349608556\\
0.22	0.585581115327071\\
0.24	0.586913626169225\\
0.26	0.588452296568502\\
0.28	0.590180418629007\\
0.3	0.59208304238739\\
0.32	0.594146750535942\\
0.34	0.59635946692274\\
0.36	0.598710293064442\\
0.38	0.601189368004871\\
0.4	0.603787747720411\\
0.42	0.606497300963727\\
0.44	0.609310618989684\\
0.46	0.612220937051773\\
0.48	0.6152220659166\\
0.5	0.618308331936021\\
0.52	0.621474524454802\\
0.54	0.624715849527205\\
0.56	0.628027889076852\\
0.58	0.631406564767346\\
0.6	0.634848105961708\\
0.62	0.638349021240727\\
0.64	0.641906073027442\\
0.66	0.645516254929592\\
0.68	0.649176771466359\\
0.7	0.652885019891718\\
0.72	0.656638573865726\\
0.74	0.66043516875822\\
0.76	0.664272688397638\\
0.78	0.668149153101861\\
0.8	0.672062708848671\\
0.82	0.676011617461195\\
0.84	0.679994247699078\\
0.86	0.684009067159337\\
0.88	0.688054634902322\\
0.9	0.692129594728169\\
0.92	0.696232669037761\\
0.94	0.700362653219764\\
0.96	0.704518410511883\\
0.98	0.708698867290247\\
1	0.71290300874587\\
};
\addplot [color=mycolor12, line width=1.0pt, forget plot]
  table[row sep=crcr]{%
0	0.606497486076447\\
0.02	0.603311775960689\\
0.04	0.600731951903582\\
0.06	0.598691242129585\\
0.08	0.597132338669452\\
0.1	0.596005778062901\\
0.12	0.595268643542285\\
0.14	0.594883516620839\\
0.16	0.594817623904051\\
0.18	0.595042137989855\\
0.2	0.595531600939241\\
0.22	0.596263445957656\\
0.24	0.597217598307691\\
0.26	0.598376140552719\\
0.28	0.59972303034987\\
0.3	0.601243861413774\\
0.32	0.602925660137894\\
0.34	0.604756711818331\\
0.36	0.606726411572299\\
0.38	0.608825135951854\\
0.4	0.611044131976966\\
0.42	0.613375420891571\\
0.44	0.615811714412888\\
0.46	0.618346341622074\\
0.48	0.620973184951475\\
0.5	0.62368662397483\\
0.52	0.626481485912842\\
0.54	0.629353001936392\\
0.56	0.632296768490254\\
0.58	0.635308712977009\\
0.6	0.638385063238241\\
0.62	0.641522320351675\\
0.64	0.64471723433141\\
0.66	0.647966782376092\\
0.68	0.651268149358697\\
0.7	0.65461871029298\\
0.72	0.658016014546799\\
0.74	0.661457771602615\\
0.76	0.664941838191078\\
0.78	0.668466206645695\\
0.8	0.672028994345466\\
0.82	0.675628434128737\\
0.84	0.679262865575594\\
0.86	0.682930727068363\\
0.88	0.686630548550378\\
0.9	0.690360944912397\\
0.92	0.694120609944105\\
0.94	0.697908310795173\\
0.96	0.701722882896464\\
0.98	0.705563225297421\\
1	0.709428296380366\\
};
\addplot [color=mycolor13, line width=1.0pt, forget plot]
  table[row sep=crcr]{%
0	0.621347063264531\\
0.02	0.617953262707027\\
0.04	0.615128289768745\\
0.06	0.612812517927608\\
0.08	0.610954372987575\\
0.1	0.609509018178177\\
0.12	0.608437288764522\\
0.14	0.607704822623954\\
0.16	0.607281345986717\\
0.18	0.607140082966755\\
0.2	0.607257264555226\\
0.22	0.607611718063974\\
0.24	0.608184522049327\\
0.26	0.608958714847564\\
0.28	0.609919047249909\\
0.3	0.611051771710446\\
0.32	0.612344461942476\\
0.34	0.613785857912197\\
0.36	0.61536573215398\\
0.38	0.61707477406226\\
0.4	0.618904489401673\\
0.42	0.62084711275037\\
0.44	0.622895530975372\\
0.46	0.625043216151614\\
0.48	0.627284166592424\\
0.5	0.629612854869745\\
0.52	0.632024181876227\\
0.54	0.634513436125436\\
0.56	0.637076257606359\\
0.58	0.639708605608515\\
0.6	0.642406730017962\\
0.62	0.645167145655072\\
0.64	0.647986609284562\\
0.66	0.65086209897865\\
0.68	0.653790795557073\\
0.7	0.656770065864131\\
0.72	0.659797447674075\\
0.74	0.662870636042817\\
0.76	0.665987470946808\\
0.78	0.669145926069673\\
0.8	0.672344098614143\\
0.82	0.675580200031572\\
0.84	0.678852547574034\\
0.86	0.682159556585092\\
0.88	0.685499733454974\\
0.9	0.688871669174294\\
0.92	0.692274033427808\\
0.94	0.695705569176175\\
0.96	0.699165087679284\\
0.98	0.702651463919761\\
1	0.706163632389593\\
};
\addplot [color=mycolor14, line width=1.0pt, forget plot]
  table[row sep=crcr]{%
0	0.636781025328072\\
0.02	0.633207312615914\\
0.04	0.630161910211746\\
0.06	0.627592441290669\\
0.08	0.62545323108313\\
0.1	0.623704268123702\\
0.12	0.622310352886474\\
0.14	0.621240395520373\\
0.16	0.620466833084479\\
0.18	0.61996514321392\\
0.2	0.619713436100109\\
0.22	0.619692110457093\\
0.24	0.619883562065427\\
0.26	0.62027193575196\\
0.28	0.620842913436414\\
0.3	0.621583532270617\\
0.32	0.622482028001059\\
0.34	0.623527699565541\\
0.36	0.624710791639841\\
0.38	0.626022392418212\\
0.4	0.627454344371354\\
0.42	0.628999166099502\\
0.44	0.630649983703986\\
0.46	0.632400470351574\\
0.48	0.634244792912813\\
0.5	0.636177564726826\\
0.52	0.638193803687379\\
0.54	0.640288894963714\\
0.56	0.64245855776906\\
0.58	0.644698815673195\\
0.6	0.647005970025802\\
0.62	0.649376576116858\\
0.64	0.651807421750746\\
0.66	0.654295507953704\\
0.68	0.656838031570848\\
0.7	0.659432369540312\\
0.72	0.662076064658933\\
0.74	0.664766812676996\\
0.76	0.667502450579491\\
0.78	0.670280945928522\\
0.8	0.673100387156437\\
0.82	0.675958974712184\\
0.84	0.678855012974655\\
0.86	0.681786902856606\\
0.88	0.684753135031317\\
0.9	0.687752283721663\\
0.92	0.690783000997863\\
0.94	0.69384401153595\\
0.96	0.696934107794114\\
0.98	0.700052145568566\\
1	0.703197039894523\\
};
\addplot [color=mycolor15, line width=1.0pt, forget plot]
  table[row sep=crcr]{%
0	0.652813429067417\\
0.02	0.649094371814432\\
0.04	0.64585964236658\\
0.06	0.643064135627612\\
0.08	0.640668178676194\\
0.1	0.638636736896711\\
0.12	0.636938755377633\\
0.14	0.635546609430311\\
0.16	0.634435643702139\\
0.18	0.63358378365514\\
0.2	0.632971206491969\\
0.22	0.63258006118121\\
0.24	0.632394229242542\\
0.26	0.632399119532783\\
0.28	0.632581491525117\\
0.3	0.632929302570353\\
0.32	0.633431575427275\\
0.34	0.634078282991903\\
0.36	0.634860247675687\\
0.38	0.635769053305789\\
0.4	0.636796967766323\\
0.42	0.637936874883227\\
0.44	0.639182214289298\\
0.46	0.640526928199534\\
0.48	0.641965414187761\\
0.5	0.643492483189657\\
0.52	0.645103322069561\\
0.54	0.646793460182776\\
0.56	0.648558739444493\\
0.58	0.650395287483699\\
0.6	0.652299493517362\\
0.62	0.654267986628643\\
0.64	0.656297616174242\\
0.66	0.658385434081305\\
0.68	0.660528678824646\\
0.7	0.662724760901122\\
0.72	0.664971249640443\\
0.74	0.667265861211144\\
0.76	0.669606447697235\\
0.78	0.671990987135664\\
0.8	0.674417574417418\\
0.82	0.676884412966164\\
0.84	0.67938980711799\\
0.86	0.681932155134297\\
0.88	0.684509942787304\\
0.9	0.68712173746414\\
0.92	0.689766182741289\\
0.94	0.692441993386154\\
0.96	0.695147950747039\\
0.98	0.697882898496786\\
1	0.700645738698803\\
};
\addplot [color=mycolor16, line width=1.0pt, forget plot]
  table[row sep=crcr]{%
0	0.669453402545614\\
0.02	0.665631162694723\\
0.04	0.662245892571952\\
0.06	0.659259680087635\\
0.08	0.65663887561911\\
0.1	0.654353509349926\\
0.12	0.652376801642255\\
0.14	0.650684749561414\\
0.16	0.649255776083622\\
0.18	0.648070431173606\\
0.2	0.647111135999231\\
0.22	0.646361963191186\\
0.24	0.645808447357711\\
0.26	0.64543742110342\\
0.28	0.6452368726352\\
0.3	0.645195821710926\\
0.32	0.64530421123226\\
0.34	0.64555281222722\\
0.36	0.645933140331987\\
0.38	0.646437382180443\\
0.4	0.647058330356779\\
0.42	0.647789325771119\\
0.44	0.648624206488333\\
0.46	0.649557262182359\\
0.48	0.650583193507514\\
0.5	0.651697075778432\\
0.52	0.652894326434817\\
0.54	0.654170675838725\\
0.56	0.655522141012817\\
0.58	0.656945001979748\\
0.6	0.658435780407007\\
0.62	0.659991220299296\\
0.64	0.661608270513037\\
0.66	0.663284068895471\\
0.68	0.665015927874959\\
0.7	0.666801321349874\\
0.72	0.668637872741596\\
0.74	0.670523344092781\\
0.76	0.67245562610572\\
0.78	0.674432729027584\\
0.8	0.676452774299702\\
0.82	0.678513986897195\\
0.84	0.680614688293297\\
0.86	0.682753289989703\\
0.88	0.684928287560546\\
0.9	0.687138255163034\\
0.92	0.689381840472642\\
0.94	0.691657760005051\\
0.96	0.693964794790792\\
0.98	0.696301786371974\\
1	0.698667633093448\\
};
\addplot [color=mycolor17, line width=1.0pt, forget plot]
  table[row sep=crcr]{%
0	0.686702901316497\\
0.02	0.682828678125546\\
0.04	0.679340915133823\\
0.06	0.676208694064777\\
0.08	0.673404309078145\\
0.1	0.670902860114214\\
0.12	0.668681906486003\\
0.14	0.666721170556649\\
0.16	0.665002283252313\\
0.18	0.663508564678188\\
0.2	0.662224834315542\\
0.22	0.661137246248702\\
0.24	0.660233145653886\\
0.26	0.659500943416427\\
0.28	0.658930006259844\\
0.3	0.658510560193071\\
0.32	0.658233605429716\\
0.34	0.658090841219965\\
0.36	0.658074599273427\\
0.38	0.65817778464882\\
0.4	0.658393823151424\\
0.42	0.658716614417479\\
0.44	0.659140489980927\\
0.46	0.659660175715932\\
0.48	0.66027075813157\\
0.5	0.660967654065469\\
0.52	0.661746583383129\\
0.54	0.662603544340801\\
0.56	0.663534791313597\\
0.58	0.664536814628084\\
0.6	0.665606322270976\\
0.62	0.666740223273403\\
0.64	0.667935612594387\\
0.66	0.669189757348046\\
0.68	0.6705000842372\\
0.7	0.671864168071882\\
0.72	0.673279721265012\\
0.74	0.674744584209579\\
0.76	0.676256716452218\\
0.78	0.677814188587321\\
0.8	0.679415174803988\\
0.82	0.681057946025296\\
0.84	0.682740863585664\\
0.86	0.684462373397733\\
0.88	0.686221000565097\\
0.9	0.688015344401626\\
0.92	0.689844073822034\\
0.94	0.691705923071794\\
0.96	0.693599687767612\\
0.98	0.695524221222429\\
1	0.697478431031364\\
};
\addplot [color=mycolor18, line width=1.0pt, forget plot]
  table[row sep=crcr]{%
0	0.704553623428966\\
0.02	0.70068938463252\\
0.04	0.697158355359085\\
0.06	0.693936271417343\\
0.08	0.69100116866537\\
0.1	0.688333116786961\\
0.12	0.685913989208812\\
0.14	0.683727263561224\\
0.16	0.681757848054309\\
0.18	0.679991929926089\\
0.2	0.67841684275689\\
0.22	0.677020949965658\\
0.24	0.675793542231673\\
0.26	0.674724746937686\\
0.28	0.673805448022363\\
0.3	0.673027214872307\\
0.32	0.672382239086165\\
0.34	0.671863278112471\\
0.36	0.671463604905064\\
0.38	0.671176962859612\\
0.4	0.670997525396039\\
0.42	0.670919859637492\\
0.44	0.670938893709482\\
0.46	0.671049887245112\\
0.48	0.671248404735563\\
0.5	0.671530291410697\\
0.52	0.671891651373883\\
0.54	0.672328827749033\\
0.56	0.672838384627055\\
0.58	0.673417090624298\\
0.6	0.674061903887502\\
0.62	0.674769958398953\\
0.64	0.67553855145217\\
0.66	0.676365132183046\\
0.68	0.677247291054112\\
0.7	0.678182750200772\\
0.72	0.679169354558196\\
0.74	0.680205063696182\\
0.76	0.68128794429698\\
0.78	0.682416163217744\\
0.8	0.683587981085297\\
0.82	0.684801746376158\\
0.84	0.686055889939471\\
0.86	0.687348919924652\\
0.88	0.688679417079295\\
0.9	0.690046030386179\\
0.92	0.691447473011188\\
0.94	0.692882518536592\\
0.96	0.694349997456523\\
0.98	0.695848793913595\\
1	0.697377842657525\\
};
\addplot [color=mycolor19, line width=1.0pt, forget plot]
  table[row sep=crcr]{%
0	0.722982749332744\\
0.02	0.719203311414488\\
0.04	0.715701756897534\\
0.06	0.712459976903628\\
0.08	0.709461402574358\\
0.1	0.706690844770779\\
0.12	0.704134353387471\\
0.14	0.701779093536948\\
0.16	0.699613236290496\\
0.18	0.69762586201725\\
0.2	0.695806874658698\\
0.22	0.694146925522005\\
0.24	0.692637345381424\\
0.26	0.691270083849971\\
0.28	0.690037655129046\\
0.3	0.68893308936675\\
0.32	0.687949888959835\\
0.34	0.687081989222883\\
0.36	0.68632372292384\\
0.38	0.685669788249643\\
0.4	0.685115219821015\\
0.42	0.684655362423122\\
0.44	0.684285847159775\\
0.46	0.684002569774262\\
0.48	0.683801670910575\\
0.5	0.683679518115389\\
0.52	0.683632689404275\\
0.54	0.683657958235819\\
0.56	0.683752279754893\\
0.58	0.683912778181764\\
0.6	0.684136735237201\\
0.62	0.68442157950565\\
0.64	0.684764876648942\\
0.66	0.68516432039225\\
0.68	0.685617724212111\\
0.7	0.686123013663546\\
0.72	0.686678219289665\\
0.74	0.687281470062814\\
0.76	0.687930987311359\\
0.78	0.688625079090649\\
0.8	0.689362134960734\\
0.82	0.690140621136937\\
0.84	0.690959075982599\\
0.86	0.691816105816159\\
0.88	0.692710381007275\\
0.9	0.693640632339017\\
0.92	0.694605647615207\\
0.94	0.695604268493844\\
0.96	0.696635387529235\\
0.98	0.697697945406963\\
1	0.698790928357158\\
};
\addplot [color=mycolor20, line width=1.0pt, forget plot]
  table[row sep=crcr]{%
0	0.741947024632614\\
0.02	0.738342555739801\\
0.04	0.734959582497297\\
0.06	0.731785479844873\\
0.08	0.728808564322176\\
0.1	0.726018007892732\\
0.12	0.723403761062111\\
0.14	0.720956484141454\\
0.16	0.718667485667016\\
0.18	0.716528667121406\\
0.2	0.714532473216768\\
0.22	0.712671847097804\\
0.24	0.710940189905835\\
0.26	0.709331324216516\\
0.28	0.707839460925061\\
0.3	0.706459169205613\\
0.32	0.705185349216862\\
0.34	0.704013207265438\\
0.36	0.702938233172667\\
0.38	0.701956179619988\\
0.4	0.701063043274102\\
0.42	0.700255047515487\\
0.44	0.699528626613581\\
0.46	0.698880411209231\\
0.48	0.698307214980102\\
0.5	0.697806022378123\\
0.52	0.697373977339726\\
0.54	0.697008372880025\\
0.56	0.696706641491244\\
0.58	0.696466346273775\\
0.6	0.696285172735487\\
0.62	0.696160921201256\\
0.64	0.696091499780403\\
0.66	0.696074917844765\\
0.68	0.696109279974672\\
0.7	0.69619278033412\\
0.72	0.696323697440055\\
0.74	0.696500389293938\\
0.76	0.696721288846625\\
0.78	0.696984899770271\\
0.8	0.697289792513246\\
0.82	0.697634600616237\\
0.84	0.698018017269558\\
0.86	0.698438792093454\\
0.88	0.698895728124724\\
0.9	0.699387678994401\\
0.92	0.699913546282509\\
0.94	0.70047227703706\\
0.96	0.701062861445523\\
0.98	0.701684330647945\\
1	0.702335754681768\\
};
\addplot [color=mycolor21, line width=1.0pt, forget plot]
  table[row sep=crcr]{%
0	0.75\\
0.06	0.75\\
0.08	0.749052359031822\\
0.1	0.746347573334519\\
0.12	0.743779246213124\\
0.14	0.741341190125955\\
0.16	0.73902758604968\\
0.18	0.736832956446488\\
0.2	0.734752140576477\\
0.22	0.732780271921248\\
0.24	0.73091275751104\\
0.26	0.729145258970758\\
0.28	0.727473675120429\\
0.3	0.725894125983351\\
0.32	0.724402938070862\\
0.34	0.722996630826365\\
0.36	0.721671904123459\\
0.38	0.720425626723769\\
0.4	0.719254825609609\\
0.42	0.718156676115082\\
0.44	0.717128492786768\\
0.46	0.716167720911814\\
0.48	0.715271928657277\\
0.5	0.714438799769857\\
0.52	0.713666126789972\\
0.54	0.712951804738386\\
0.56	0.712293825237426\\
0.58	0.711690271032302\\
0.6	0.711139310881104\\
0.62	0.710639194784858\\
0.64	0.710188249531523\\
0.66	0.709784874530095\\
0.68	0.709427537913006\\
0.7	0.709114772886888\\
0.72	0.708845174313416\\
0.74	0.708617395503496\\
0.76	0.708430145209421\\
0.78	0.70828218480087\\
0.8	0.708172325611796\\
0.82	0.708099426446226\\
0.84	0.708062391232012\\
0.86	0.708060166812365\\
0.88	0.708091740865846\\
0.9	0.708156139946165\\
0.92	0.708252427633825\\
0.94	0.708379702792228\\
0.96	0.708537097921422\\
0.98	0.708723777603174\\
1	0.7089389370315\\
};
\addplot [color=mycolor22, line width=1.0pt, forget plot]
  table[row sep=crcr]{%
0	0.75\\
0.26	0.75\\
0.28	0.749269261029426\\
0.3	0.747604178892658\\
0.32	0.746006245830294\\
0.34	0.744473251422788\\
0.36	0.743003081214675\\
0.38	0.741593711562704\\
0.4	0.740243204812345\\
0.42	0.738949704778501\\
0.44	0.737711432508283\\
0.46	0.736526682305508\\
0.48	0.735393817998269\\
0.5	0.734311269432398\\
0.52	0.733277529175039\\
0.54	0.732291149413786\\
0.56	0.731350739037983\\
0.58	0.730454960889822\\
0.6	0.729602529173829\\
0.62	0.728792207014182\\
0.64	0.728022804150122\\
0.66	0.727293174760433\\
0.68	0.726602215408626\\
0.7	0.725948863101099\\
0.72	0.725332093451079\\
0.74	0.724750918941695\\
0.76	0.724204387281984\\
0.78	0.723691579850087\\
0.8	0.723211610218288\\
0.82	0.722763622754907\\
0.84	0.722346791298435\\
0.86	0.721960317899579\\
0.88	0.721603431627199\\
0.9	0.721275387434377\\
0.92	0.720975465081127\\
0.94	0.720702968110456\\
0.96	0.720457222874736\\
0.98	0.720237577609505\\
1	0.720043401552051\\
};
\addplot [color=mycolor23, line width=1.0pt, forget plot]
  table[row sep=crcr]{%
0	0.75\\
0.66	0.75\\
0.68	0.749190115897933\\
0.7	0.748306586263624\\
0.72	0.747450175767248\\
0.74	0.746620295240122\\
0.76	0.745816372457267\\
0.78	0.745037851532665\\
0.8	0.74428419234022\\
0.82	0.743554869959183\\
0.84	0.742849374142821\\
0.86	0.742167208809219\\
0.88	0.74150789155313\\
0.9	0.740870953177888\\
0.92	0.740255937246399\\
0.94	0.739662399650337\\
0.96	0.739089908196669\\
0.98	0.738538042210717\\
1	0.738006392154975\\
};
\addplot [color=mycolor24, line width=1.0pt, forget plot]
  table[row sep=crcr]{%
0	0.75\\
1	0.75\\
};
\end{axis}
\end{tikzpicture}%

%% file: tikz/restriction_diamond.tex
%
\definecolor{mycolor1}{rgb}{0.00000,0.44700,0.74100}%
\definecolor{mycolor2}{rgb}{0.85000,0.32500,0.09800}%
\begin{tikzpicture}

\begin{axis}[%
width=0.742\textwidth,
height=0.585\textwidth,
at={(0\textwidth,0\textwidth)},
scale only axis,
xmin=-8,
xmax=8,
xlabel style={font=\color{white!15!black}, font=\small},
xlabel={\(\lambda\)},
xtick={-8,-6,-4,-2,0,2,4,6,8},
ymin=-0.012,
ymax=0.006,
ylabel style={font=\color{white!15!black}, font=\small},
ylabel={\(H_p(\lambda, \mu_{\eta}(\lambda))\)},
ticklabel style = {font=\scriptsize},
axis background/.style={fill=white},
xmajorgrids,
ymajorgrids,
legend style={at={(0.97,0.03)}, anchor=south east, legend cell align=left, align=left, draw=white!15!black, font=\small}
]
\addplot [very thick, color=blue!70!black, forget plot]
  table[row sep=crcr]{%
-8	-0.0104587382089354\\
-7.6734693877551	-0.00982516548093937\\
-7.34693877551021	-0.00921326167244629\\
-7.02040816326531	-0.00862302983688323\\
-6.69387755102041	-0.00805447339710819\\
-6.36734693877551	-0.00750759619750285\\
-6.04081632653061	-0.00698240256354232\\
-5.71428571428571	-0.00647889736951157\\
-5.38775510204082	-0.0059970861147427\\
-5.06122448979592	-0.00553697500794215\\
-4.73469387755102	-0.00509857105814859\\
-4.40816326530612	-0.00468188216819598\\
-4.08163265306122	-0.00428691722242647\\
-3.75510204081633	-0.00391368615261811\\
-3.42857142857143	-0.00356219995207763\\
-3.25651016550287	-0.00338582536632259\\
-3.10204081632653	-0.00323247058262588\\
-2.77551020408163	-0.00292451067287658\\
-2.44897959183673	-0.00263833281974613\\
-2.12244897959184	-0.00237394814031688\\
-1.79591836734694	-0.00213136339819364\\
-1.46938775510204	-0.00191057537460537\\
-1.16912122638553	-0.00172677920813236\\
-1.14285714285714	-0.0017115597781534\\
-0.816326530612246	-0.00153477702311022\\
-0.489795918367347	-0.0013798500810065\\
-0.163265306122449	-0.00124319392308969\\
-0.148174094712861	-0.00123790132984247\\
-0.130415842526355	-0.00123172023354989\\
-0.116490303593758	-0.00122690912930162\\
-0.105277515339946	-0.00122305866403583\\
-0.0960551395330221	-0.00121990773225233\\
-0.0883363854619134	-0.00121728195977155\\
-0.081781234893005	-0.0012150604635881\\
-0.0761450026586322	-0.00121315680692646\\
-0.0712471137767983	-0.00121150754850241\\
-0.0669513972972172	-0.00121006507389776\\
-0.0631532278961445	-0.0012087929507878\\
-0.059770891064392	-0.00120766282893425\\
-0.0567396373549111	-0.00120665231715655\\
-0.0540074968970846	-0.00120574349601554\\
-0.0515322747761268	-0.00120492185475296\\
-0.0492793560454194	-0.00120417551763702\\
-0.047220076784608	-0.00120349467191422\\
-0.045330497920175	-0.00120287113852413\\
-0.0435904702246505	-0.00120229804573535\\
-0.0419829128957204	-0.00120176957800311\\
-0.0404932508839493	-0.00120128078058102\\
-0.0391089716561939	-0.0012008274059454\\
0.0476969255283841	-0.00116800648508431\\
0.05006443548992	-0.00116731573735862\\
0.0526718681682996	-0.00116655424615431\\
0.0555576450862335	-0.00116571055720782\\
0.0587688542130908	-0.00116477060276843\\
0.0623638395208886	-0.00116371691124151\\
0.0664157780748257	-0.00116252751158896\\
0.0710177138449282	-0.00116117438557772\\
0.0762897907354006	-0.00115962123399083\\
0.0823898938136191	-0.00115782017409871\\
0.0895297320176238	-0.00115570672002363\\
0.0979999110903319	-0.0011531919062282\\
0.108210462229899	-0.00115014945912773\\
0.120759214445032	-0.00114639396174709\\
0.136553213553876	-0.00114164166108921\\
0.157038334284306	-0.00113543537795824\\
0.163265306122449	-0.0011335393234555\\
0.489795918367346	-0.00102888886865019\\
0.571686915355847	-0.00100092077039893\\
0.816326530612244	-0.000913215912753174\\
1.14285714285714	-0.000786443814295268\\
1.16560000993187	-0.000777202745573469\\
1.46938775510204	-0.000648622201577798\\
1.79591836734694	-0.000499939998245935\\
1.84186442182413	-0.000478172703857638\\
2.12244897959184	-0.000340743671868537\\
2.44897959183673	-0.000171563355159776\\
2.65239276919637	-6.14790787814012e-05\\
2.77551020408163	6.84935548278531e-06\\
3.10204081632653	0.00019346613911253\\
3.42857142857143	0.000386899996255689\\
3.75510204081632	0.000585491169132979\\
3.76882424856355	0.000594016441333167\\
4.08163265306122	0.000796316550611298\\
4.40816326530612	0.00102374441164343\\
4.73469387755102	0.00126773810540293\\
5.06122448979592	0.00152828952432138\\
5.2020603076096	0.00164576652164072\\
5.38775510204082	0.00180541134889545\\
5.71428571428571	0.00209913145523544\\
6.04081632653061	0.00240948915614148\\
6.36734693877551	0.00273682668347419\\
6.58353237862553	0.00297289005030699\\
6.69387755102041	0.00309848766506349\\
7.02040816326531	0.00348817928890455\\
7.3469387755102	0.00390446888341245\\
7.61633237779898	0.00426699297867827\\
7.6734693877551	0.00434623589072913\\
8	0.00481254774456374\\
};

\addplot[only marks, mark=*, mark options={}, mark size=2.5000pt, color=red!60!black, fill=red!60!black] table[row sep=crcr]{%
x	y\\
2.76325917006675	0\\
};
\addlegendentry{root \( \lambda^{\ast} \)}

\end{axis}
\end{tikzpicture}%

%% file: tex/conclusion.tex
\section{Conclusion, Discussion and Future Work}
In this paper we have analyzed the steady state problem of hydrogen blended natural gas transport in pipeline networks. We have clarified the differences to pure natural gas transport and discussed the existence and uniqueness of solutions on networks. For tree-shaped networks, we can compute unique flows by the mass conservation a priori, which guarantees the existence of a (unique) solution. 
For networks including a cycle, we have applied a cutting edge approach in order to show the existence of a solution. 
Cutting an edge in the cycle results in a tree-shaped graph with two new nodes and edges, and two additional coupling equations. We have shown that admissible solutions on the cut network exists which also meets the two additional coupling equations. This implies the existence of a solution for the graph with a cycle. 

The uniqueness of a solution for tree-structured networks follows by the fact, that the flow direction is a priori given and thus, it is also a priori clear which node determines the composition on an edge. 
For networks including cycles, the composition of the gas on the edges in the cycle depends on the flow direction, leading to discontinuities in the  mixing condition as discussed in \cref{sec: existence single cycle}. A natural extension of this work includes analyzing the discontinuities and thus also the uniqueness of a solution on networks with a cycle. 


Additionally, considering optimization problems that involve hydrogen blends offer a valuable direction for future research. Such problems could explore optimal blending ratios, operational settings, or cost efficiencies while ensuring safety and regulatory standards. 
Even though the simulation and numerical optimization results presented in \cite{Zlotnik2022,Pfetsch2024,Zlotnik2024,Zlotnik2024b,Zlotnik2024c,Zlotnik2023} are promising, the research on the theoretical analysis of optimization problems that involve hydrogen blended natural gas transport is still very limited. 


\section{Acknowledgment}
Simone G\"ottlich was supported by the German Research Foundation (DFG) under grant GO 1920/11-1 and 12-1.